\documentclass[11pt]{article}
\usepackage[utf8]{inputenc}
\usepackage[T1]{fontenc}
\usepackage[frenchb,english]{babel} 
\usepackage{textcomp}
\usepackage{amsmath,amssymb}
\usepackage{amsthm}
                  
\usepackage{lmodern}

\usepackage[a4paper,margin=2cm]{geometry}

\usepackage{graphicx}             
\usepackage{xcolor}               
\usepackage{microtype,stmaryrd}          

\usepackage{hyperref}
\hypersetup{pdfstartview=XYZ}

\usepackage{bbm}
\usepackage{mathrsfs}
\usepackage{subfig}

\def\build#1_#2^#3{\mathrel{
\mathop{\kern 0pt#1}\limits_{#2}^{#3}}}
\def\llbracket{[\hspace{-.10em} [ }
\def\rrbracket{ ] \hspace{-.10em}]}

\def\midd{\,|\,}

\newtheorem{theorem}{Theorem}
\newtheorem{proposition}[theorem]{Proposition}
\newtheorem{definition}[theorem]{Definition}
\newtheorem{lemma}[theorem]{Lemma}

\newtheorem{corollary}[theorem]{Corollary}

\def\w{\mathrm{w}}
\def\t{\mathcal{T}}

\def\W{\mathcal{W}}
\def\S{\mathcal{S}}

\def\N{\mathbb{N}}
\def\M{\mathbb{M}}

\def\D{\mathbb{D}}
\def\P{\mathbb{P}}
\def\E{\mathbb{E}}

\def\R{\mathbb{R}}
\def\z{\mathcal{Z}}

\def\n{\mathcal{N}}
\def\cc{\mathcal{C}}
\def\ee{\mathcal{E}}
\def\ve{{\varepsilon}}
\def\la{\longrightarrow}
\def\da{\downarrow}
\def\ov{\overline}

\def\dd{\mathrm{d}}
\def\wh{\widehat}
\def\wt{\widetilde}
\def\tr{\mathrm{tr}}

\def\nn{\mathcal{N}}
\def\pp{\mathcal{P}}

\def\ii{\mathcal{I}}
\def\jj{\mathcal{J}}

\def\bm{\mathbf{m}}

\def\EE{\mathbf{E}}
\def\SS{\mathbf{S}}
\def\HH{\mathbf{H}}
\def\RR{\mathbf{R}}

\def\rems{\noindent{\bf Remarks. }}

\author{Jean-Fran\c cois Le Gall}
\title{Geodesic stars in random geometry\footnote{Supported by the ERC Advanced Grant 740943 {\sc GeoBrown}}}
\date{\small Universit\'e Paris-Saclay}

\begin{document}

\maketitle

\begin{abstract}
A point of a metric space is called a geodesic star with $m$ arms if it
is the endpoint of $m$ disjoint geodesics. For every $m\in\{1,2,3,4\}$, we prove that the set 
of all geodesic stars with $m$ arms in the Brownian sphere has dimension $5-m$.
This complements recent results of Miller and Qian, who proved that this dimension
is smaller than or equal to $5-m$.
\end{abstract}


\section{Introduction}

This work is concerned with the continuous models of random geometry that have been
studied extensively in the recent years. In particular, we consider the Brownian sphere or
Brownian map, which is the scaling limit in the Gromov-Hausdorff sense of
triangulations or quadrangulations of the sphere with $n$ faces chosen uniformly
at random, and of much more general random planar maps (see in particular \cite{Abr,AA,BJM,Uniqueness,Mar,Mie-Acta}). We are primarily interested in the
study of geodesics in the Brownian sphere, but our main result remains valid in the
related models called the Brownian plane \cite{Plane,CLG} and the Brownian disk \cite{BM,Disks}.

Recall that a geodesic in a metric space $(E,d)$ is a continuous path $(\gamma(t))_{t\in[0,\delta]}$,
where $\delta>0$, such that $d(\gamma(s),\gamma(t))=|s-t|$ for every $s,t\in[0,\delta]$. For every $t\in(0,\delta)$, we say that
$\gamma(t)$ is an interior point of the geodesic (whereas $\gamma(0)$ and $\gamma(\delta)$ are its endpoints).
If $m\geq 1$ is an integer, we then say that a point $x$ is a geodesic star with $m$ arms (in short, an $m$-geodesic star) if there exist $\delta>0$
and $m$ geodesics $(\gamma_1(t))_{t\in[0,\delta]},\ldots,(\gamma_m(t))_{t\in[0,\delta]}$
such that $\gamma_1(0)=\gamma_2(0)=\cdots=\gamma_m(0)=x$ and the sets $\{\gamma_j(t):t\in(0,\delta]\}$, for $j\in\{1,\ldots,m\}$, are disjoint.
 If $(E,d)$
is a geodesic space (having more than one point), any pair of distinct points is connected by a (possibly not unique) geodesic, and it is then immediate 
that every point is a $1$-geodesic star. Our main result is the following theorem.

\begin{theorem}
\label{main-th}
Let $(\bm_\infty,D)$ denote the Brownian sphere. For every integer $m\in\{1,2,3,4\}$, let $\mathfrak{S}_m$
be the set of all $m$-geodesic stars in $(\bm_\infty,D)$. Then the Hausdorff dimension of
$\mathfrak{S}_m$ is a.s.\ equal to $5-m$. 
\end{theorem}

The upper bound $\mathrm{dim}(\mathfrak{S}_m)\leq 5-m$ has been obtained by Miller and Qian in \cite[Theorem 1.4]{MQ}. So the contribution of the present work is to prove the corresponding lower bound. We note that $m$-geodesic stars in the Brownian sphere were first discussed by Miermont \cite[Definition 7]{Mie-Acta}, who conjectured that 
they exist for $m\leq 4$ but not for $m\geq 6$ (see the concluding remarks of \cite{Mie-Acta}). In fact the
non-existence of $m$-geodesic stars when $m\geq 6$ has been proved by Miller and Qian \cite[Theorem 1.4]{MQ}.

Let us briefly comment on Theorem \ref{main-th}. The Brownian sphere is a geodesic space, and
thus $\mathfrak{S}_1=\bm_\infty$, so that  in the case $m=1$ the result follows  from the known fact \cite{Invent} that
$\mathrm{dim}(\bm_\infty)=4$. Next we may observe that any interior point of a geodesic
is a $2$-geodesic star, and therefore $\mathfrak{S}_2$ contains the set of all interior points of all
geodesics. However, Miller and Qian \cite[Corollary 1.3]{MQ} proved that the Hausdorff dimension of the latter set is $1$ (it is obviously greater than
or equal to $1$), thus confirming a conjecture of Angel, Kolesnik and Miermont \cite{AKM}. Since $\mathrm{dim}(\mathfrak{S}_2)=3$, this implies,
at least informally, that typical $2$-geodesic stars are not interior points of geodesics.
Let us then consider $m=3$. It is relatively easy to construct $3$-geodesic stars in the Brownian sphere. Indeed, write $x_*$ for the (first) distinguished point of $\bm_\infty$
(see Section \ref{Brown-sphere} below), and suppose that $x$ and $y$ are two points chosen independently according to
the volume measure on $\bm_\infty$. From the results of \cite{Acta}, the (a.s.\ unique) geodesics from $x$ to $x_*$ and from $y$ to $x_*$
coalesce  before hitting $x_*$, and the point at which they coalesce is a $3$-geodesic star. Again one expects that such points
are not typical $3$-geodesic stars. Finally, to the best of our knowledge, the existence of $4$-geodesic stars had not been established before.
We state an open question in the case $m=5$.

\medskip
\noindent{\bf Open problem.} {\it Prove or disprove the existence of $5$-geodesic stars in the Brownian map}. 

\medskip
\noindent If $5$-geodesic stars do exist, \cite[Theorem 1.4]{MQ} implies that $\mathrm{dim}(\mathfrak{S}_5)=0$. 

\medskip
Let us discuss the earlier work about geodesics in the Brownian sphere. The paper \cite{Acta} provides a 
complete description of geodesics ending at the distinguished point $x_*$ , by showing that all such geodesics must
be  ``simple geodesics'' (see Section \ref{hulls-sphere} for the definition of a simple geodesic 
in a slightly different model). It follows from this description that any two geodesics ending at $x_*$
must coalesce before hitting $x_*$. This is the so-called confluence of geodesics phenomenon, which implies that $x_*$ is (a.s.) not a $2$-geodesic star.
These results still hold if $x_*$ is replaced by a point chosen according to the volume measure
of $\bm_\infty$, by the symmetry properties of the Brownian sphere (see Section \ref{symm} below), and they
imply the uniqueness of the geodesic between two points chosen independently according to the volume measure.
A different approach to the latter property was given by Miermont \cite{Mie-AENS}, in the more general setting
of scaling limits of random planar maps in arbitrary genus.

Miermont's approach \cite{Mie-Acta} to the uniqueness of the Brownian sphere as the scaling limit of random quadrangulations 
makes heavy use of the notion of geodesic stars. We also note  that both \cite{Mie-Acta} and the alternative approach
to the uniqueness of the Brownian sphere developed in \cite{Uniqueness} strongly rely on the
characterization of geodesics to $x_*$. 

The paper \cite{AKM} by Angel, Kolesnik and Miermont goes further in the study of geodesics in the Brownian sphere.
In particular, it is proved in \cite{AKM} that $2$-geodesic stars form a set of first Baire category. Moreover, \cite{AKM}
contains a thorough discussion of the so-called geodesic networks: for each pair $(x,y)$ of distinct points in the
Brownian map, the geodesic network between $x$ and $y$ is the union of all geodesics from $x$ to $y$. 
If $y=x_*$ (more generally, if $y$ is a ``typical'' point) the results of \cite{Acta} show that the geodesic network consists of the union of
at most $3$ geodesics, but things may be much more complicated if $x$ and $y$ are both exceptional points. The paper \cite{AKM} studies the
possible ``normal'' geodesic networks (normality implies that there 
is a common point other than $x$ and $y$  to all geodesics between $x$ and $y$), and obtains in particular that there is a dense
set of pairs $(x,y)\in\bm_\infty\times\bm_\infty$ such that the geodesic network between $x$ and $y$ consists of $9$ distinct geodesics
--- in that case, both $x$ and $y$ must be $3$-geodesic stars. 
Even deeper results (without the normality assumption) are derived in the recent paper \cite{MQ} of Miller and Qian, which shows 
that $9$ is indeed the maximal number of geodesics between two points of $\bm_\infty$, and moreover 
computes the Hausdorff dimension of the set of pairs $(x,y)$ such that there are exactly $j$ geodesics between 
$x$ and $y$ (see \cite[Theorem 1.6]{MQ}). As already mentioned, \cite{MQ} also gives the upper bound
$\mathrm{dim}(\mathfrak{S}_m)\leq 5-m$. Both \cite{AKM} and \cite{MQ} make heavy use of strong forms of
the confluence of geodesics phenomenon, see in particular \cite[Proposition 12]{AKM} and \cite[Theorem 1.1]{MQ}.
It is worth pointing that certain analogs of the results of \cite{AKM} and \cite{MQ} have been derived
in the related setting of Liouville quantum gravity surfaces in 
the very recent papers \cite{Gwy} and \cite{GPS} (the confluence of geodesics phenomenon in that setting \cite{GM1}
played a major role in the proof of the uniqueness of the Liouville quantum gravity metric, see \cite{GM2}).

\smallskip
Let us outline the main steps of the proof of the lower bound $\mathrm{dim}(\mathfrak{S}_m)\geq 5-m$ when $m\geq 2$.
The key ideas are very similar to those that have been used in the study of exceptional points of 
Brownian motion. It is convenient to deal with the so-called free Brownian sphere, which 
means that $\bm_\infty$ is defined under the infinite Brownian snake excursion measure $\N_0$ (see Section \ref{Brown-sphere} below). For every $\ve\in(0,1)$ we introduce a set $\mathfrak{S}^\ve_m$ of ``$\ve$-approximate'' $m$-geodesic stars.
A point $x$ of $\bm_\infty$ belongs to $\mathfrak{S}^\ve_m$ if $1<D(x_*,x)<2$ and if there exist $m$ geodesics to $x$
that start at distance $1$ from $x$ and are disjoint up to the time when they arrive at distance $\ve$ from $x$. 
More precisely, we require for technical reasons that these geodesics start from the boundary of the hull of radius $1$ centered at $x$ relative to $x_*$
(roughly speaking, this hull is obtained by filling in the holes of
the ball of radius $1$ centered at $x$, except for the one containing $x_*$, see Section \ref{Brown-sphere}). Write $\mathrm{Vol}(\cdot)$ for the volume measure on
$\bm_\infty$. Using the symmetry properties of the Brownian sphere,
it is not hard to verify that 
\begin{equation}
\label{intro1}
\N_0\Big(\mathrm{Vol}(\mathfrak{S}^\ve_m)\Big)\geq c_m\,\ve^{m-1},
\end{equation}
with a positive constant $c_m$ independent of $\ve$. Then, if $\delta\in(0,1)$, we rely on a two-point estimate to get
the bound
\begin{equation}
\label{intro2}
\N_0\Bigg(\int\!\!\int\mathbf{1}_{\mathfrak{S}^\ve_m\times \mathfrak{S}^\ve_m}(x,y)\,D(x,y)^{-(5-m-\delta)}\,\mathrm{Vol}(\dd x)\mathrm{Vol}(\dd y)\Bigg)\leq c_{\delta,m}\,\ve^{2(m-1)},
\end{equation}
with a constant $c_{\delta,m}$ independent of $\ve$. From \eqref{intro1} and \eqref{intro2}, standard arguments
show that, at least on a set of positive $\N_0$-measure, the volume measure restricted to $\mathfrak{S}^\ve_m$ and
scaled by the factor $\ve^{-(m-1)}$ converges when 
$\ve$ tends to $0$, along a suitable subsequence, to a limiting random measure $\mu$ satisfying
$$\int\!\!\int D(x,y)^{-(5-m-\delta)}\,\mu(\dd x)\mu(\dd y)<\infty.$$
If we know that $\mu$ is supported on $\mathfrak{S}_m$, the classical Frostman lemma 
gives $\mathrm{dim}(\mathfrak{S}_m)\geq 5-m-\delta$. However, it is not obvious that $\mu$ is supported on $\mathfrak{S}_m$,
because, even if a sequence $(x_n)_{n\in\N}$ of $\ve_n$-approximate $m$-geodesic stars (with $\ve_n\to 0$)
converges, it does not necessarily follow that the limit belongs to $\mathfrak{S}_m$. To overcome this difficulty, we need
to modify the definition of $\mathfrak{S}^\ve_m$ by imposing that the geodesics to $x$ in this definition are not only disjoint
but sufficiently far apart from each other. Another delicate point is to prove that the desired property holds 
$\N_0$-a.e.\ and not only on a set of positive $\N_0$-measure. As usual, we rely on a kind
of zero-one law, which requires considering first the (scale invariant) Brownian plane and then
using a strong coupling between the Brownian plane and the Brownian sphere. 

\smallskip
The paper is organized as follows. Section \ref{sec-preli} is devoted to a number of 
preliminaries, including the Brownian snake construction of the Brownian sphere
as a measure metric space with two distinguished points denoted by $x_*$ and $x_0$,
and a discussion of the symmetry properties of the Brownian sphere, which roughly speaking say that
$x_*$ and $x_0$ play the same role as two points chosen independently
according to the (normalized) volume measure. Section \ref{sec-hull} starts with the construction
of the random metric space corresponding to the hull of radius $r>0$ centered at $x_*$ relative to $x_0$,
under $\N_0(\cdot\midd D(x_*,x_0)>r)$. This construction yields an explicit calculation 
of the probability that there are $m$ geodesics from the boundary of the hull to 
$x_*$ that stay disjoint until they hit the ball of radius $\ve$ centered at $x_*$. 
Then Theorem \ref{deco-hull}, which is a result of independent interest, shows that the
hull of radius $r>0$ centered at $x_*$ and relative to $x_0$ is independent of
its complement conditionally on its boundary size, and the complement itself
is a Brownian disk --- this is in fact an analog of a result proved in \cite{Spine} for the Brownian
plane. One then derives a two point-version saying that, under
the conditional probability measure $\N_0(\cdot\midd D(x_*,x_0)>2r)$, the hull
of radius $r$ centered at $x_*$ (relative to $x_0$) and the hull
of radius $r$ centered at $x_0$ (relative to $x_*$) are independent conditionally on their boundary sizes
(Corollary \ref{two-hulls}). Section \ref{sec-first-mom} is devoted to the proof
of the version of \eqref{intro1} where the definition of $\mathfrak{S}_m^\ve$ is modified as
explained above to ensure that geodesics stay ``sufficiently far apart'' from each other. 
An important ingredient here is the notion of a slice, which roughly speaking separates
two successive disjoint geodesics from the hull boundary to the ball of radius $\ve$ (slices 
also played a key role in the characterization of the distribution of 
Brownian disks in \cite{BM}). Section \ref{sec:key-est},
which is the most technical part of the paper, uses the results of Section \ref{sec-hull} to derive the key estimate (Lemma \ref{key-lem}) that
eventually leads to the bound \eqref{intro2}.  Section \ref{sec-proof-main} then gives the proof of Theorem \ref{main-th}
along the lines of the preceding discussion. The Appendix contains the proof
of a couple of technical lemmas, including the strong coupling between the Brownian plane and the Brownian sphere
that is used to justify the zero-one law argument. 

We finally mention that Jason Miller and Wei Qian \cite{MQ1} have independently developed 
a different approach to Theorem \ref{main-th}. We believe that our general strategy and the
intermediate steps of our proof are of independent interest and should prove useful to investigate other sets 
of exceptional points in the Brownian sphere. We hope to pursue this matter in
the future.

\smallskip
\noindent{\bf Acknowledgement.} I thank Jason Miller and Wei Qian for 
keeping me informed of their ongoing work.

\section{Preliminaries}
\label{sec-preli}

\subsection{Measure metric spaces}
\label{sec-GHP}

A (compact) measure metric space is a compact metric space $(X,d)$ equipped with a finite Borel
 measure $\mu$ which is often called the volume measure. We write $\M$ for the set of all measure metric spaces, where two such spaces
$(X,d,\mu)$ and $(X',d',\mu')$ are identified if there exists an isometry $\phi$ from $X$ onto 
$X'$ such that $\phi_*\mu=\mu'$. 

For our purposes, it will be important to consider measure metric spaces given together
with two distinguished closed subsets that we call the boundaries for reasons that
will become clear later. We say that $(X,d,\mu,F_1,F_2)$ is a two-boundary measure metric space
if $(X,d,\mu)$ is a measure metric space and if $F_1$ and $F_2$
are two closed subsets of $X$ (the order between $F_1$ and $F_2$ is important). We write $\M^{bb}$
for the set of all two-boundary measure metric spaces modulo isometries (of course we
now consider only isometries that preserve both the volume measures and the ``boundaries'').

The Gromov-Hausdorff-Prokhorov distance on $\M^{bb}$ is then defined by
\begin{align*}
&d_{\mathrm{GHP}}((X,d,\mu,F_1,F_2),(X',d',\mu',F'_1,F'_2))\\
&= \inf_{\phi:X\to E,\phi':X'\to E} \Big\{ d^E_\mathrm{H}(\phi(X),\phi'(X)) \vee d^E_\mathrm{H}(\phi(F_1),\phi'(F'_1))\vee 
d^E_\mathrm{H}(\phi(F_2),\phi'(F'_2))\vee d^E_\mathrm{P}(\phi_*\mu,\phi'_*\mu')\Big\},
\end{align*}
where the infimum is over all isometric embeddings $\phi$ and $\phi'$ of $X$ and $X'$ into a compact
metric space $(E,d^E)$, $d^E_\mathrm{H}$ stands for the Hausdorff distance between compact subsets of $E$ and
$d^E_\mathrm{P}$ is the Prokhorov distance on the space of finite measures on $E$. 
Then, by an easy generalization of \cite[Theorem 2.5]{AD}, one verifies that $d_{\mathrm{GHP}}$
is a distance on $\M^{bb}$, and $(\M^{bb},d_{\mathrm{GHP}})$ is a Polish space.

We also  let  $\M^{\bullet b}$, resp. $\M^{\bullet\bullet}$, denote the closed subset of $\M^{bb}$ that consists of all 
(isometry classes of) two-boundary measure metric spaces $(X,d,\mu,F_1,F_2)$ such that 
$F_1$ is a singleton, resp. both $F_1$ and $F_2$ are singletons. Note that
$\M^{\bullet\bullet}$ is just the space of two-pointed measure metric spaces as considered in 
\cite[Section 2.1]{Disks}. 

\smallskip
\noindent{\bf Remark.} Gwynne and Miller \cite{GM0} consider the closely
related notion of a curve-decorated measure metric space.

\subsection{Snake trajectories}
\label{sna-tra}

We will use the formalism of snake trajectories as developed in \cite{ALG}. First recall that
a finite path $\w$ is a continuous mapping $\w:[0,\zeta]\la\R$, where the
number $\zeta=\zeta_{(\w)}\geq 0$ is called the lifetime of $\w$. We let 
$\W$ denote the space of all finite paths, which is a Polish space when equipped with the
distance
$$d_\W(\w,\w')=|\zeta_{(\w)}-\zeta_{(\w')}|+\sup_{t\geq 0}|\w(t\wedge
\zeta_{(\w)})-\w'(t\wedge\zeta_{(\w')})|.$$
The endpoint or tip of the path $\w$ is denoted by $\wh \w=\w(\zeta_{(\w)})$.
For $x\in\R$, we
set $\W_x=\{\w\in\W:\w(0)=x\}$. The trivial element of $\W_x$ 
with zero lifetime is identified with the point $x$ of $\R$.

\begin{definition}
\label{def:snakepaths}
Let $x\in\R$. 
A snake trajectory with initial point $x$ is a continuous mapping $s\mapsto \omega_s$
from $\R_+$ into $\W_x$ 
which satisfies the following two properties:
\begin{enumerate}
\item[\rm(i)] We have $\omega_0=x$ and the number $\sigma(\omega):=\sup\{s\geq 0: \omega_s\not =x\}$,
called the duration of the snake trajectory $\omega$,
is finite (by convention $\sigma(\omega)=0$ if $\omega_s=x$ for every $s\geq 0$). 
\item[\rm(ii)] {\rm (Snake property)} For every $0\leq s\leq s'$, we have
$\omega_s(t)=\omega_{s'}(t)$ for every $t\in[0,\displaystyle{\min_{s\leq r\leq s'}} \zeta_{(\omega_r)}]$.
\end{enumerate} 
\end{definition}

We will write $\S_x$ for the set of all snake trajectories with initial point $x$
and $\S=\bigcup_{x\in\R}\S_x$ for the set of all snake trajectories. If $\omega\in \S$, we often write $W_s(\omega)=\omega_s$ and $\zeta_s(\omega)=\zeta_{(\omega_s)}$
for every $s\geq 0$. The set $\S$ is a Polish space for the distance
$d_{\S}(\omega,\omega')= |\sigma(\omega)-\sigma(\omega')|+ \sup_{s\geq 0} \,d_\W(W_s(\omega),W_{s}(\omega'))$.
A snake trajectory $\omega$ is completely determined 
by the knowledge of the lifetime function $s\mapsto \zeta_s(\omega)$ and of the tip function $s\mapsto \wh W_s(\omega)$: See \cite[Proposition 8]{ALG}.

Let $\omega\in \S$ be a snake trajectory and $\sigma=\sigma(\omega)$. The lifetime function $s\mapsto \zeta_s(\omega)$ codes a
compact $\R$-tree, which will be denoted 
by $\t_{(\omega)}$ and called the {\it genealogical tree} of the snake trajectory. This $\R$-tree is the quotient space $\t_{(\omega)}:=[0,\sigma]/\!\sim$ 
of the interval $[0,\sigma]$
for the equivalence relation
$$s\sim s'\ \hbox{if and only if }\ \zeta_s(\omega)=\zeta_{s'}(\omega)= \min_{s\wedge s'\leq r\leq s\vee s'} \zeta_r(\omega),$$
and $\t_{(\omega)}$ is equipped with the distance induced by
$$d_{(\omega)}(s,s')= \zeta_s(\omega)+\zeta_{s'}(\omega)-2 \min_{s\wedge s'\leq r\leq s\vee s'} \zeta_r(\omega).$$
(notice that $d_{(\omega)}(s,s')=0$ if and only if $s\sim s'$).
We write $p_{(\omega)}:[0,\sigma]\la \t_{(\omega)}$
for the canonical projection. By convention, $\t_{(\omega)}$ is rooted at the point
$\rho_{(\omega)}:=p_{(\omega)}(0)$,
and the volume measure on $\t_{(\omega)}$ is defined as the pushforward of
Lebesgue measure on $[0,\sigma]$ under $p_{(\omega)}$. The mapping $s\mapsto p_{(\omega)}(s)$
may be interpreted as a (clockwise) {\it cyclic exploration} of $\t_{(\omega)}$. 

It will be useful to define also intervals on the tree $\t_{(\omega)}$. For $s,s'\in[0,\sigma]$, 
we use the convention $[s,s']=[s,\sigma]\cup [0,s']$ if $s>s'$
(and of course, $[s,s']$ is the usual interval if $s\leq s'$). If $a,b\in \t_{(\omega)}$ are distinct, then we can find
$s,s'\in[0,\sigma]$ in a unique way  so that $p_{(\omega)}(s)=a$ and $p_{(\omega)}(s')=b$ and 
the interval $[s,s']$ is as small as possible, and we define $[a,b]:=p_{(\omega)}([s,s'])$. 
Informally, $[a,b]$ is the set of all points that are visited when going from $a$ to $b$
in ``clockwise order'' around the tree. We also take $[a,a]=\{a\}$.

By property (ii) in the definition of  a snake trajectory, the condition $p_{(\omega)}(s)=p_{(\omega)}(s')$ implies that 
$W_s(\omega)=W_{s'}(\omega)$. So the mapping $s\mapsto W_s(\omega)$ could be viewed as defined on the quotient space $\t_{(\omega)}$.
For $a\in\t_{(\omega)}$, we set $\ell_a(\omega):=\wh W_s(\omega)$ for any  $s\in[0,\sigma]$ such that $a=p_{(\omega)}(s)$  (by the previous observation, this does not
depend on the choice of $s$). We then interpret $\ell_a(\omega)$ as a label assigned to the point $a$ of $\t_{(\omega)}$, and we observe that, if $p_{(\omega)}(s)=a$, the
path $[0,\zeta_s]\ni t\mapsto W_s(t)$ records the labels along the line segment from $\rho_{(\omega)}$ to $a$ in $\t_{(\omega)}$. 
We also note that the mapping $a\mapsto \ell_a(\omega)$ is continuous on $\t_{(\omega)}$.
We will use the notation 
$$W_*(\omega):=\min\{\ell_a(\omega):a\in\t_{(\omega)}\}=\min\{\wh W_s(\omega):0\leq s\leq \sigma\}.$$

Let us introduce a truncation operation on snake trajectories. Let $x,y\in\R$ with $y<x$. If $\w\in \mathcal{W}_x$, we set $\tau_y(\w):=\inf\{t\geq 0:\w(t)=y\}$,
with the usual convention $\inf\varnothing=\infty$. Then, if 
$\omega\in \S_x$, we set, for every $s\geq 0$,
$$\eta_s(\omega)=\inf\Big\{t\geq 0:\int_0^t \mathrm{d}u\,\mathbf{1}_{\{\zeta_{(\omega_u)}\leq\tau_y(\omega_u)\}}>s\Big\}.$$
Note that the condition $\zeta_{(\omega_u)}\leq\tau_y(\omega_u)$ holds if and only if $\tau_y(\omega_u)=\infty$ or $\tau_y(\omega_u)=\zeta_{(\omega_u)}$.
Then, setting $\omega'_s=\omega_{\eta_s(\omega)}$ for every $s\geq 0$ defines an element $\omega'$ of $\S_x$,
which will be denoted by  $\tr_y(\omega)$ and called the truncation of $\omega$ at $y$
(see \cite[Proposition 10]{ALG}). The effect of the time 
change $\eta_s(\omega)$ is to ``eliminate'' those paths $\omega_s$ that hit $y$ and then survive for a positive
amount of time. We can then also define the excursions of $\omega$ ``below'' level $y$. To this end, 
we let $(\alpha_j,\beta_j)$, $j\in J$, be the connected components of the open set
$$\{s\in[0,\sigma]:\tau_y(\omega_s)<\zeta_{(\omega_s)}\},$$
and notice that $\omega_{\alpha_j}=\omega_{\beta_j}$ for every $j\in J$.
For every $j\in J$ we define a snake trajectory $\omega^j\in\S_0$ by setting
$$\omega^j_{s}(t):=\omega_{(\alpha_j+s)\wedge\beta_j}(\zeta_{(\omega_{\alpha_j})}+t)-y\;,\hbox{ for }0\leq t\leq \zeta_{(\omega^j_s)}
:=\zeta_{(\omega_{(\alpha_j+s)\wedge\beta_j})}-\zeta_{(\omega_{\alpha_j})}\hbox{ and } s\geq 0.$$
We say that $\omega_j$, $j\in J$ are the excursions of $\omega$ below level $y$. 

We finally introduce the re-rooting operation on snake trajectories (see \cite[Section 2.2]{ALG}). Let $\omega\in \S_0$ and
$r\in[0,\sigma(\omega)]$. Then $\omega^{[r]}$ is the snake trajectory in $\S_0$ such that
$\sigma(\omega^{[r]})=\sigma(\omega)$ and for every $s\in [0,\sigma(\omega)]$,
\begin{align*}
\zeta_s(\omega^{[r]})&= d_{(\omega)}(r,r\oplus s),\\
\wh W_s(\omega^{[r]})&= \wh W_{r\oplus s}(\omega)-\wh W_r(\omega),
\end{align*}
where we use the notation $r\oplus s=r+s$ if $r+s\leq \sigma(\omega)$, and $r\oplus s=r+s-\sigma(\omega)$ otherwise. 
These prescriptions completely determine $\omega^{[r]}$.
The genealogical tree $\t_{(\omega^{[r]})}$ may be interpreted as the tree $\t_{(\omega)}$ re-rooted at the vertex $p_{(\omega)}(r)$,
and vertices of $\t_{(\omega^{[r]})}$
receive the same labels as in $\t_{(\omega)}$, shifted so that the label of the (new) root is still $0$.

\subsection{The Brownian snake excursion 
measure on snake trajectories}
\label{sna-mea}

Let $x\in\R$. The Brownian snake excursion 
measure $\N_x$ is the $\sigma$-finite measure on $\S_x$ that satisfies the following two properties: Under $\N_x$,
\begin{enumerate}
\item[(i)] the distribution of the lifetime function $(\zeta_s)_{s\geq 0}$ is the It\^o 
measure of positive excursions of linear Brownian motion, normalized so that, for every $\ve>0$,
$$\N_x\Big(\sup_{s\geq 0} \zeta_s >\ve\Big)=\frac{1}{2\ve};$$
\item[(ii)] conditionally on $(\zeta_s)_{s\geq 0}$, the tip function $(\wh W_s)_{s\geq 0}$ is
a Gaussian process with mean $x$ and covariance function 
$$K(s,s'):= \min_{s\wedge s'\leq r\leq s\vee s'} \zeta_r.$$
\end{enumerate}
Informally, the lifetime process $(\zeta_s)_{s\geq 0}$ evolves under $\N_x$ like a Brownian excursion,
and conditionally on $(\zeta_s)_{s\geq 0}$, each path $W_s$ is a linear Brownian path started from $x$ with lifetime $\zeta_s$, which
is ``erased'' from its tip when $\zeta_s$ decreases and is ``extended'' when $\zeta_s$ increases. We note that the
density of $\sigma$ under $\N_0$ is $(2\sqrt{2\pi s^3})^{-1}$. 

For every $y<x$, we have
\begin{equation}
\label{hittingpro}
\N_x(W_*\leq y)={\displaystyle \frac{3}{2(x-y)^2}}.
\end{equation}
See e.g. \cite[Section VI.1]{Zurich} for a proof. Additionally, one can prove that $\N_x(\dd \omega)$ a.e. there is
a unique $s_*\in(0,\sigma)$ such that $\wh W_{s_*}=W_*$ (see \cite[Proposition 2.5]{LGW})
and we set $a_*=p_{(\omega)}(s_*)$ so that $\ell_{a_*}=W_*$.  

The following scaling property is often useful. For $\lambda>0$, for every 
$\omega\in \S_x$, we define $\theta_\lambda(\omega)\in \S_{x\sqrt{\lambda}}$
by $\theta_\lambda(\omega)=\omega'$, with
$$\omega'_s(t):= \sqrt{\lambda}\,\omega_{s/\lambda^2}(t/\lambda)\;,\quad
\hbox{for } s\geq 0\hbox{ and }0\leq t\leq \zeta'_s:=\lambda\zeta_{s/\lambda^2}.$$
Then it is a simple exercise to verify that the pushforward of $\N_x$ under $\theta_\lambda$ is  $\lambda\, \N_{x\sqrt{\lambda}}$. 

For every $t>0$, we can make sense of the conditional probability measure $\N_0^{(t)}:=\N_0(\cdot\midd \sigma =t)$.
If $s\in[0,t]$, $\N_0^{(t)}$ is invariant under the re-rooting operation $\omega\mapsto \omega^{[s]}$ (see e.g. \cite[Theorem 2.3]{LGW}).

\medskip
\noindent{\bf Exit measures.} Let $x,y\in\R$, with $y<x$. Under the measure $\N_x$, one can make sense of the ``quantity'' of paths 
$W_s$ that hit level $y$. One shows \cite[Proposition 34]{Disks} that the limit
\begin{equation}
\label{formu-exit}
L^y_t:=\lim_{\ve \da 0} \frac{1}{\ve^2} \int_0^t \dd s\,\mathbf{1}_{\{\tau_y(W_s)=\infty,\, \wh W_s<y+\ve\}}
\end{equation}
exists uniformly for $t\geq 0$, $\N_x$ a.e., and defines a continuous nondecreasing function, which is 
obviously constant on $[\sigma,\infty)$. 
The process $(L^y_t)_{t\geq 0}$ is called the exit local time at level $y$, and the exit measure 
$\z_y$ is defined by $\z_y=L^y_\infty=L^y_\sigma$. Then, $\N_x$ a.e., the topological support of the measure 
$\dd L^y_t$ is exactly the set $\{s\in[0,\sigma]:\tau_y(W_s)=\zeta_s\}$, and, in particular, $\z_y>0$ if and only if one of the paths $W_s$ hits $y$. The definition of $\z_y$
is a special case of the theory of exit measures (see \cite[Chapter V]{Zurich} for this general theory). 
Notice that the quantities in the right-hand side of \eqref{formu-exit} are functions of $\tr_y(\omega)$. 

The special Markov property of the Brownian snake states that, under $\N_x(\dd\omega)$ and conditionally on 
the truncation $\tr_y(\omega)$, the excursions of $\omega$ below $y$ form a Poisson measure
with intensity $\z_y\N_0$ (see the Appendix of \cite{subor} for a more precise statement).

\subsection{Decomposing the Brownian snake at its minimum}
\label{spine-decomp}

Let $u>0$. We will use the description of the conditional measure $\N_0(\cdot\midd W_*=-u)$, which can be found
in \cite{Bessel}. Let $(\alpha_i,\beta_i)$, $i\in I$, be the connected components of 
$\{s\in[0,s_*]:\zeta_s>\min_{[s,s_*]}\zeta_r\}$, and similarly let $(\alpha_i,\beta_i)$, $i\in J$, be the connected components of 
$\{s\in[s_*,\sigma]:\zeta_s>\min_{[s_*,s]}\zeta_r\}$ (the indexing sets $I$ and $J$ are disjoint). Notice that $\zeta_{\alpha_i}=\zeta_{\beta_i}$ and
$\wh W_{\alpha_i}=\wh W_{\beta_i}=W_{s*}(\zeta_{\alpha_i})$ by the snake property. For every 
$i\in I\cup J$, we write $\omega^{(i)}$ for the unique snake trajectory such that 
$$\zeta_s(\omega^{(i)})=\zeta_{(\alpha_i+s)\wedge \beta_i}-\zeta_{\alpha_i}\;,\quad \wh W_s(\omega^{(i)})=\wh W_{(\alpha_i+s)\wedge \beta_i}.$$
Then, under $\N_0(\cdot\midd W_*=-u)$, the finite path $(u+W_{s_*}(\zeta_{s_*}-t))_{0\leq t\leq \zeta_{s_*}}$ is distributed as a nine-dimensional Bessel process 
started at $0$ and stopped at its last passage time at $u$, and, conditionally on $W_{s_*}$, the point measures
$$\sum_{i\in I} \delta_{(\zeta_{\alpha_i},\omega^{(i)})}(\dd t\dd \omega')\hbox{ and }\sum_{i\in J} \delta_{(\zeta_{\alpha_i},\omega^{(i)})}(\dd t\dd \omega')$$
are independent Poisson measures with intensity
$$2\,\mathbf{1}_{[0,\zeta_{s_*}]}(t)\,\mathbf{1}_{\{W_*(\omega')>-u\}}\,\dd t\,\N_{W_{s_*}(t)}(\dd\omega').$$
Informally, this corresponds to a spine decomposition of the labeled tree $\t_{(\omega)}$
under $\N_0(\cdot\midd W_*=-u)$: $W_{s_*}$ records the labels along a spine isometric to $[0,\zeta_{s_*}]$, and each 
(labeled) tree $\t_{(\omega^{(i)})}$ for $i\in I$, resp.~for $i\in J$,
is grafted to the left side of the spine, resp.~to the right side of the spine, at level $\zeta_{\alpha_i}$. 

Let $r>0$. The preceding decomposition can be used to make sense of the exit measure $\z_{W_*+r}$
under $\N_0(\cdot\midd W_*<-r)$. Notice that the exit measure $\z_{y}$  was defined in the previous section for
a {\it deterministic} level $y$, whereas here $W_*+r$ is random. Nonetheless, we may 
argue under $\N_0(\cdot\midd W_*=-u)$ for every fixed $u>r$, and then define
\begin{equation}
\label{def-exit2}
\z_{W_*+r}:=\sum_{\{i\in I\cup J:\,\zeta_{\alpha_i}>\tau_{-u+r}(W_{s_*})\}} \z_{-u+r}(\omega^{(i)})
\end{equation}
with the notation above. Moreover, the special Markov property implies that the distribution 
of $\tr_{-u+r}(\omega)$ under $\N_0(\cdot\midd W_*=-u)$ is absolutely continuous with
respect to its distribution 
under $\N_0(\cdot\midd W_*<-u+r)$, so that
we can apply  \eqref{formu-exit} (with $y=-u+r$ and $x=0$) under $\N_0(\cdot\midd W_*=-u)$, and get that
\begin{equation}
\label{formu-exit2}
\z_{W_*+r}=\lim_{\ve\to 0} \ve^{-2}\int_0^\sigma \dd s\,\mathbf{1}_{\{\tau_{W_*+r}(W_s)=\infty, \wh W_s<W_*+r+\ve\}}\;,
\end{equation}
$\N_0(\cdot\midd W_*=-u)$ a.e. and thus also $\N_0(\cdot\midd W_*<-r)$ a.e.

\subsection{The Brownian sphere}
\label{Brown-sphere}

Let us argue under the excursion measure $\N_0(\dd \omega)$,
and recall the notation $\ell_a=\ell_a(\omega)$ for the label
assigned to $a\in\t_{(\omega)}$, and the definition of intervals on $\t_{(\omega)}$. 
We define, for every $a,b\in\t_{(\omega)}$,
\begin{equation}
\label{Dzero}
D^\circ(a,b):=\ell_a + \ell_b -2\max\Big(\min_{c\in[a,b]} \ell_c, \min_{c\in[b,a]} \ell_c\Big).
\end{equation}
We record two easy but important properties of $D^\circ$. First, for every $a,b\in\t_\zeta$,
\begin{equation}
\label{low-bd-D}
D^\circ(a,b)\geq |\ell_a-\ell_b|.
\end{equation}
Then,  recalling that $a_*$ is the unique point such that $\ell_{a_*}=W_*$, we have for every $a\in\t_\zeta$,
\begin{equation}
\label{equal-D}
D^\circ(a_*,a)=\ell_a-\ell_{a_*}.
\end{equation}

We let $D(a,b)$ be the largest symmetric function of the pair $(a,b)$ that
is bounded above by $D^\circ(a,b)$ and satisfies the triangle inequality: For every
$a,b\in\t_\zeta$,
\begin{equation}
\label{formulaD}
D(a,b) = \inf\Big\{ \sum_{i=1}^k D^\circ(a_{i-1},a_i)\Big\},
\end{equation}
where the infimum is over all choices of the integer $k\geq 1$ and of the
elements $a_0,a_1,\ldots,a_k$ of $\t_{(\omega)}$ such that $a_0=a$
and $a_k=b$. We note that $D$ is a pseudo-metric on $\t_{(\omega)}$, and 
we let $\bm_\infty(\omega):=\t_{(\omega)}/\{D=0\}$ be the associated quotient space (that is, the quotient space of
$\t_{(\omega)}$ for the equivalence relation $a\approx b$ if and only if $D(a,b)=0$), which is
equipped with the distance induced by $D$, for which we keep the same notation $D$.
Then $(\bm_\infty,D)$ is a compact metric space and also a geodesic space. 
The
canonical projection from $\t_\zeta$ onto $\bm_\infty$ is denoted by $\Pi$, and
the 
volume measure $\mathrm{Vol}$ on $\bm_\infty$ is defined as the pushforward of the volume measure 
on $\t_{(\omega)}$ under $\Pi$. Note that the total mass of $\mathrm{Vol}$ is $\sigma$. 

We view $(\bm_\infty,D,\mathrm{Vol})$
as a random  two-pointed measure metric space, or equivalently as a random variable with values
in the space $\M^{\bullet\bullet}$ of Section \ref{sec-GHP}: the first distinguished point of $\bm_\infty$ is $x_*=\Pi(a_*)$ and the
second distinguished point is $x_0:=\Pi(\rho_{(\omega)})$. As a consequence of \eqref{equal-D}, we have
$$D(x_*,x_0)=-W_*.$$

The free Brownian sphere is the (two-pointed measure) metric space $\bm_\infty$
under the measure $\N_0$. We note that $\bm_\infty$ is a geodesic space.  It makes sense, and is also of interest, to consider
$\bm_\infty$ under conditional measures. The space $\bm_\infty$
under the probability measure $\N^{(1)}_0:=\N_0(\cdot\midd \sigma=1)$ is the standard
Brownian sphere (or Brownian map). For every $r>0$, we will also consider
$\bm_\infty$ under the conditional measure
$\N_0^{[r]}:=\N_0(\cdot\midd W_*<-r)$. 
This corresponds to conditioning the free Brownian sphere on the event that the distance 
between the two distinguished points is greater than $r$. 

The metric space $\bm_\infty$ is homeomorphic to the usual two-dimensional sphere,
$\N_0$ a.e. (or a.s. for any of the conditional measures introduced above). 

For $x\in\bm_\infty$
and $r>0$, let $B_r(x)$ denote the closed ball of radius $r$ centered at $x$ in $\bm_\infty$.
If $x$ and $y$ are distinct points of $\bm_\infty$ and $r\in (0,D(x,y))$, the hull $B^{\bullet(y)}_r(x)$ 
is the closed subset of $\bm_\infty$ such that $\bm_\infty\backslash B^{\bullet(y)}_r(x)$ is  the connected component of the complement
of $B_r(x)$ that contains $y$. We say that $B^{\bullet(y)}_r(x)$ is
the hull of radius $r$ centered at $x$ relative to $y$. Obviously $B_r(x)\subset B^{\bullet(y)}_r(x)$, and every 
point of the topological boundary $\partial B^{\bullet(y)}_r(x)$ is at distance $r$ from $x$. 

The following fact known as the {\it cactus bound} \cite[Proposition 3.1]{Acta} is useful to study hulls centered at $x_*$. Let $x=\Pi(a)$ and $y=\Pi(b)$ be two points of
$\bm_\infty$, and let $(\gamma(t))_{0\leq t\leq 1}$ be a continuous path in $\bm_\infty$ such that $\gamma(0)=x$
and $\gamma(1)=y$. Then,
\begin{equation}
\label{cactus-bd}
\min_{0\leq t\leq 1} D(x_*,\gamma(t))\leq \min_{c\in\llbracket a,b\rrbracket} \ell_c - W_*,
\end{equation}
where $\llbracket a,b\rrbracket$ denotes the line segment between $a$ and $b$ in $\t_{(\omega)}$.

\subsection{Symmetry properties}
\label{symm}

As mentioned above, $\bm_\infty$
is viewed as a two-pointed measure metric space, and in this section we 
will write $(\bm_\infty,x_*,x_0)$ to make it explicit that $x_*$ and $x_0$
are the distinguished points. Our goal here is to
observe that $x_*$ and $x_0$ can be replaced by points chosen uniformly
according to the volume measure without changing the distribution
of $(\bm_\infty,x_*,x_0)$. This will be very important for our
applications.

\begin{proposition}
\label{interchange}
Let $F$ be a nonnegative measurable function on the space
$\M^{\bullet\bullet}$. Then,
$$\N_0(F(\bm_\infty,x_*,x_0))=\N_0\Big(\int\!\!\int \frac{\mathrm{Vol}(\dd x)}{\sigma}\,\frac{\mathrm{Vol}(\dd y)}{\sigma}
\,F(\bm_\infty,x,y)\Big).$$
The same identity holds if $\N_0$ is replaced by $\N^{(s)}_0=\N_0(\cdot\midd \sigma=s)$, for any $s>0$. 
\end{proposition}

\proof It is enough to treat the case of $\N^{(s)}_0$. We fix $s>0$ and recall that, for every
$t\in[0,s]$, the measure $\N_0^{(s)}$ is invariant under the re-rooting operation $\omega\mapsto \omega^{[t]}$. 
On the other hand, it is easy to verify that the measure metric space $\bm_\infty$ is left unchanged if
$\omega$ is replaced by $\omega^{[t]}$, and that the first distinguished point also remains the same
(the minimal label is attained at the ``same'' point of $\t_{(\omega)}$ and $\t_{(\omega^{[t]})}$). However,
the second distinguished point $x_0$ is replaced by $\Pi(p_{(\omega)}(t))$. It follows from these
considerations and the definition of the volume measure that
$$\N_0^{(s)}(F(\bm_\infty,x_*,x_0))= \N_0^{(s)}\Big(\frac{1}{s}\int \mathrm{Vol}(\dd y)\,F(\bm_\infty,x_*,y\Big).$$
Then an application of \cite[Theorem 8.1]{Acta} shows that
the right-hand side is also equal to
$$\N_0^{(s)}\Big(\int\!\!\int \frac{\mathrm{Vol}(\dd x)}{s}\,\frac{\mathrm{Vol}(\dd y)}{s}
\,F(\bm_\infty,x,y)\Big).$$
To be precise, \cite[Theorem 8.1]{Acta} considers $\bm_\infty$ as a metric space, and so we 
need a slight extension of this result, when $\bm_\infty$ is viewed as a {\it measure} metric space.
This extension is obtained by the very same arguments as in \cite{Acta}, using the convergence of 
rescaled quadrangulations to the Brownian sphere in the Gromov-Hausdorff-Prokhorov sense,
as stated in \cite[Theorem 7]{Disks}. This completes the proof. \endproof

Let us mention some immediate consequences of Proposition \ref{interchange}. First, we have
\begin{equation}
\label{inter-2}
\N_0(F(\bm_\infty,x_*,x_0))=\N_0(F(\bm_\infty,x_0,x_*)).
\end{equation}
Since the conditioning defining $\N_0^{[r]}=\N_0(\cdot\midd D(x_0,x_*)>r)$ depends on 
$x_*$ and $x_0$ in a symmetric manner, we get for every $r>0$,
$$\N^{[r]}_0(F(\bm_\infty,x_*,x_0))=\N^{[r]}_0(F(\bm_\infty,x_0,x_*)).$$

The following consequence of Proposition \ref{interchange} will also be useful. If $F$ is now defined
on the space of three-pointed compact metric spaces (see e.g.\ \cite[Section 2.1]{Disks}), we have 
\begin{equation}
\label{inter-3}
\N_0\Big(\int  \frac{\mathrm{Vol}(\dd x)}{\sigma}F(\bm_\infty,x_*,x_0,x))\Big)
=\N_0\Big(\int\!\!\int\!\!\int \frac{\mathrm{Vol}(\dd x)}{\sigma}\,\frac{\mathrm{Vol}(\dd y)}{\sigma}\,\frac{\mathrm{Vol}(\dd z)}{\sigma}
\,F(\bm_\infty,x,y,z)\Big).
\end{equation}

\subsection{Moments of exit measures and volumes of balls}

We start with a lemma providing bounds on moments of the volume
of balls centered at $x_*$.

\begin{lemma}
\label{volume-ball}
Let $p\geq 1$ be an integer. There exists a constant $C_p$ such that,
for every $r>0$,
$$\N_0\Big(\mathrm{Vol}(B_r(x_*))^p\Big)= C_p\,r^{4p-2}.$$
Consequently, for every integer $p\geq 1$, and every $\eta\in(0,1)$,
there exists a constant $C_{p,\eta}$ such that, for every $r\in(0,1)$,
$$\N_0^{[1]}\Big(\mathrm{Vol}(B_r(x_*))^p\Big) \leq C_{p,\eta}\,r^{4p-\eta}.$$
\end{lemma}

\proof Using \eqref{equal-D} and the scaling property of the measures $\N_x$, we get
$$\N_0\Big(\mathrm{Vol}(B_r(x_*))^p\Big)
=\N_0\Big(\Big(\int_0^\sigma \dd s\,\mathbf{1}_{\{\wh W_s-W_*\leq r\}}\Big)^p\Big)
=r^{4p-2}\,\N_0\Big(\Big(\int_0^\sigma \dd s\,\mathbf{1}_{\{\wh W_s-W_*\leq 1\}}\Big)^p\Big).$$
So we only need to verify that
\begin{equation}
\label{volume-tech}
\N_0\Big(\Big(\int_0^\sigma \dd s\,\mathbf{1}_{\{\wh W_s-W_*\leq 1\}}\Big)^p\Big)<\infty.
\end{equation}
To this end, we note that, for every $\delta\in(0,1)$, \cite[Lemma 6.1]{Acta} gives
$$\N^{(1)}_0\Big(\Big(\int_0^1 \dd s\,\mathbf{1}_{\{\wh W_s-W_*\leq r\}}\Big)^p\Big)\leq c_{p,\delta}\,r^{4p-\delta},$$
with a constant $c_{p,\delta}$ independent of $r>0$. By scaling, we get for $t\geq 1$,
$$\N_0^{(t)}\Big(\Big(\int_0^t \dd s\,\mathbf{1}_{\{\wh W_s-W_*\leq 1\}}\Big)^p\Big)=t^p\,
\N^{(1)}_0\Big(\Big(\int_0^1 \dd s\,\mathbf{1}_{\{\wh W_s-W_*\leq t^{-1/4}\}}\Big)^p\Big)\leq c_{p,\delta}\,t^{\delta/4},$$
and \eqref{volume-tech} follows by writing
\begin{align*}
\N_0\Big(\Big(\int_0^\sigma \dd s\,\mathbf{1}_{\{\wh W_s-W_*\leq 1\}}\Big)^p\Big)&= \int_0^\infty \frac{\dd t}{2\sqrt{2\pi t^3}} \N_0^{(t)}\Big(\Big(\int_0^t \dd s\,\mathbf{1}_{\{\wh W_s-W_*\leq 1\}}\Big)^p\Big)\\
&\leq \int_0^1 \frac{\dd t}{2\sqrt{2\pi t^3}}\,t^p
+\int_1^\infty \frac{\dd t}{2\sqrt{2\pi t^3}} c_{p,\delta}\,t^{\delta/4}.
\end{align*}

For the second assertion, let $q\geq 2$ be an integer. Since $\N_0(W_*<-1)=\frac{3}{2}$, we have
$$\N_0^{[1]}\Big(\mathrm{Vol}(B_r(x_*))^p\Big)\leq \Big(\N_0^{[1]}\Big(\mathrm{Vol}(B_r(x_*))^{qp}\Big)\Big)^{1/q}\leq \Big(\frac{2}{3}\N_0\Big(\mathrm{Vol}(B_r(x_*))^{qp}\Big)\Big)^{1/q}
\leq (C_{qp})^{1/q}\,r^{4p-2/q},$$
which gives the desired result by taking $q$ such that $2/q<\eta$. \endproof

We will need bounds on the moments of $\z_{W_*+r}$ under $\N^{[r]}_0=\N_0(\cdot\midd W_*<-r)$, 
and to this end we will use 
a coupling with the Brownian plane. Let $\pp$ stand for the Brownian plane as defined in \cite{Plane,CLG}. According to
\cite[Proposition 1.1]{CLG}, we can make sense of a quantity $Z_r$
which corresponds to the boundary size of the hull of radius $r$ in $\pp$ (see the introduction of \cite{CLG}
for more details).

\begin{lemma} 
\label{coupling-exit}
Let $r>0$ and $u\geq r$. The random variable $\z_{W_*+r}$ under $\N^{[u]}_0$ is stochastically dominated by
$Z_r$.
\end{lemma}

This lemma is a straightforward consequence of a coupling between the Brownian plane $\pp$
and the Brownian sphere $\bm_\infty$ under $\N_0(\cdot\mid W_*=-u)$, which relies on the spine
decomposition of Section \ref{spine-decomp} and is described
in the proof of another technical lemma in the Appendix below. For this reason, we
also postpone the proof of Lemma \ref{coupling-exit} to the Appendix.

According to \cite[Proposition 1.2]{CLG}, the variable $Z_r$
follows the Gamma distribution with parameter $3/2$ and mean $r^2$.
In particular, for every $p\geq 1$, there exists a constant $c_p$ such that
$\EE[(Z_r)^p]=c_p\,r^{2p}$, and then Lemma \ref{coupling-exit} implies that,
for every $0<r\leq u$,
\begin{equation}
\label{bd-mom-exit}
\N^{[u]}_0\Big((\z_{W_*+r})^p\Big) \leq c_p\,r^{2p}.
\end{equation}

\section{Hulls}
\label{sec-hull}

\subsection{The construction of hulls}
\label{cons-hull}

Let us start by defining the random compact metric space which will
correspond to the conditional distribution of 
a hull of radius $r$ in the free Brownian sphere, given its boundary size.
Throughout this section,  $r>0$ and $z>0$ are fixed. We let
$$\n=\sum_{i\in I} \delta_{(t_i,\omega_i)}$$
be a Poisson point measure on $[0,z]\times \S_0$ with intensity
$$\dd t\otimes \N_0(\dd \omega \cap\{W_*>-r\}).$$
Futhermore, we let $\omega_*$ be a random snake trajectory
distributed according to $\N_0(\cdot\midd W_*=-r)$, and 
we let $U_*$ be uniformly distributed over $[0,z]$. We assume that
$\omega_*,U_*$ and $\nn$ are independent. We also set
$$\Sigma:=\sigma(\omega_*)+\sum_{i\in I} \sigma(\omega_i).$$

We let $\mathbf{H}$ be derived from the disjoint union 
$$[0,z] \cup \Big(\bigcup_{i\in I} \t_{(\omega_{i})}\Big) \cup \t_{(\omega_*)}$$
by identifying $0$ with $z$, the root of $\t_{(\omega_*)}$
with the point $U_*$ of $[0,z]$ and, for every $i\in I$, the root of $\t_{(\omega_{i})}$ 
with the point $t_i$ of $[0,z]$. The volume measure on $\mathbf{H}$
is defined by saying that it puts no mass on $[0,z]$ and that its restriction
to each of the trees $\t_{(\omega_i)}$, resp. to $\t_{(\omega_*)}$, is
the volume measure on this tree. 

We then assign labels $(\Lambda_a)_{a\in\mathbf{H}}$ to the 
points of $\mathbf{H}$. We take $\Lambda_a=0$ if $a\in [0,z]$, and, for every $a\in\t_{(\omega_i)}$,
resp. $a\in\t_{(\omega_*)}$, we take $\Lambda_a=\ell_{a}(\omega_i)$, resp. $\Lambda_a=\ell_a(\omega_*)$.
We let $b_*\in \t_{(\omega_*)}$
be the unique point of $\mathbf{H}$ with label $-r$. 

We finally define a cyclic exploration $(\ee_s)_{s\in [0,\Sigma]}$ of $\mathbf{H}$. Informally, this
cyclic exploration is obtained by concatenating the cyclic explorations of the $\t_{(\omega_i)}$'s, and of 
$\t_{(\omega_*)}$, in the order prescribed by the reals $t_i$, and $U_*$. To get a more precise
description, set 
$$A_u=\sum_{i\in I,t_i\leq u} \sigma(\omega_i)+\mathbf{1}_{\{u\geq U_*\}}\sigma(\omega_*)$$
for every $u\in[0,z]$, and use the notation $A_{u-}$ for the left limit at $u$. Then:
\begin{description}
\item[$\bullet$] If $0\leq s\leq A_{U_*-}$ or $A_{U_*}\leq s\leq \Sigma$, either there is a (unique) $i\in I$ 
such that $A_{t_i-}\leq s\leq A_{t_i}$, and we set $\ee_s=p_{(\omega_i)}(s-A_{t_i-})$, or
there is no such $i$ and we set $\ee_s=\sup\{t_i:A_{t_i}\leq s\}\in[0,z]$.
\item[$\bullet$] If $A_{U_*-}<s<A_{U_*}$, we set $\ee_s=p_{(\omega_*)}(s-A_{U_*-})$.
\end{description}
There is a unique $s_*\in[0,\Sigma]$ such that $\ee_{s_*}=b_*$. 

If $s,s'\in[0,\Sigma]$ and $s<s'$, we make the convention that $[s',s]=[s',\Sigma]\cup[0,s]$
(and of course $[s,s']$ is the usual interval). Then, for every $a,b\in \mathbf{H}$, 
we can find $s,s'\in[0,\Sigma]$ such that $\ee_s=a$, $\ee_{s'}=b$ and 
$[s,s']$ is as small as possible (note that we may have $s'<s$). We let the
interval $[a,b]$ of $\mathbf{H}$ be defined by $[a,b]:=\{\ee_r:r\in [s,s']\}$. We then define,
for every $a,b\in\mathbf{H}$,
$$D^\circ_{\mathbf{H}}(a,b):=\Lambda_a+\Lambda_b-2\,\max\Big(\inf_{c\in[a,b]} \Lambda_c,\inf_{c\in[b,a]} \Lambda_c\Big),$$
and 
\begin{equation}
\label{dist-hull}
D_{\mathbf{H}}(a,b)=\inf_{a_0=a,a_1,\ldots,a_{k-1},a_k=b} \sum_{i=1}^k D_{\mathbf{H}}^\circ(a_{i-1},a_i).
\end{equation}
Then, almost surely for every $a,b\in\mathbf{H}$, the property
$D_{\mathbf{H}}(a,b)=0$ holds if and only if $D_{\mathbf{H}}^\circ(a,b)=0$. This is
derived from the analogous property for the Brownian map, which is proved in \cite{Invent}.
The bound $D^\circ_{\mathbf{H}}(a,b)\geq |\Lambda_a -\Lambda_b|$ is immediate from the
definition, and it follows that we have also $D_{\mathbf{H}}(a,b)\geq |\Lambda_a -\Lambda_b|$.

The mapping $(a,b)\mapsto D_{\mathbf{H}}(a,b)$ is a pseudo-metric on $\mathbf{H}$. We writte
$\mathfrak{H}:=\mathbf{H}/\{D_{\mathbf{H}}=0\}$ for the associated quotient space, which is
equipped with the metric (induced by) $D_{\mathbf{H}}$. We write $\Pi_\mathfrak{H}$ for the canonical 
projection from $\mathbf{H}$ onto $\mathfrak{H}$.
The restriction of $\Pi_\mathfrak{H}$ to $[0,z)$ is continuous and one-to-one, and 
its range is a simple loop denoted by $\partial\mathfrak{H}=\Pi_\mathfrak{H}([0,z))$ 
(one proves that $\mathfrak{H}$ is homeomorphic to the closed unit disk, and 
via such a homeomorphism $\partial\mathfrak{H}$ indeed corresponds to the
unit circle).
The volume measure on $\mathfrak{H}$ is the pushforward of the volume measure 
on $\mathbf{H}$ under $\Pi_\mathfrak{H}$. 
We view $\mathfrak{H}$ as a random variable with values in the space $\M^{\bullet b}$
of Section \ref{sec-GHP}:
the distinguished point is $\Pi_\mathfrak{H}(b_*)$, and the distinguished ``boundary'' is $\partial\mathfrak{H}$.
Without risk of confusion, we identify $b_*$ and $\Pi_\mathfrak{H}(b_*)$ in what follows. 

We call the random two-boundary measure metric space
$\mathfrak{H}$  the  {\bf standard hull of radius $r$ and perimeter $z$}. This terminology will be justified
below by relations with the Brownian sphere. If $x\in\mathfrak{H}$, we can set $\Lambda_x=\Lambda_a$,
for $a$ such that $\Pi_\mathfrak{H}(a)=x$ (this does not depend on the choice of $a$) and it easily follows from the definitions that
\begin{equation}
\label{dis-root}
D_{\mathbf{H}}(x,b_*)=\Lambda_x+r
\end{equation}
for every $x\in\mathfrak{H}$. In particular, all points of $\partial\mathfrak{H}$ are at distance $r$ from
$b_*$.

Let us turn to geodesics in $\mathfrak{H}$. More precisely, we are interested in geodesics 
between an arbitrary point of $\mathfrak{H}$ and $b_*$. Let $x\in\mathfrak{H}$, and 
let $a\in\mathbf{H}$ such that $\Pi_\mathfrak{H}(a)=x$, and $s\in[0,\Sigma]$ such that
$\ee_s=a$. Consider first the case where $s\in[0,s_*]$. We then set, for every $t\in[0,\Lambda_a+r]$,
$$\gamma_s(t):=\Pi_\mathfrak{H}\Big(\ee_{\inf\{u\geq s:\Lambda_{\ee_u}=\Lambda_a-t\}}\Big).$$
If $s\in[s_*,\Sigma]$, we define similarly, for every $t\in[0,\Lambda_a+r]$,
$$\gamma_s(t):=\Pi_\mathfrak{H}\Big(\ee_{\sup\{u\leq s:\Lambda_{\ee_u}=\Lambda_a-t\}}\Big).$$
Then, using \eqref{dis-root} and the bound  $D_{\mathbf{H}}(a,b)\geq |\Lambda_a -\Lambda_b|$, it is straightforward to verify that $(\gamma_s(t))_{0\leq t\leq \Lambda_a+r}$ is a geodesic
from $x$ to $b_*$ in $\mathfrak{H}$. Such a geodesic is called a {\it simple geodesic}.

\begin{proposition}
\label{simple-geo}
All geodesics in $\mathfrak{H}$ that end at $b_*$ are simple geodesics.
\end{proposition}

The analog of this result for the Brownian sphere is proved in \cite{Acta}. The proposition
can be derived from this analog by using the relations with the Brownian sphere that
are discussed below. 

\begin{proposition}
\label{distinct-geo}
Let $\ve\in(0,r)$. 
Define an integer-valued random
variable $N_\ve$ by saying that $N_\ve\geq k$ if and only if there
exist $k$ geodesics $\phi_1,\phi_2,\ldots,\phi_k$
from $\partial\mathfrak{H}$ to $b_*$ such that 
the sets $\{\phi_1(t):0\leq t\leq r-\ve\}$, $\{\phi_2(t):0\leq t\leq r-\ve\},
\ldots,\{\phi_k(t):0\leq t\leq r-\ve\}$ are disjoint. Then $N_\ve-1$
follow a Poisson distribution with parameter
$$\frac{3z}{2}\Big( (r-\ve)^{-2}-r^{-2})\Big).$$
\end{proposition}

\proof Let $u,v\in[0,z]$, and $s,s'\in[0,\Sigma]$ such that $\ee_s=u$
and $\ee_{s'}=v$. Consider the simple geodesics $\gamma_{s}$
and $\gamma_{s'}$ as defined above. Then it follows 
from this definition (using also the fact that $D_{\mathbf{H}}(a,b)=0$ if and only if $D_{\mathbf{H}}^\circ(a,b)=0$)
that the sets $\{\gamma_s(t):0\leq t\leq r-\ve\}$ and $\{\gamma_{s'}(t):0\leq t\leq r-\ve\}$
are disjoint if and only if we have both
$$\min_{t\in[s,s']} \Lambda_{\ee_t} <-r+\ve\hbox{ \ and \ }\min_{t\in[s',s]} \Lambda_{\ee_t} <-r+\ve.$$
Using also Proposition \ref{simple-geo}, we get that
\begin{equation}
\label{nber-exc}
N_\ve-1=\#\{i\in I: W_*(\omega_i)<-r+\ve\}.
\end{equation}
which follows a Poisson distribution with parameter
$$z\,\N_0(-r<W_*<-r+\ve)= \frac{3z}{2}\Big( (r-\ve)^{-2}-r^{-2})\Big).$$
This completes the proof. \endproof

\subsection{Hulls in the Brownian sphere} 
\label{hulls-sphere}

We now consider the free Brownian sphere $(\bm_\infty,D)$,
which is defined under the measure $\N_0$. Recall that the
 two distinguished points of $\bm_\infty$ are $x_*$ and $x_0$.

We will be 
interested in hulls centered at $x_*$ relative to $x_0$ in $\bm_\infty$. We
write $B^\bullet_r(x_*)=B^{\bullet(x_0)}_r(x_*)$ to simplify notation: this hull is defined
on the event where $D(x_0,x_*)=-W_*>r$. 
We notice the following useful fact, which is obtained from the cactus bound \eqref{cactus-bd}
in a way similar to the proof of formulas (16) and (17) in \cite{CLG}. A point $x=\Pi(p_{(\omega)}(s))$, for
$s\in [0,\sigma]$, belongs to $B^\bullet_r(x_*)$ if and only if $\tau_{W_*+r}(W_s)\leq \zeta_s$, and 
to the interior of $B^\bullet_r(x_*)$ if and only if $\tau_{W_*+r}(W_s)<\zeta_s$. 

The exit measure $\z_{W_*+r}$,
which was introduced in Section \ref{spine-decomp} and is also defined on the event $\{W_*<-r\}$, 
can be
interpreted as the boundary size of $B^\bullet_r(x_*)$. This interpretation is justified 
by the following approximation, which is a reformulation of \eqref{formu-exit2} using
the preceding observations,
\begin{equation}
\label{approx-exit}
\z_{W_*+r}=\lim_{\ve\to 0} \ve^{-2}\mathrm{Vol}((B^\bullet_r(x_*))^c\cap B_{r+\ve}(x_*))\;,
\hbox{ a.e. on }\{D(x_*,x_0)>r\}.
\end{equation}

If $O$ is a connected open subset of $\bm_\infty$, the intrinsic distance $D^O_{\rm int}$
on $O$ is defined as follows. For $x,y\in O$, $D^O_{\rm int}(x,y)$ is the infimum 
of the lengths of paths connecting $x$ and $y$ and staying inside $O$.

Recall that $\N_0^{[r]}=\N_0(\cdot\midd D(x_*,x_0)>r)$. We write $\ov{\bm_\infty\backslash B^\bullet_r(x_*)}$
for the closure of $\bm_\infty\backslash B^\bullet_r(x_*)$.

\begin{theorem}
\label{deco-hull}
With $\N_0^{[r]}$-probability one, the intrinsic distance on the interior of $B^\bullet_r(x_*)$ has a continuous
extension to $B^\bullet_r(x_*)$, which is a metric on $B^\bullet_r(x_*)$, and similarly the intrinsic distance on  $\bm_\infty\backslash B^\bullet_r(x_*)$
has a continuous extension to $\ov{\bm_\infty\backslash B^\bullet_r(x_*)}$, which is a metric on this space.
Consider both $B^\bullet_r(x_*)$ and $\ov{\bm_\infty\backslash B^\bullet_r(x_*)}$ as metric spaces for these
(extended) intrinsic metrics.
The metric space  $B^\bullet_r(x_*)$ equipped with the restriction of the volume measure on $\bm_\infty$, with the distinguished point $x_*$
and with the distinguished boundary $\partial B^\bullet_r(x_*)$ is a random
element of $\M^{\bullet b}$ and the same holds for the metric space $\ov{\bm_\infty\backslash B^\bullet_r(x_*)}$ 
equipped with the restriction of the volume measure, with the distinguished point $x_0$
and with the distinguished boundary $\partial B^\bullet_r(x_*)$. 
Then, for any nonnegative measurable functions $F$ and $G$ defined on $\M^{\bullet b}$, for every $z>0$, we have
$$\N_0^{[r]}\Big( F(B^\bullet_r(x_*))\, G\big(\ov{\bm_\infty\backslash B^\bullet_r(x_*)}\big)\,\Big|\, \z_{W_*+r}=z\Big)
= \E[F(\mathfrak{H}_{r,z})]\,\E[G(\D^\bullet_z)],$$
where $\mathfrak{H}_{r,z}$ stands for the standard hull of radius $r$ and perimeter $z$, and $\D^\bullet_z$
denotes a free pointed Brownian disk of perimeter $z$.
\end{theorem}

\rems (i) The theorem implies in particular that $B^\bullet_r(x_*)$ and $\ov{\bm_\infty\backslash B^\bullet_r(x_*)}$ are independent under $\N^{[r]}_0$ conditionally on $\z_{W_*+r}$. If now we take $r'>r$, it remains 
true that $B^\bullet_r(x_*)$ and $\ov{\bm_\infty\backslash B^\bullet_r(x_*)}$ are independent under $\N^{[r']}_0$ conditionally on $\z_{W_*+r}$. The point is that, if one already knows that $D(x_*,x_0)>r$, the event $D(x_*,x_0)>r'$ occurs if and
only if the distance from the distinguished point of $\ov{\bm_\infty\backslash B^\bullet_r(x_*)}$ to the boundary
is greater than $r'-r$, which only depends on $\ov{\bm_\infty\backslash B^\bullet_r(x_*)}$. 

\smallskip
\noindent(ii) The free pointed Brownian disk $\D^\bullet_z$ in the theorem is also viewed as a random variable with values in $\M^{\bullet b}$
(the ``boundary'' is of course the usual boundary $\partial\D^\bullet_z$, see e.g.\ \cite{BM}). 

\proof This is very similar to the proof of Theorems 29 and 31 of \cite{Spine}, which give the
analogous statements for the Brownian plane (then $\D^\bullet_z$ is replaced by an 
infinite Brownian disk with perimeter $z$). For this reason, we only outline certain arguments, especially
in the final part of the proof, and we
refer to \cite{Spine} for more details. Throughout the proof, we argue 
under the measure $\N^{[r]}_0(\dd \omega)$. We consider the truncation $\tr_{W_*+r}(\omega)$ and
 the excursions $(\omega^j)_{j\in J}$ 
below level $W_*+r$ as defined in Section \ref{sna-tra}. According to Proposition 12 of \cite{Disks}, we know that, conditionally
on $\z_{W_*+r}=z$, $\sum_{j\in J}\delta_{\omega^j}$ is independent of $\tr_{W_*+r}(\omega)$ and distributed as
a Poisson point measure with intensity $z\,\N_0$ conditioned to have a minimum equal to
$-r$. We need in fact a slightly more precise result
involving also the exit local time at level $W_*+r$, which 
is denoted by $(L^{W_*+r}_s)_{s\geq 0}$. Notice that this exit local time is easily defined by using the spine decomposition in
Section \ref{spine-decomp}, and that $L^{W_*+r}_\sigma=\z_{W_*+r}$.

For every $j\in J$, let $(\alpha_j,\beta_j)$ be the time interval 
corresponding to the excursion $\omega^j$, and let $l_j=L^{W_*+r}_{\alpha_j}=L^{W_*+r}_{\beta_j}$.
Then an application of the special Markov property in the form given in the appendix
of \cite{subor} (using again the spine decomposition in
Section \ref{spine-decomp}) shows that, conditionally on $\z_{W_*+r}=z$, the point measure
\begin{equation}
\label{point-mea}
\sum_{j\in J} \delta_{(l_j,\omega^j)}
\end{equation}
is independent of $\tr_{W_*+r}(\omega)$ and has the same distribution as $\nn+\delta_{(U_*,\omega_*)}$
with the notation of Section \ref{cons-hull}. According to Section \ref{cons-hull}, the hull $\mathfrak{H}_{r,z}$
is constructed as a measurable function of $\nn+\delta_{(U_*,\omega_*)}$. 
We will then verify that, if we perform the construction of Section \ref{cons-hull}
from the point measure \eqref{point-mea}, 
we get a pointed compact metric space isometric to 
$B^\bullet_r(x_*)$ (equipped with its intrinsic metric). We let 
$(\mathfrak{H}^\star,D^\star)$ stand for the metric space obtained from the point measure in \eqref{point-mea}
by the construction of Section \ref{cons-hull} --- notice that this construction makes sense even with a random perimeter $z$.
We also use the notation
$\mathbf{H}^\star$ for the space constructed from the 
point measure \eqref{point-mea} in a way similar to $\mathbf{H}$
in Section \ref{cons-hull}, and $\Pi^\star$ for the canonical projection from $\mathbf{H}^\star$ 
onto $\mathfrak{H}^\star$. Note that $D^\star=D_{\mathbf{H}^\star}$ in the notation of Section \ref{cons-hull}.
We also define $D_{\mathbf{H}^\star}^\circ$ as
in Section \ref{cons-hull} replacing $\mathbf{H}$ by $\mathbf{H^\star}$.

We first explain that
the set $\mathfrak{H}^\star$ can be identified to $B^\bullet_r(x_*)$. To this end, set
$$F_r:=\{p_{(\omega)}(s):0\leq s\leq \sigma, \tau_{W_*+r}(W_s)\leq \zeta_s\}\;,\ 
\partial F_r:=\{p_{(\omega)}(s):0\leq s\leq \sigma, \tau_{W_*+r}(W_s)=\zeta_s\}.$$
As we already noticed, we have $B^\bullet_r(x_*)=\Pi(F_r)$ and $\partial B^\bullet_r(x_*)=\Pi(\partial F_r)$. We define a mapping $\ii$
from $F_r$ onto $\mathbf{H}^\star$ by the following prescriptions. If $p_{(\omega)}(s)\in F_r\backslash \partial F_r$, then 
$s$ belongs to $(\alpha_j,\beta_j)$ for some $j\in J$, and we take $\ii(p_{(\omega)}(s))=p_{(\omega^j)}(s-\alpha_j)\in \t_{(\omega^j)}\subset \mathbf{H}^\star$.
On the other hand, if $p_{(\omega)}(s)\in \partial F_r$, we take $\ii(p_{(\omega)}(s))=L^{W_*+r}_s\in [0,\z_{W_*+r}]=\partial \mathbf{H}^\star$. 
The reader will easily check that $\ii(a)$ is well-defined for every $a\in F_r$ independently of the choice of $s$ such that $p_{(\omega)}(s)=a$. 
Moreover, we have $D^\circ_{\mathbf{H}^\star}(\ii(a),\ii(b))=D^\circ(a,b)$ for every $a,b\in F_r$ (we omit the details). The mapping
$a\mapsto \ii(a)$ is not one-to-one (though its restriction to $F_r\backslash \partial F_r$ is one-to-one) but the latter property
shows that $\Pi(a)=\Pi(b)$ implies $\Pi^\star(\ii(a))=\Pi^\star(\ii(b))$. So one can define a mapping $\jj$ from $\Pi(F_r)=B^\bullet_r(x_*)$ onto $\mathfrak{H}^\star$
by declaring that $\jj(\Pi(a))=\Pi^\star(\ii(a))$ for every $a\in F_r$. This mapping $\jj$ is one-to-one since $\Pi^\star(\ii(a))=\Pi^\star(\ii(b))$
is only possible if $D^\circ_{\mathbf{H}^\star}(\ii(a),\ii(b))=0$, which implies $D^\circ(a,b)=0$ and $\Pi(a)=\Pi(b)$. The mapping $\jj$
provides the desired identification of $B^\bullet_r(x_*)$ with $\mathfrak{H}^\star$, and we also observe that $\partial B^\bullet_r(x_*)$ is 
identified with $\partial\mathfrak{H}^\star$. From now one we make these identifications, and we notice that
we have $D(x,y)\leq D^\star(x,y)$ for $x,y\in B^\bullet_r(x_*)$, by comparing formulas \eqref{formulaD} and \eqref{dist-hull},
using the equality $D^\circ_{\mathbf{H}^\star}(\ii(a),\ii(b))=D^\circ(a,b)$ for $a,b\in F_r$ (the point is that there
are more choices for the intermediate points $a_1,\ldots,a_{k-1}$ in \eqref{formulaD} than in \eqref{dist-hull}).

Then, we need to
verify that the restriction of $D^\star$ to the interior of $B^\bullet_r(x_*)$ coincides with
the intrinsic distance, which we denote by $D^{\rm intr}$. The bound $D^\star\leq D^{\rm intr}$
is easy. If $(\gamma(t))_{0\leq t\leq 1}$ is a path connecting two points $x$ and $y$ of 
the interior of $B^\bullet_r(x_*)$ that stays in this interior, the length of $\gamma$ is bounded 
below by $\sum_{k=1}^n D(\gamma(t_{k-1}),\gamma(t_k))$, where $0=t_0<t_1<\cdots<t_n=1$
is a subdivision of $[0,1]$. For every $k\in\{0,\ldots, n\}$, let $a_k\in\t_{(\omega)}$ be such that
$\Pi(a_k)=\gamma(t_k)$.  Then,
\begin{equation}
\label{deco-tech12}
D(\gamma(t_{k-1}),\gamma(t_k))=\build{\inf_{c_0=a_{k-1},c_1,\ldots,c_p=a_k}}_{c_1,\ldots,c_{p-1}\in\t_{(\omega)}}^{} \sum_{j=1}^p D^\circ(c_{j-1},c_j)
\end{equation}
but if the mesh of the subdivision is sufficiently small (so that all $D(\gamma(t_{k-1}),\gamma(t_k))$
are small) we can assume that the infimum of the previous display is attained by
considering only points $c_j$ such that $\Pi(c_j)\in B^\bullet_r(x_*)$ (otherwise the sum in the right-hand side of
\eqref{deco-tech12} is bounded below by the $D$-distance between the range of $\gamma$ and $\partial B^\bullet_r(*)$). For such a choice of the $c_j$'s,
we have $D^\circ(c_{j-1},c_j)=D^\circ_{\mathbf{H}^\star}(\ii(c_{j-1}),\ii(c_j))$. It
follows that $D(\gamma(t_{k-1}),\gamma(t_k))\geq D^\star(\gamma(t_{k-1}),\gamma(t_k))$, and finally that
the length of $\gamma$ is bounded below by $D^\star(x,y)$ as desired.

The reverse bound $D^{\rm intr}\leq D^\star$ is slightly more delicate, and we only sketch the
argument. The difficulty comes from the following observation. In formula \eqref{dist-hull}
giving $D^\star$ in terms of $D^\circ_{\mathbf{H}^\star}$, even 
if $a$ and $b$ do not belong to the boundary $\partial\mathbf{H}^\star$, we need a priori to consider points $a_1,\ldots,a_{k-1}\in \mathbf{H}^\star$
that may belong to this boundary. However, we leave it as an exercise for the reader to check that
the infimum remains the same even if we impose that all points $a_1,\ldots,a_{k-1}$
do not lie on $\partial\mathbf{H}^\star$. In that case, $D^\circ_{\mathbf{H}^\star}(a_{i-1},a_i)$ can be interpreted as the
length of a path from $\Pi^\star(a_{i-1})$ to $\Pi^\star(a_i)$ made of the concatenation of
two simple geodesics started respectively from $\Pi^\star(a_{i-1})$ and from $\Pi^\star(a_i)$, see
e.g.\ the end of \cite[Section 4.1]{Spine} for a very similar argument. In the 
identification of $\mathfrak{H}^\star$ with $B^\bullet_r(x_*)$, these simple geodesics remain
geodesics for the distance $D$ and stay in the interior of $B^\bullet_r(x_*)$. Summarizing, we
obtain that $D^\star(x,y)$ is obtained as an infimum of quantities that are lengths of 
paths connecting $x$ to $y$ and staying in the interior of $B^\bullet_r(x_*)$. This gives the 
desired bound $D^{\rm intr}\leq D^\star$.

Once we know that $D^{\rm intr}= D^\star$ in the interior of $B^\bullet_r(x_*)$, the fact that
$D^{\rm intr}$ has a continuous extension to the boundary is easy. Suppose that $(x_n)_{n\in\N}$
and $(y_n)_{n\in\N}$ are two sequences in the interior of $B^\bullet_r(x_*)$ that converge 
to $x$ and $y$ respectively (for the metric $D$). We have to verify that $D^{\rm intr}(x_n,y_n)$ has a limit
as $n\to\infty$. We observe that the convergence of the sequence $(x_n)_{n\in\N}$ for the 
metric $D$ also implies that $D^\star(x_n,x)\la 0$ as $n\to\infty$. Indeed, by compactness, we can find a 
subsequence $(x_{n_k})_{k\in\N}$ and $x'\in B^\bullet_r(x_*)$ such that $D^\star(x_{n_k},x') \la 0$.
But since $D\leq D^\star$ on $B^\bullet_r(x_*)$, this readily implies that $x'=x$
and we get that $D^\star(x_n,x)\la 0$ as $n\to\infty$. We have similarly $D^\star(y_n,y)\la 0$ as $n\to\infty$,
and we conclude that $D^{\rm intr}(x_n,y_n)=D^\star(x_n,y_n)$ converges to $D^\star(x,y)$. 

At this stage we have proved that the intrinsic distance on the interior of $B^\bullet_r(x_*)$ has a continuous
extension to $B^\bullet_r(x_*)$, and that, conditionally on $\z_{W_*+r}=z$, the resulting random (pointed)
metric space $B^\bullet_r(x_*)$ has the same distribution as $\mathfrak{H}_{r,z}$ and is independent
of $\tr_{W_*+r}(\omega)$. To complete the proof, we need to verify that that the intrinsic metric on
$\bm_\infty\backslash B^\bullet_r(x_*)$
has a continuous extension to $\overline{\bm_\infty\backslash B^\bullet_r(x_*)}$ and that the (two-boundary measure)
metric space
$\overline{\bm_\infty\backslash B^\bullet_r(x_*)}$ is a function of $\tr_{W_*+r}(\omega)$, which,
conditionally on $\z_{W_*+r}=z$, is distributed as $\D_z$. The proof of the first assertion proceeds by minor
modifications of the proofs of \cite[Theorem 28]{Disks} or of \cite[Theorem 29]{Spine}, and we omit the details. 
As for the second assertion, we rely on \cite[Proposition 12]{Disks}, which gives the conditional 
distribution of $\tr_{W_*+r}(W)$ given $\z_{W_*+r}=z$: For every nonnegative measurable function
$F$ on $\S_0$,
\begin{equation}
\label{law-trunca}
\N^{[r]}_0\Big( F(\mathrm{tr}_{W_*+r}(\omega)) \,\Big|\, \z_{W_*+r}=z\Big)
= z^{-2}\, \N^{*,z}_0\Big( \int_0^\sigma \mathrm{d}s\,F(\omega^{[s]})\Big),
\end{equation}
where we recall that $\omega^{[s]}$ stands for $\omega$ re-rooted at $s$, and $\N^{*,z}_0$ is the law
of a positive Brownian snake excursion with boundary size $z$, as
defined in \cite{ALG}. Write $\wt{\mathbb{N}}^{*,z}_0$
for the probability measure on $\S_0$ such that $\wt{\mathbb{N}}^{*,z}_0(F)$
is the right-hand side of \eqref{law-trunca}. From the beginning of Section 4.3 in \cite{Spine},
the free pointed Brownian disk
$\D_z$ can be constructed as a measurable function of a random snake trajectory distributed 
according to $\wt{\mathbb{N}}^{*,z}_0$. We leave it to the reader to check that the same measurable function applied to
$\mathrm{tr}_{W_*+r}(\omega)$ yields the space
$\overline{\bm_\infty\backslash B^\bullet_r(x_*)}$
equipped with its extended intrinsic metric --- here again arguments are very similar to the
proof of \cite[Theorem 29]{Spine}. It thus follows from \eqref{law-trunca}
that $\overline{\bm_\infty\backslash B^\bullet_r(x_*)}$ under $\N^{[r]}_0( \cdot \midd\z_{W_*+r}=z)$
has the distribution of $\D_z$. This completes the proof.
 \endproof

We will need a ``two-point version'' of Theorem \ref{deco-hull}, which we now state 
as a corollary. It  is convenient
to write
$$\z^{x_*(x_0)}_r=\z_{W_*+r}$$
for the boundary size of the hull $B^\bullet_r(x_*)$, which makes sense on the event where
$D(x_*,x_0)>r$ or equivalently $W_*<-r$. 
By interchanging the roles of
$x_*$ and $x_0$ and relying on the symmetry properties of the Brownian sphere (cf.~\eqref{inter-2}), we can also define on the same event the quantity
$\z^{x_0(x_*)}_r$ now corresponding to the
boundary size of the hull $B^{\bullet(x_*)}_r(x_0)$ --- one may use the analog of the approximation formula \eqref{approx-exit}.

On the event where $D(x_*,x_0)>2r$, the hulls $B^{\bullet(x_0)}_r(x_*)$ and $B^{\bullet(x_*)}_r(x_0)$ are disjoint, and we set 
$$\mathcal{C}^{x_*,x_0}_r:= \bm_\infty\backslash \big(B^{\bullet(x_0)}_r(x_*)\cup B^{\bullet(x_*)}_r(x_0)\big).$$

In the next corollary, we view both $B^{\bullet(x_0)}_r(x_*)$ and $B^{\bullet(x_*)}_r(x_0)$ equipped with their (extended) intrinsic metrics as random variables 
with values in $\M^{\bullet b}$ as stated in Theorem \ref{deco-hull} (the fact that this is also
legitimate for $B^{\bullet(x_*)}_r(x_0)$ is a consequence of \eqref{inter-2}).
We use the notation $\Theta_{r,z}$ for the distribution
of the standard hull $\mathfrak{H}_{r,z}$ of Section \ref{cons-hull}, so that $\Theta_{r,z}$ is a probability 
measure on the space $\mathbb{M}^{\bullet b}$.

\begin{corollary}
\label{two-hulls}
A.e.~under $\N_0(\cdot\cap \{D(x_*,x_0)>2r\})$,
the intrinsic metric on $\mathcal{C}^{x_*,x_0}_r$ has a continuous extension to its closure $\overline{\mathcal{C}^{x_*,x_0}_r}$,
which is a metric on this space. We equip this metric space with the restriction of the volume measure of $\bm_\infty$
and with the boundaries $\partial  B^{\bullet(x_0)}_r(x_*)$ and $\partial B^{\bullet(x_*)}_r(x_0)$, so that we view $\overline{\mathcal{C}^{x_*,x_0}_r}$
as a random variable with values in $\M^{bb}$.
Then, if
$F_1$ and $F_2$ are two nonnegative measurable functions on $\mathbb{M}^{\bullet b}$,
and $G$ is  a nonnegative measurable function on $\mathbb{M}^{bb}$, we have
$$\N_0^{[2r]}\Big(F_1(B^{\bullet(x_0)}_r(x_*))\,F_2(B^{\bullet(x_*)}_r(x_0))\,G\big(\overline{\mathcal{C}^{x_*,x_0}_r}\big)\Big)
=\N_0^{[2r]}\Big(\Theta_{r,\z^{x_*(x_0)}_r}(F_1)\,\Theta_{r,\z^{x_0(x_*)}_r}(F_2)\,G\big(\overline{\mathcal{C}^{x_*,x_0}_r}\big)\Big).$$
\end{corollary}

\proof
Write $D^{\mathrm{intr},\cc}$ for the intrinsic distance on $\mathcal{C}^{x_*,x_0}_r$. We need to verify that,
if $(x_n)_{n\in\N}$ and $(y_n)_{n\in\N}$ are two sequences in $\mathcal{C}^{x_*,x_0}_r$ that converge 
respectively to $x$ and $y$ belonging to $\overline{\mathcal{C}^{x_*,x_0}_r}$, then the sequence
$D^{\mathrm{intr},\cc}(x_n,y_n)$ converges. To this end, it suffices to prove that $D^{\mathrm{intr},\cc}(x_n,x_p)$
converges to $0$ as $n,p\to\infty$ (and similarly for $D^{\mathrm{intr},\cc}(y_n,y_p)$). If $x\in \mathcal{C}^{x_*,x_0}_r$
this is trivial, so we can suppose that $x\in\partial B^{\bullet(x_0)}_r(x_*)$ --- the case $x\in\partial B^{\bullet(x_*)}_r(x_0)$
is treated in a symmetric manner. However, writing $D^{\mathrm{intr},x_*}$ for the intrinsic distance on $\bm_\infty\backslash 
B^{\bullet(x_0)}_r(x_*)$, and recalling that this distance is extended 
continuously to the closure of $\bm_\infty\backslash 
B^{\bullet(x_0)}_r(x_*)$, we already know from Theorem \ref{deco-hull} that 
$D^{\mathrm{intr},x_*}(x_n,x_p)$
converges to $0$ as $n,p\to\infty$. Then the desired result follows from the fact that we have $D^{\mathrm{intr},\cc}(x_n,x_p)
=D^{\mathrm{intr},x_*}(x_n,x_p)$ as soon as $n,p$ are large enough. Indeed, for any $\ve>0$, if both
$D^{\mathrm{intr},x_*}(x_n,x)$ and $D^{\mathrm{intr},x_*}(x_p,x)$ are smaller than $\ve$, we have
$$D(x_n,B^{*(x_*)}(x_0))+D(x_p,B^{*(x_*)}(x_0))> 
2\min\{D(u,v):u\in B^{\bullet(x_0)}_r(x_*), v\in B^{\bullet(x_*)}_r(x_0)\}-2\ve,$$ and therefore, taking $\ve$ small enough, we see that,
when $n$ and $p$ are sufficiently large, the
 length of a path from $x_n$ to $x_p$ that hits
$B^{\bullet(x_*)}_r(x_0)$ is bounded below by a positive quantity. Since $D^{\mathrm{intr},x_*}(x_n,x_p)
\leq D^{\mathrm{intr},x_*}(x_n,x)+D^{\mathrm{intr},x_*}(x_p,x)$, which tends to $0$ as $n,p\to\infty$, it follows that
the infimum that gives $D^{\mathrm{intr},x_*}(x_n,x_p)$ must be attained for paths that do not hit
$B^{\bullet(x_*)}_r(x_0)$, and thus $D^{\mathrm{intr},\cc}(x_n,x_p)
=D^{\mathrm{intr},x_*}(x_n,x_p)$ when $n$ and $p$ are large enough. So $D^{\mathrm{intr},\cc}$ can be extended by 
continuity to $\overline{\mathcal{C}^{x_*,x_0}_r}$, and a similar argument shows that $D^{\mathrm{intr},\cc}(x,y)>0$
if $x$ and $y$ are distinct points of $\partial B^{\bullet(x_0)}_r(x_*)$ (resp., of $\partial B^{\bullet(x_*)}_r(x_0)$) since otherwise this
would imply $D^{\mathrm{intr},x_*}(x,y)=0$. 

Let us turn to the second assertion of the corollary. It follows from Theorem \ref{deco-hull},
that, for functions $F$ and $G$ as in this statement, 
$$
\N_0^{[r]}\Big( F(B^\bullet_r(x_*))\, G\big(\ov{\bm_\infty\backslash B^\bullet_r(x_*)}\big)\Big)\\
=\N_0^{[r]}\Big( \Theta_{r,\z^{x_*(x_0)}_r}(F)\, G\big(\ov{\bm_\infty\backslash B^\bullet_r(x_*)}\big)\Big)
$$
Up to replacing $G\big(\ov{\bm_\infty\backslash B^\bullet_r(x_*)}\big)$ by 
$\mathbf{1}_{\{D(x_0,\partial B^\bullet_r(x_*))>r\}}\,G\big(\ov{\bm_\infty\backslash B^\bullet_r(x_*)}\big)$,
we see that the preceding display remains valid if $\N_0^{[r]}$ is replaced by $\N^{[2r]}_0$. Next, under the
assumptions of the corollary, the quantity 
$$F_2(B^{\bullet(x_*)}_r(x_0))\,G\big(\overline{\mathcal{C}^{x_*,x_0}_r}\big)$$
is equal $\N_0^{[2r]}$ a.e.~to a measurable function of $\ov{\bm_\infty\backslash B^\bullet_r(x_*)}$, and so 
we get from the preceding considerations that
$$
\N_0^{[2r]}\Big(F_1(B^{\bullet(x_0)}_r(x_*))\,F_2(B^{\bullet(x_*)}_r(x_0))\,G\big(\overline{\mathcal{C}^{x_*,x_0}_r}\big)
)\Big)
=\N_0^{[2r]}\Big( \Theta_{r,\z^{x_*(x_0)}_r}(F_1)\,F_2(B^{\bullet(x_*)}_r(x_0))\,G\big(\overline{\mathcal{C}^{x_*,x_0}_r}\big)
\Big).
$$
At this stage, we use \eqref{inter-2} to interchange the roles
of $x_*$ and $x_0$. It follows that the preceding quantity is equal to
$$\N_0^{[2r]}\Big( \Theta_{r,\z^{x_0(x_*)}_r}(F_1)\,F_2(B^{\bullet(x_0)}_r(x_*))\,G\big(\overline{\mathcal{C}^{x_0,x_*}_r}\big)
\Big),$$
which, by the same application of Theorem \ref{deco-hull}, is also equal to
$$\N_0^{[2r]}\Big( \Theta_{r,\z^{x_0(x_*)}_r}(F_1)\,\Theta_{r,\z^{x_*(x_0)}_r}(F_2)\,G\big(\overline{\mathcal{C}^{x_0,x_*}_r}\big)
)\Big).$$
Finally, by interchanging once again the roles of $x_0$ and $x_*$ in the last display, we obtain the desired result.
\endproof

\section{The first-moment estimate}
\label{sec-first-mom}

\subsection{Slices}
\label{sec:slice}

In this section, we provide a brief description of the random 
compact metric spaces called slices, which apear as 
scaling limits of the planar maps with geodesic boundaries
considered in \cite{BM,Uniqueness} (notice that \cite{MQ} uses a slightly different 
definition of slices).
We fix $h>0$ and argue under $\N_0(\dd\omega\midd W_*=-h)$. Under
this probability measure, the Brownian sphere $\bm_\infty(\omega)$
is conditioned on the event that the distance between the two distinguished
points is equal to $h$. According to \cite{Acta},
the unique geodesic from $x_0$ to $x_*$ is the path $(\gamma(r))_{0\leq r\leq h}$ defined by
\begin{equation}
\label{main-geo}
\gamma(r)=\Pi\circ p_{(\omega)}(\inf\{s\in[0,\sigma]:\wh W_s=-r\})=\Pi\circ p_{(\omega)}(\sup\{s\in[0,\sigma]:\wh W_s=-r\}),
\end{equation}
for every $r\in[0,h]$.

We will now define another quotient space of 
$\t_{(\omega)}$ (later called a slice), which roughly speaking 
corresponds to cutting the Brownian sphere $\bm_\infty(\omega)$ along
the geodesic from $x_0$ to $x_*$ (see the end of Section 3.2 in \cite{Uniqueness}
for more details about this interpretation). We start by setting, for every
$s,t\in[0,\sigma]$, 
$$\wt{D}^\circ(s,t):= \wh W_s + \wh W_t -2\,\min_{s\wedge t\leq r\leq s\vee t} \wh W_r,$$
and then, for every $a,b\in \t_{(\omega)}$,
\begin{equation}
\label{D0-slice}
\wt{D}^\circ(a,b):=\min\{\wt{D}^\circ(s,t):s,t\in[0,\sigma], p_{(\omega)}(s)=a,\,p_{(\omega)}(t)=b\}.
\end{equation}
We finally let $\wt D(a,b)$ be the maximal pseudo-metric on $\t_{(\omega)}$
that is bounded above by $\wt{D}^\circ(a,b)$. We define the slice $\SS{(\omega)}$
as the quotient space $\t_{(\omega)}/\{\wt D=0\}$, which is equipped with the
metric induced by $\wt D$. We write $\wt\Pi$ for the canonical projection from 
$\t_{(\omega)}$ onto $\SS{(\omega)}$. We may view $\SS{(\omega)}$ as a $2$-pointed measure 
metric space, with the two distinguished points $\wt x_*:=\wt\Pi(a_*)$ and $\wt x_0:=\wt\Pi(\rho_{(\omega)})$ and
the volume measure which is the pushforward of the volume measure on $\t_{(\omega)}$. 

It is immediate to verify that $\wt{D}(a,b)\geq D(a,b)$, and therefore 
$\wt{D}(a,b)=0$ implies that $D(a,b)=0$. Conversely, suppose that
$D(a,b)=0$ and $a\not = b$. We know that $D^\circ(a,b)=0$, and thus
(up to interchanging $a$ and $b$) we can assume that
$$\Lambda_a=\Lambda_b=\min_{c\in[a,b]}\Lambda_c.$$
Pick $s,t\in[0,\sigma]$ such that $p_{(\omega)}(s)=a,\,p_{(\omega)}(t)=b$, and $[s,t]$ is as
small as possible, where we use the convention $[s,t]=[s,\sigma]\cup[0,t]$ if $s>t$. 
Then the equalities of the last display are equivalent to 
$$\wh W_s=\wh W_t= \min_{r\in[s,t]} \wh W_r.$$
If $s\leq t$, this implies $\wt D^\circ(s,t)=0$, and thus $\wt D(a,b)=0$. On the other hand, 
if $s>t$, we obtain that necessarily 
$$\wh W_s=\wh W_t,\qquad \wh W_s=\min_{s\leq r\leq \sigma}\wh W_r,\qquad\wh W_t=\min_{0\leq r\leq t}\wh W_r.$$
These equalities imply that $\Pi(a)=\Pi(b)=\gamma(-\Lambda_a)$ belongs to the range of the geodesic $\gamma$. 
On the other hand, if, for every $r\in(0,h)$, we
take $a=p_{(\omega)}(\sup\{s\in[0,\sigma]:\wh W_s=-r\})$ and $b=p_{(\omega)}(\inf\{s\in[0,\sigma]:\wh W_s=-r\})$
we have $\Pi(a)=\Pi(b)=\gamma(r)$, but $\wt D(a,b)>0$ (see Lemma 12 (ii) in \cite{BM}) and therefore
$\wt \Pi(a)\not =\wt \Pi(b)$. 

The preceding considerations show that every point of $\bm_\infty(\omega)$ that does not belong to the
geodesic $\gamma$ corresponds to a 
single point of $\SS{(\omega)}$, but every point of the geodesic $\gamma$ (other than $x_0$ and $x_*$) corresponds to
two points of $\SS{(\omega)}$. More precisely, if we set for every $r\in [0,h]$,
$$\gamma'(r):=\wt\Pi\circ p_{(\omega)}(\inf\{s\in[0,\sigma]:\wh W_s=-r\}),\quad \gamma''(r):=
\wt\Pi\circ p_{(\omega)}(\sup\{s\in[0,\sigma]:\wh W_s=-r\}),$$
then $\gamma'$ and $\gamma''$ are two geodesics in $\mathbf{S}{(\omega)}$ from $\wt x_0$
to $\wt x_*$ that are disjoint except at their initial and terminal times. We  call 
$\gamma'$ and $\gamma''$ the left and right boundary geodesics of $\SS{(\omega)}$.

We will use the fact that, for every $\kappa\in(0,h/2)$, one has a.s.
\begin{equation}
\label{dist-geo}
\inf_{s,t\in[\kappa,h-\kappa]} \wt D(\gamma'(s),\gamma''(t)) >0.
\end{equation}
This is immediate since the function $(s,t)\mapsto \wt D(\gamma'(s),\gamma''(t))$
is continuous and does not vanish on $[\kappa,1-\kappa]\times[\kappa,1-\kappa]$.

Taking $h=1$, \eqref{dist-geo} allows us to find a sequence $(\delta_k)_{k\geq 1}$ of positive reals such that the probability (under $\N_0(\cdot\midd W_*=-1)$)
of the event where 
$$\inf_{s,t\in[1-2^{-k},1-2^{-k-4}]} \wt D(\gamma'(s),\gamma''(t)) > 2\delta_k\,,\quad \hbox{for every }k\geq 1,$$
is at least $9/10$. Without loss of generality, we may and will assume that $\delta_k<2^{-k-5}$ for every $k\geq 1$. By scaling, we also obtain that, for every $h\in [3/4,1]$, the probability
under $\N_0(\cdot\midd W_*=-h)$ of the
event where 
$$\inf_{s,t\in[h(1-2^{-k}),h(1-2^{-k-4})]} \wt D(\gamma'(s),\gamma''(t)) > \delta_k\,,\quad \hbox{for every }k\geq 1,$$
is at least $9/10$. Finally, fix $\ve\in(0,1/4)$, and observe that, if $h\in[1-\ve,1]$ and if the integer $k\geq 1$ is such that $2^{-k-4}>\ve$, we
have $[1-2^{-k},1-2^{-k-3}]\subset [h(1-2^{-k}),h(1-2^{-k-4})]$. We arrive at the following lemma.

\begin{lemma}
\label{preli-bd-dist}
Let $\ve\in(0,1/4)$.
Then, for every $h\in[1-\ve,1]$, the probability
under $\N_0(\cdot\midd W_*=-h)$ of the event where
\begin{equation}
\label{low-bd-dist}
\inf_{s,t\in[1-2^{-k},1-2^{-k-3}]} \wt D(\gamma'(s),\gamma''(t)) > \delta_k\,,\quad \hbox{for every }k\geq 1\hbox{ such that
}2^{-k-4}>\ve,
\end{equation}
is at least $9/10$. 
\end{lemma}

\subsection{Slices in hulls}
\label{slice-hull}

We consider the standard hull $\mathfrak{H}$ of radius $r=1$ and perimeter $z$ as constructed in Section \ref{cons-hull}
from $(U_*,\omega_*)$ and the point measure $\sum_{i\in I}\delta_{(t_i,\omega_i)}$, and
we keep the notation of this section. To avoid confusion, we write $[a,b]_\mathbf{H}$ for the
intervals of $\mathbf{H}$ as defined in Section \ref{cons-hull}. We fix $m\in\{1,2,3\}$ and $\ve\in(0,1/4)$, and we will now condition on the event 
\begin{equation}
\label{m-slice}
E^m_\ve:=\{\#\{i\in I:W_*(\omega_i)< -1+\ve\} \geq m\}.
\end{equation}
Under this conditioning, there are (at least) $m$ indices $i\in I$
such that $W_*(\omega_i)\leq -1+\ve$, and we write $i_1,\ldots,i_m$ for these indices ranked in
such a way that $t_{i_1}<\cdots<t_{i_m}$ (if there are more than $m$ indices with the desired property,
we keep those corresponding to the smaller values of $t_i$). We then set
$$\RR_j:=\Pi_\mathfrak{H}(\t_{(\omega_{i_j})})\,,\ \hbox{for }1\leq j\leq m\,,\qquad \RR_*=\Pi_\mathfrak{H}(\t_{(\omega_*)}).$$
Let $r_1\in[1-\ve,1]$. We note that, conditionally on $W_*(\omega_{i_1})=-r_1$, $\omega_{i_1}$
is distributed according according to $\N_0(\cdot\midd W_*=-r_1)$, and so we can consider 
the slice $\SS{(\omega_{i_1})}$ constructed in the previous section (with $h=r_1$).

We observe that $\RR_1$ and $\SS{(\omega_{i_1})}$ are canonically identified as sets. Indeed both
$\RR_1$ and $\SS{(\omega_{i_1})}$ are quotient spaces of $\t_{(\omega_{i_1})}$, and the point is to verify that,
for $a,b\in \t_{(\omega_{i_1})}$, we have $D^\circ_\HH(a,b)=0$ if and only if $\wt D^\circ(a,b)=0$
(we abuse notation by still writing $\wt D^\circ(a,b)$ for the function defined in \eqref{D0-slice}
when $\omega=\omega_{i_1}$). Suppose first that $\wt D^\circ(a,b)=0$. Then there
exist $s,t\in[0,\sigma(\omega_{i_1})]$ such that $p_{(\omega_{i_1})}(s)=a$, $p_{(\omega_{i_1})}(t)=b$, and
$$\wh W_s(\omega_{i_1})=\wh W_t(\omega_{i_1})= \min_{s\wedge t\leq u\leq s\vee t} \wh W_u(\omega_{i_1}).$$
Suppose that $s\leq t$ for definiteness. Then, writing $[u_1,u_1+\sigma(\omega_{i_1})]$
for the time interval corresponding to $\t_{(\omega_{i_1})}$ in the time scale of the cyclic exploration
$\ee$ of $\HH$, we have $[a,b]_\mathbf{H}\subset \{\ee_{u_1+u}:s\leq u\leq t\}$, so that we get
\begin{equation}
\label{tech-slice-hull}
\Lambda_a=\Lambda_b=\min_{c\in[a,b]_\mathbf{H}} \Lambda_c
\end{equation}
and therefore $D^\circ_\HH(a,b)=0$. Conversely, if $D^\circ_\HH(a,b)=0$, we can assume without
loss of generality that \eqref{tech-slice-hull} holds, and this is only possible if $[a,b]_\mathbf{H}\subset \t_{(\omega_{i_1})}$
(otherwise we would have $b_*\in[a,b]_\mathbf{H}$). But then, using the definition of intervals in $\HH$, we have
$$\min_{c\in[a,b]_\mathbf{H}} \Lambda_c= \max\Big\{\min_{u\in[s,t]}\wh W_u: s\leq t\leq \sigma(\omega_{i_1}),\,p_{(\omega_{i_1})}(s)=a,
\,p_{(\omega_{i_1})}(t)=b\Big\},$$
and \eqref{tech-slice-hull} implies that $\wt D^\circ(a,b)=0$. Similarly, $\RR_j$ and $\SS{(\omega_{i_j})}$ are canonically identified as sets, for every $1\leq j\leq m$ (note
however that $\SS{(\omega_*)}$, which also makes sense since $\omega_*$
is distributed according to $\N_0(\cdot\midd W_*=-1)$, is {\it not} identified to $\RR_*$). 

We also observe that we have, for every $a,b\in\t_{(\omega_{i_1})}$,
\begin{equation}
\label{bound-dista}
\wt D(a,b)\geq D_\HH(a,b).
\end{equation}
Indeed, it is immediately seen that $\wt D(a,b)$ is obtained by considering the same infimum as in
\eqref{dist-hull}, with the additional restriction that we require $a_1,\ldots,a_{k-1}\in\t_{(\omega_{i_1})}$.

Since both metric spaces $(\RR_1,D_\HH)$ and $(\SS{(\omega_{i_1})},\wt D)$ are compact, the bound \eqref{bound-dista} implies that the
identity mapping from $\SS{(\omega_{i_1})}$ (equipped with $\wt D$) onto $\RR_1$ (equipped with $D_\HH$)
is a homeomorphism. See Fig.~1 below for a schematic representation of $\mathbf{R}_1,\mathbf{R}_2,\mathbf{R}_*$
in the case where $E^2_\ve$ holds. 

Consider then the two geodesics $\gamma'$ and $\gamma''$ from $\wt x_0$ to $\wt x_*$ in $\SS{(\omega_{i_1})}$,
as defined in the previous section. Let $\Gamma':=\{\gamma'(t):0\leq t\leq -W_*(\omega_{i_1})\}$
denote the range of $\gamma'$, and similarly let  $\Gamma''$ denote the range of $\gamma''$. 
Via the identification of $\RR_1$ and $\SS{(\omega_{i_1})}$, we may and will view $\Gamma'$ and $\Gamma''$
as subsets of $\RR_1$. Then, it is easily checked that $\Gamma'\cup \Gamma''$ is the topological boundary
of $\RR_1$ in $\HH$ (note that a point $a$ of $\t_{(\omega_{i_1})}$ different from the root is identified to a point of $\HH\backslash \t_{(\omega_{i_1})}$
if and only if $\Pi_\mathfrak{H}(a)$ belongs to $\Gamma'\cup \Gamma''$). Furthermore, the restriction of $D_\HH$ to $\Gamma'$ (resp. to $\Gamma''$)
clearly coincides with the restriction of $\wt D$. We write $\mathrm{Int}(\RR_1)=\RR_1\backslash (\Gamma'\cup \Gamma'')$
for the interior of $\RR_1$. We call $\Gamma'$ and $\Gamma''$ the left and right boundaries of the slice $\RR_1$
and write $\Gamma'=\partial_\ell\mathbf{R}_1$ and $\Gamma''=\partial_r\mathbf{R}_1$. We use
a similar terminology for the slices $\mathbf{R}_j$, $2\leq j\leq m$. 

\begin{lemma}
\label{length-path}
Let $(\phi(t))_{0\leq t\leq 1}$ be a continuous path in $\RR_1$ satisfying the condition:
\begin{description}
\item[\rm(H)] There exist two reals $u$ and $v$ such that $0\leq u\leq v\leq 1$ such that $\phi(t)\in\Gamma'$ for every $t\in[0,u]$,
$\phi(t)\in\Gamma''$ for every $t\in[v,1]$, and $\phi(t)\in\mathrm{Int}(\RR_1)$ for every $t\in(u,v)$. 
\end{description}
Then, the length 
of $\phi$ with respect to the distance $D_\HH$ coincides with its length with respect to the distance $\wt D$.
\end{lemma}

\proof Since the restriction of $D_\HH$ to $\Gamma'$ or to $\Gamma''$
coincides with the restriction of $\wt D$, we may assume that condition (H) holds with $u=0$ and $v=1$. By a continuity argument, we 
may even replace (H) by 
\begin{description}
\item[\rm(H')] $\phi(t)\in\mathrm{Int}(\RR_1)$ for every $t\in[0,1]$. 
\end{description}
Then, we can find $\delta>0$ such that 
$D_\HH(\phi(t),(\mathrm{Int}(\RR_1))^c)\geq \delta$ for every $t\in [0,1]$. The conclusion of
the lemma follows from the fact that we have $\wt D(\phi(t),\phi(t'))=D_\HH(\phi(t),\phi(t'))$
as soon as $t,t'\in[0,1]$ are such that $D_\HH(\phi(t),\phi(t'))\leq \delta/2$. Indeed, 
in the definition \eqref{dist-hull} of $D_\HH(\phi(t),\phi(t'))$, we may restrict the infimum
of the quantities $\sum D_\HH^\circ(a_{i-1},a_i)$ to the case where $a_1,\ldots,a_{k-1}$ all belong
to $\mathrm{Int}(\RR_1)$, because if for instance $a_j\notin \mathrm{Int}(\RR_1)$, we have
$$\sum_{i=1}^j D^\circ_\HH(a_{i-1},a_i)\geq D_\HH(\phi(t),a_j)\geq \delta.$$
If $a_1,\ldots,a_{k-1}$ all belong to $\mathrm{Int}(\RR_1)$, we have 
$\sum D_\HH^\circ(a_{i-1},a_i)=\sum \wt D^\circ(a_{i-1},a_i)$, and thus we
obtain that $\wt D(\phi(t),\phi(t'))\leq D_\HH(\phi(t),\phi(t'))$. The reverse bound follows from \eqref{bound-dista}. \endproof

We can now combine Lemma \ref{length-path}
with the discussion of the end of Section \ref{sec:slice}. 
For every $0\leq u\leq v\leq -W_*(\omega_{i_1})$, we use the notation
$$\Gamma'_{[u,v]}=\{\gamma'(t):u\leq t\leq v\},$$
and we similarly define $\Gamma''_{[u,v]}$.
We let the sequence $(\delta_k)_{k\geq 1}$ be as in Lemma \ref{preli-bd-dist}, and we recall the 
definition \eqref{m-slice} of the event $E^m_\ve$. 

\begin{lemma}
\label{key-bd-dist}
Let $\ve\in(0,1/4)$. The following property holds with probability at least $9/10$
under $\P(\cdot\midd E^m_\ve)$. For every integer $k\geq 1$
such that $2^{-k-4}>\ve$, for every continuous path 
$(\phi(t))_{0\leq t\leq 1}$ in $\RR_1$ such that $\phi(0)\in \Gamma'_{[1-2^{-k},1-2^{-k-3}]}$
and $\phi(1)\in \Gamma''_{[1-2^{-k},1-2^{-k-3}]}$, the length (with respect to $D_\HH$)
of $\phi$ is at least $\delta_k$. 
\end{lemma}

\proof First note that we may restrict our attention to paths $\phi$ satisfying condition (H) 
of Lemma \ref{length-path}: Indeed, we may set $u=\sup\{t\in[0,1]:\phi(t)\in\Gamma'\}$ and
$v=\inf\{t\in[u,1]:\phi(t)\in\Gamma''\}$, and observe that replacing 
$\phi$ by a portion of the geodesic $\gamma'$ (resp. of the geodesic $\gamma''$) on the interval $[0,u]$
(resp. on $[v,1]$) will only decrease its length. By Lemma \ref{length-path}, if $\phi$ satisfies (H), its length
with respect to $D_\HH$ is the same as its length with respect to $\wt D$, provided we view $\phi$
as a path in $\SS{(\omega_{i_1})}$. In particular, this length is bounded below by $\wt D(\phi(0),\phi(1))$.
The desired result now follows from Lemma \ref{preli-bd-dist}. \endproof

We also need an analog of Lemma \ref{key-bd-dist} when $\RR_1$ is replaced by $\RR_*$. The situation in that case 
is a bit different since $\SS{(\omega_*)}$ is no longer identified bijectively with $\RR_*$. Instead, $\RR_*$ appears as 
a quotient space of $\SS{(\omega_*)}$, where, for 
$u\in(0,1]$, the points $\gamma'(u)$ and $\gamma''(u)$ of the left and right boundary geodesics 
in $\SS{(\omega*)}$ are identified in $\RR_*$ if (and only if) $u\geq \mu_0$, where
\begin{equation}
\label{coales-time}
\mu_0:=\sup\{-W_*(\omega_i):i\in I\}\geq 1-\ve.
\end{equation}
We can nonetheless define the range $\Gamma'$ (resp. $\Gamma''$) of the left boundary geodesic 
(resp. of the right boundary geodesic) as a closed subset of $\RR_*$, and make sense of the sets
$\Gamma'_{[u,v]}$ and $\Gamma''_{[u,v]}$ for $0\leq u\leq v\leq 1$. We have $\Gamma'_{[u,v]}\cap \Gamma''_{[u,v]}=\varnothing$
if $v<\mu_0$. 
The reader will easily check that an analog of Lemma \ref{length-path} remains valid
when $\RR_1$ is replaced by $\RR_*$, provided that we add the constraint $D_\HH(\phi(t),b_*)> \ve$
for every $t\in[0,1]$ in condition (H). Then, the proof of Lemma \ref{key-bd-dist} when $\RR_1$ is
replaced by $\RR_*$ goes through almost without change: recall that we have assumed $\delta_k<2^{-k-5}$, and therefore a 
path $(\phi(t))_{0\leq t\leq 1}$ in $\RR_*$ such that $\phi(0)\in \Gamma'_{[1-2^{-k},1-2^{-k-3}]}$
will not visit the set $\{x\in \RR_*:D(x,x_*)\leq\ve\}$
if its length is bounded by $\delta_k
$. The preceding discussion is summarized in the following statement.

\begin{lemma}
\label{key-bd-dist2}
The statement of Lemma \ref{key-bd-dist} remains valid if $\RR_1$ is replaced by $\RR_*$.
\end{lemma}

\subsection{The first-moment estimate}
\label{sec:first-mom}

We again fix $m\in\{1,2,3\}$
($m+1$ will correspond to what we called $m$ in the introduction). For $x,y\in\bm_\infty$, we recall our notation $B^{\bullet(y)}_r(x)$ for the hull of radius $r$ centered at $x$ relative to
$y$, which makes sense if $D(x,y)>r$. 

For $x,y\in\bm_\infty$ and $\ve\in(0,r)$, we let $F^{(m)}_{\ve,r}(x,y)$ 
be equal to $1$ if $D(x,y)>r$ and there exist $m+1$ geodesic paths $(\xi_0(t))_{0\leq t\leq r},
(\xi_1(t))_{0\leq t\leq r},\ldots,(\xi_m(t))_{0\leq t\leq r}$
with $\xi_i(0)\in \partial B^{\bullet(y)}_r(x)$ and $\xi_i(r)=x$, for every $i\in\{0,\ldots,m\}$, 
such that the sets $\{\xi_i(t):0\leq t\leq r-\ve\}$, for $i\in\{0,\ldots,m\}$,  are disjoint. If these 
conditions do not hold, we take $F^{(m)}_{\ve,r}(x,y)=0$.
By convention, if $\ve\geq r$, we take  $F^{(m)}_{\ve,r}(x,y)=1$ if $D(x,y)>r$ and $F_{\ve,r}^{(m)}(x,y)=0$ otherwise.

In the case $r=1$, we also define $\wt F^{(m)}_{\ve,1}(x,y)$ exactly as $F^{(m)}_{\ve,1}(x,y)$,
except that we require the additional property that the geodesics $\xi_0,\xi_1,\ldots,\xi_m$ satisfy
\begin{equation}
\label{condit-dist}
D(\xi_i(t),\xi_j(t))\geq \delta_k
\end{equation}
for every $t\in[1-2^{-k-1},1-2^{-k-2}]$, for every $k\geq 1$ such that $2^{-k-4}>\ve$, and for
every $0\leq i<j\leq m$. 
Here the sequence $(\delta_k)_{k\geq 1}$ is fixed as in Lemma \ref{key-bd-dist}. 

We let $\wt\N_0$ denote the measure with density $\frac{1}{\sigma}$ with respect to $\N_0$.

\begin{proposition}
\label{prop-first-mom}
There exists a constant $c>0$ such that, for every $\ve\in(0,1/4)$, we have
$$\wt\N_0\Big(\int \mathrm{Vol}(\dd x)\,\mathbf{1}_{\{D(x,x_*)<2\}}\,\wt F^{(m)}_{\ve,1}(x,x_*)\Big) \geq c\,\ve^m.$$
\end{proposition}

\proof By the symmetry properties of the Brownian sphere (Section \ref{symm}), we have
\begin{align*}
\wt\N_0\Big(\int \mathrm{Vol}(\dd x)\,\,\mathbf{1}_{\{D(x,x_*)< 2\}}\,\wt F^{(m)}_{\ve,1}(x,x_*)\Big)&=\N_0\Big(
\,\mathbf{1}_{\{D(x_*,x_0)<2\}}\,\wt F^{(m)}_{\ve,1}(x_*,x_0)\Big)\\
&= \frac{3}{2}\,\N^{[1]}_0\Big(\,\mathbf{1}_{\{D(x_*,x_0)< 2\}}\,\wt F^{(m)}_{\ve,1}(x_*,x_0)\Big).
\end{align*}
We then use Theorem \ref{deco-hull}, noting that the indicator function $\mathbf{1}_{\{D(x_*,x_0)< 2\}}$ is 
a function of $\ov{\bm_\infty \backslash B^{\bullet(x_0)}_1(x_*)}$ and that $\wt F^{(m)}_{\ve,1}(x_*,x_0)$ 
can be written as a function of the hull $B^{\bullet(x_0)}_1(x_*)$: for this last fact, observe that
geodesics from the boundary of $B^{\bullet(x_0)}_1(x_*)$
to the center $x_*$ are the same for $D$ and for the intrinsic distance of the hull, and that in condition \eqref{condit-dist}
we can also replace $D$ by the intrinsic distance, since if $D(\xi_i(t),\xi_j(t))$ is smaller than the intrinsic distance
between $\xi_i(t)$ and $\xi_j(t)$, this means that a geodesic (for $D$) from $\xi_i(t)$ to $\xi_j(t)$ 
has to intersect the boundary of the hull, and thus $D(\xi_i(t),\xi_j(t))$ is at least $1/2$. It follows from Theorem \ref{deco-hull}
and
these observations that
$$\N^{[1]}_0\Big(\,\mathbf{1}_{\{D(x_*,x_0)< 2\}}\,\wt F^{(m)}_{\ve,1}(x_*,x_0)\Big)
=\N^{[1]}_0\Big(\theta(\z_{W_*+1})\,\N^{[1]}_0(\wt F^{(m)}_{\ve,1}(x_*,x_0)\midd \z_{W_*+1})),$$
where $\theta(z)$ is the probability for a free pointed Brownian disk of perimeter $z$ that
the distance from the distinguished point to the boundary is smaller than $1$. 
Since we also know that, under $\N^{[1]}_0(\cdot \midd \z_{W_*+1}=z)$, $B^{\bullet(x_0)}_1(x_*)$ is distributed as a 
standard hull of radius $1$ and perimeter $z$ (as defined in Section \ref{cons-hull}), the proof of Proposition \ref{prop-first-mom} reduces to
establishing the following claim.

\medskip

\noindent{\it Claim.} Let $\mathfrak{H}$ be a standard hull of radius $1$ and perimeter $z$ as in Section \ref{cons-hull}, and, for every $\ve\in(0,1/4)$, let
$\mathcal{A}^m_\ve$ be the event where there exist $m+1$ geodesics $(\eta_0(t))_{0\leq t\leq 1}$,
$(\eta_1(t))_{0\leq t\leq 1},\ldots,(\eta_m(t))_{0\leq t\leq 1}$ from the boundary $\partial\mathfrak{H}$
to the center of the hull, such that the sets $\{\eta_j(t):0\leq t\leq 1-\ve\}$ are disjoint, and moreover
\begin{equation}
\label{condit-dist2}
D_\mathbf{H}(\eta_i(t),\eta_j(t))\geq \delta_k
\end{equation}
for every $t\in[1-2^{-k-1},1-2^{-k-2}]$, for every $k\geq 1$ such that $2^{-k-4}>\ve$, and for
every $0\leq i<j\leq m$. Then there exists a constant $c>0$ (depending on $z$) such that $\P(\mathcal{A}^m_\ve)\geq c\,\ve^m$.

\medskip
In the remaining part of the proof, we establish the preceding claim.
We consider the event $E^m_\ve$ defined in \eqref{m-slice}. Then $\P(E^m_\ve)$ is just the probability that
a Poisson variable with parameter $\frac{3z}{2}((1-\ve)^{-2}-1)$ is greater than or equal to $m$, and
thus there exists a constant $c'$ such that $\P(E^m_\ve)\geq c'\ve^m$. To get the desired lower
bound, it remains to verify that we have also $\P(\mathcal{A}^m_\ve\midd E^m_\ve)\geq c''$
with another constant $c''>0$. To this end, we rely on Lemma \ref{key-bd-dist} and Lemma \ref{key-bd-dist2}.
For every $1\leq j\leq m$ (resp. for $j=0$), we write $\mathcal{B}^{m,j}_\ve$
for the intersection of $E^m_\ve$ with the event where the property of Lemma \ref{key-bd-dist} holds when $\mathbf{R}_1$ is replaced by $\mathbf{R}_j$
(resp. by $\mathbf{R}_*$). We then know that
$$\P((\mathcal{B}^{m,j}_\ve)^c\midd E^m_\ve)\leq \frac{1}{10},\quad\hbox{for every }0\leq j\leq m,$$
and, if
$$\mathcal{B}^m_\ve=\bigcap_{j=0}^m \mathcal{B}^{m,j}_\ve,$$
 it follows that $\P(\mathcal{B}^m_\ve\midd E^m_\ve)\geq 1/2$. So to get our claim we only
 need to verify that $\mathcal{B}^m_\ve\subset \mathcal{A}^m_\ve$. 
 
 From now on, we argue on the event $E^m_\ve$, and we let $i_1,\ldots,i_m$ be the indices such that
 $W_*(\omega_i)< -1+\ve$, as defined at the beginning of Section \ref{slice-hull}. Then $t_{i_1},\ldots,t_{i_m}$
 are elements of $[0,z]\subset \mathbf{H}$, and, recalling the definition of the exploration process $(\ee_s)_{s\in[0,\Sigma]}$,
 we set
 $$\begin{array}{rll}
 s'_j\!\!\!&:=\inf\{s\in[0,\Sigma]: \ee_s=t_{i_j}\}\;,\quad s''_j\!\!\!&:=\sup\{s\in[0,\Sigma]: \ee_s=t_{i_j}\},\quad\hbox{for }1\leq j\leq m,\\
 \noalign{\smallskip}
 s'_0\!\!\!&:=\inf\{s\in[0,\Sigma]: \ee_s=U_*\}\;,\quad s''_0\!\!\!&:=\sup\{s\in[0,\Sigma]: \ee_s=U_*\}.
 \end{array}
 $$
 We can then consider the simple geodesics $(\gamma_{s'_j}(t))_{0\leq t\leq 1}$ and $(\gamma_{s''_j}(t))_{0\leq t\leq 1}$, for every $0\leq j\leq m$, as defined
 in Section \ref{cons-hull}. We note that $\gamma_{s'_j}(0)=\gamma_{s''_j}(0)$ and $\gamma_{s'_j}(r)\not=\gamma_{s''_j}(r)$ if
 $0< r < \mu_j$, where $\mu_0$ was defined in \eqref{coales-time} and $\mu_j:=-W_*(\omega_j)\in(1-\ve,1)$ if $1\leq j\leq m$. On the other hand,
 $\gamma_{s'_j}(r)=\gamma_{s''_j}(r)$ if $\mu_j\leq r\leq 1$. The set
 $$\{\gamma_{s'_j}(r):0\leq r\leq \mu_j\}\cup \{\gamma_{s''_j}(r):0\leq r\leq \mu_j\}$$
 is the range of a simple cycle, and (the closure of) the connected component of the complement of this simple cycle
 that does not contain the point $\Pi_\mathfrak{H}(0)$ of $\partial\mathfrak{H}$ coincides with the slice
 $\mathbf{R}_j=\Pi_\mathfrak{H}(\t_{(\omega_{i_j})})$ (resp. with $\mathbf{R}_*=\Pi_\mathfrak{H}(\t_{(\omega_*)})$ when $j=0$).
 More precisely, for $1\leq j\leq m$, the set $\{\gamma_{s'_j}(r):0\leq r\leq \mu_j\}$ is the left boundary of $\mathbf{R}_j$ , and the set $\{\gamma_{s''_j}(r):0\leq r\leq \mu_j\}$ is the right boundary of $\mathbf{R}_j$. Similarly, we can interpret $\{\gamma_{s'_0}(r):0\leq r\leq \mu_0\}$ as the left boundary of $\mathbf{R}_*$,
 and $\{\gamma_{s''_0}(r):0\leq r\leq \mu_0\}$ as the right boundary of $\mathbf{R}_*$.
 Since $\gamma_{s'_j}$ and $\gamma'_{s'_j}$ are geodesics to $b_*$, we have
$$D_{\mathbf{H}}(b_*,\gamma_{s'_j}(r))=1-r=D_{\mathbf{H}}(b_*,\gamma_{s''_j}(r))$$
 for every $r\in[0,1]$. 
 
 We now verify that, if $\mathcal{B}^{m}_\ve$ holds, then $\mathcal{A}^m_\ve$ also holds, and we can take
 $\eta_j=\gamma_{s'_j}$ for $0\leq j\leq m$. We first observe that, thanks to the fact that $W_*(\omega_{i_j})< -1+\ve$
 for $1\leq j\leq m$, the sets $\{\gamma_{s'_j}(t):0\leq t\leq 1-\ve\}$, $0\leq j\leq 1-\ve$, are disjoint. Then, 
 we use the defining property of $\mathcal{B}^{m}_\ve$ to get that, for every integer $k\geq 1$ such that $2^{-k-4}>\ve$, 
 the length of any continuous path starting  from $\{\gamma_{s'_j}(t): 1-2^{-k}\leq t\leq 1-2^{-k-3}\}$,
ending in $\{\gamma_{s''_j}(t): 1-2^{-k}\leq t\leq 1-2^{-k-3}\}$ and staying in the slice $\mathbf{R}_j$ (in $\mathbf{R}_*$
if $j=0$) is bounded below by $\delta_k$. Let us explain why this implies that
\begin{equation}
\label{first-mom1}
D_{\mathbf{H}}(\gamma_{s'_i}(t),\gamma_{s'_j}(t))\geq \delta_k
\end{equation}
for every $t\in[1-2^{-k-1},1-2^{-k-2}]$, for every $k\geq 1$ such that $2^{-k-4}>\ve$, and for
every $0\leq i<j\leq m$. 

\begin{figure}[!h]
\label{hull-picture}
 \begin{center}
 \includegraphics[width=10cm]{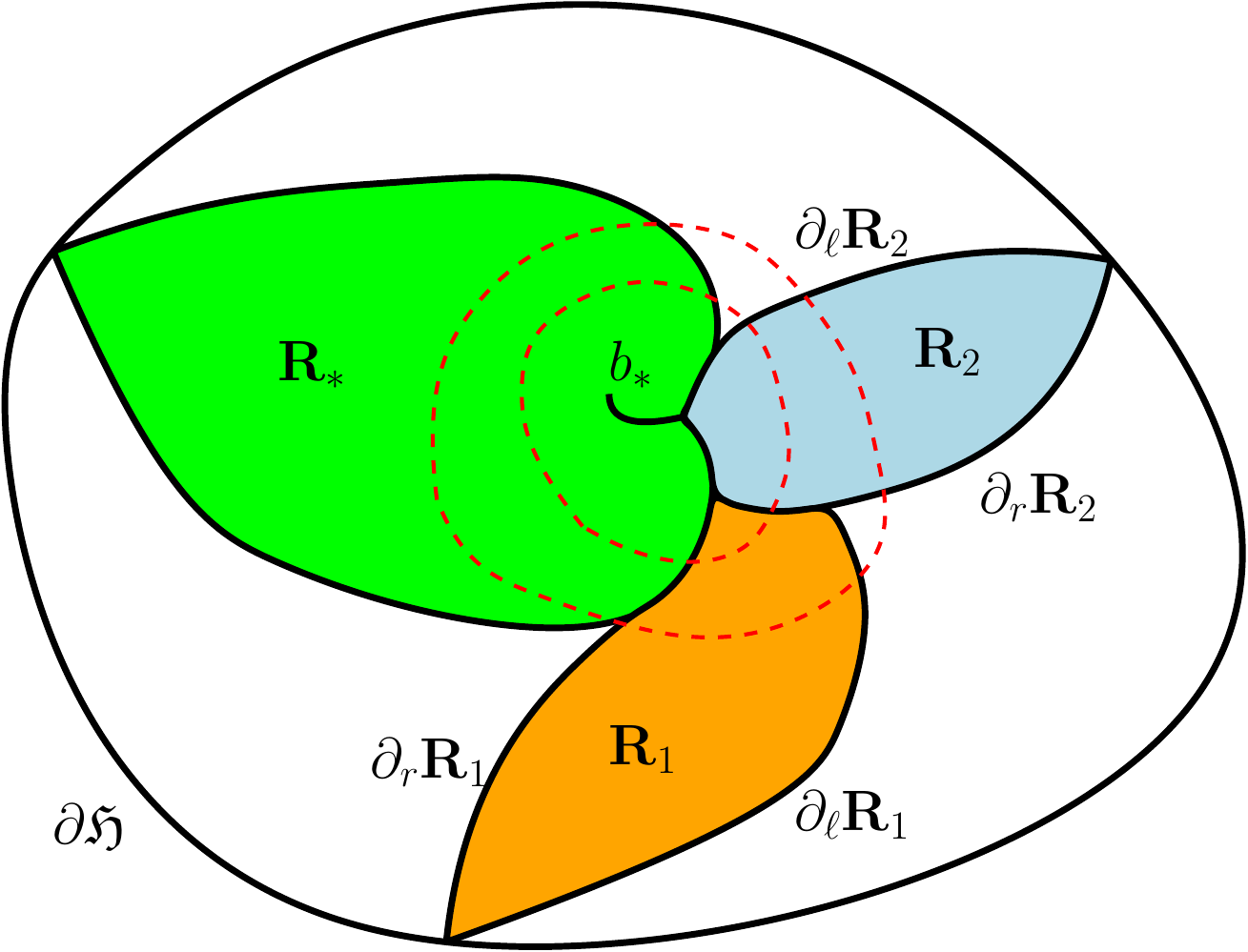}
 \caption{A schematic representation of the standard hull $\mathfrak{H}$ and the slices
 $\mathbf{R}_*$, $\mathbf{R}_1$, $\mathbf{R}_2$ in the case when $E^m_\ve$ holds with $m=2$. 
 The two dashed cycles are the boundaries of the hulls of respective radii $1-2^{-k}$ and $1-2^{-k-3}$ 
 centered at $b_*$. The end of the proof of Proposition \ref{prop-first-mom} relies on the fact that a continuous 
 path starting on the left boundary $\partial_\ell\mathbf{R}_1$ and ending on $\partial_\ell\mathbf{R}_2$
 and staying in the annulus delimited by the two dashed cycles
 will have to cross one of the two slices $\mathbf{R}_1$ and $\mathbf{R}_2$.}
 \end{center}
 \vspace{-5mm}
 \end{figure}

For simplicity we take $i=1$ and $j=2$ but the same argument works as well for any choice of 
$i$ and $j$. Fix $k\geq 1$ such that $2^{-k-4}>\ve$ and $t_0\in[1-2^{-k-1},1-2^{-k-2}]$. 
Consider a continuous path $(\phi(r))_{0\leq r\leq 1}$
such that $\phi(0)=\gamma_{s'_1}(t_0)\in\partial_\ell\mathbf{R}_1$ and $\phi(1)=\gamma_{s'_2}(t_0)
\in\partial_\ell\mathbf{R}_2$. We need
to show that the length of $\phi$ is bounded below by $\delta_k$. We may and will assume
that $D_\mathbf{H}(b_*,\phi(r))\in [2^{-k-3},2^{-k}]$ for every $r\in[0,1]$, because
otherwise the length of $\phi$ will be bounded below by $2^{-k-3}\geq \delta_k$. 

We then set 
$$
r_0:=\sup\{r\in[0,1]: \phi(r)\in \partial_\ell\mathbf{R}_1\}\,,\quad
r_1:=\inf\{r\in [r_0,1]:\phi(r)\in \partial_r\mathbf{R}_1\cup\partial_r\mathbf{R}_2\}.
$$
We note that the set $\{r\in[0,1]: \phi(r)\in \partial_r\mathbf{R}_1\cup\partial_r\mathbf{R}_2\}$ is not empty, so that
the definition of $r_1$ makes sense. Indeed, from the construction of $\mathfrak{H}$ (or via a planarity argument), it is easy to verify that
a path $(\psi(r))_{r\in[0,1]}$ starting on $\partial_\ell \mathbf{R}_1$ and ending on $\partial_\ell\mathbf{R}_2$, such
that $D_{\mathbf{H}}(b_*,\psi(r))\geq\ve$ for every $r\in[0,1]$,
must hit $\partial_r\mathbf{R}_1\cup\partial_r\mathbf{R}_2$ before (or at the same time) it hits $\partial_\ell\mathbf{R}_2$ --- see Fig.~1 for an illustration.

Then we consider two cases:
\begin{description}
\item[$\bullet$] Either $\phi(r_1)\in\partial_r\mathbf{R}_1$, and then the restriction of $\phi$
to $[r_0,r_1]$ stays in $\mathbf{R}_1$, and is a path to which we can apply the definition of the event $\mathcal{B}^{m,1}_\ve$. It 
follows that the length of this restricted path is at least $\delta_k$.
\item[$\bullet$] Or $\phi(r_1)\in\partial_r\mathbf{R}_2$, then we set
$$r_2:=\sup\{r\in[r_1,1]: \phi(r)\in \partial_r\mathbf{R}_2\}\;,\quad r_3:=\inf\{r\in[r_2,1]:\phi(r)\in\partial_\ell \mathbf{R}_2\},$$
and the restriction of $\phi$ to $[r_2,r_3]$ is a path that stays in $\mathbf{R}_2$ 
and to which we can apply the definition of the event $\mathcal{B}^{m,2}_\ve$. Again it follows that the
length of this restricted path is bounded below by $\delta_k$.
\end{description}

In both cases, we conclude that the length of $\phi$ is at least $\delta_k$.  Hence the lower bound 
\eqref{first-mom1} holds (for $i=1$ and  $j=2$). We have thus obtained that $\mathcal{B}^m_\ve\subset \mathcal{A}^m_\ve$,
which completes the proof of the proposition. \endproof

\section{The key estimate}
\label{sec:key-est}

As previously in Corollary \ref{two-hulls}, it will be convenient to use the notation $\z^{x_*(x_0)}_r=\z_{W_*+r}$, for 
$r>0$ such that $W_*<-r$. It will also be useful to consider the boundary size
 $\z^{z(x_*)}_r$ of the hull $B^{\bullet(x_*)}_r(z)$, for every $z\in\bm_\infty$ such that $D(x_*,z)>r$. To this end, we can set
\begin{equation}
\label{exit-approx2}
\z^{z(x_*)}_r=\liminf_{\ve\to 0} \ve^{-2}\mathrm{Vol}((B^{\bullet( x_*)}_r(z))^c\cap B_{r+\ve}(z))\;,
\end{equation}
and we know from 
\eqref{approx-exit} and the symmetry properties of the
Brownian sphere that the liminf in the last display is indeed a limit for a.e.~$z$ (such 
that $D(x_*,z)>r$), with respect to the volume measure on $\bm_\infty$. 
We can similarly consider the boundary size
$\z^{x_*(z)}_r$ of the hull $B^{\bullet(z)}_r(x_*)$. 

As in the previous section, the integer $m\in\{1,2,3\}$ is fixed and we recall the definition of $F^{(m)}_{\ve,r}(x,y)$ and $\wt\N_0$ from the beginning of Section \ref{sec:first-mom}. Also recall that $F^{(m)}_{\ve,r}(x,y)=0$ if $D(x,y)\leq r$.

\begin{lemma}
\label{key-lem}
Let $\delta\in(0,1)$. There exists a constant $C_{(\delta)}$ such that, for every 
$\ve\in(0,1/8)$ and every integer $k\geq 0$ such that $2^{-k}>2\ve$, 
\begin{equation}
\label{key-tech1}
\wt\N_0\Bigg(\int\!\!\int \mathrm{Vol}(\dd x)\mathrm{Vol}(\dd y)\,\mathbf{1}_{\{D(x,y)\in[2^{-k+2},2^{-k+3}]\}}
F^{(m)}_{\ve,1}(x,x_*)\,F^{(m)}_{\ve,1}(y,x_*)\Bigg)
 \leq C_{(\delta)}\,2^{-(4-m)k+\delta k}\,\ve^{2m}.
\end{equation}
\end{lemma}

\proof To simply notation, we write $F_{\ve,r}(x,y)$ instead of $F^{(m)}_{\ve,r}(x,y)$ in the proof.
 Let $A_{\ve,k}$ denote the left-hand side of \eqref{key-tech1}.
We can write $A_{\ve,k}$ in the form
$$A_{\ve,k}= \N_0\Big(\sigma\,\int\!\!\int \frac{\mathrm{Vol}(\dd x)}{\sigma}\,\frac{\mathrm{Vol}(\dd y)}{\sigma}\,
\Gamma_{\ve,k}(x_*,x,y)\Big),$$
with an appropriate function $\Gamma_{\ve,k}$. Thanks to \eqref{inter-3}, we have also
$$A_{\ve,k}= \N_0\Big(\sigma\,\int \frac{\mathrm{Vol}(\dd z)}{\sigma}\,
\Gamma_{\ve,k}(z,x_*,x_0)\Big),$$
which leads to
$$A_{\ve,k}=\N_0\Bigg(\mathbf{1}_{\{D(x_*,x_0)\in[2^{-k+2},2^{-k+3}]\}}\,\int \mathrm{Vol}(\dd z)\,
F_{\ve,1}(x_*,z)\,F_{\ve,1}(x_0,z)\Bigg).$$
We write $A_{\ve,k}=A'_{\ve,k}+A''_{\ve,k}$, where $A'_{\ve,k}$ is obtained from the right-hand side of the last display
by restricting the integral with respect to $\mathrm{Vol}(\dd z)$ to the set
$(B^{\bullet(x_0)}_{2^{-k}}(x_*)\cup B^{\bullet(x_*)}_{2^{-k}}(x_0))^c$.

\medskip

\noindent{\bf First step.} We start by estimating $A'_{\ve,k}$. Note that the property $z\notin B^{\bullet(x_0)}_{2^{-k}}(x_*)\cup B^{\bullet(x_*)}_{2^{-k}}(x_0)$ means that
$z$ and $x_0$ are in the same connected component of $(B_{2^{-k}}(x_*))^c$, and similarly 
$z$ and $x_*$ are in the same connected component of $(B_{2^{-k}}(x_0))^c$. Consequently, under the
condition $z\notin B^{\bullet(x_0)}_{2^{-k}}(x_*)\cup B^{\bullet(x_*)}_{2^{-k}}(x_0)$, we have
\begin{equation}
\label{equa-hull}
B^{\bullet(z)}_{2^{-k}}(x_*)=B^{\bullet(x_0)}_{2^{-k}}(x_*)\,,\qquad B^{\bullet(z)}_{2^{-k}}(x_0)=B^{\bullet(x_*)}_{2^{-k}}(x_0).
\end{equation}
We next observe that
$$F_{\ve,1}(x_*,z)\leq F_{\ve,2^{-k}}(x_*,z)\times F_{2^{-k+4},1}(x_*,z),$$
where we recall our convention for $F_{2^{-k+4},1}(x_*,z)$ when $2^{-k+4}\geq 1$. We have also
$F_{\ve,1}(x_0,z)\leq F_{\ve,2^{-k}}(x_0,z)$, so that $A'_{\ve,k}$ is bounded above by
$$\N_0\Bigg(\mathbf{1}_{\{D(x_*,x_0)\in[2^{-k+2},2^{-k+3}]\}}\int_{(B^{\bullet(x_0)}_{2^{-k}}(x_*)\cup B^{\bullet(x_*)}_{2^{-k}}(x_0))^c
} \mathrm{Vol}(\dd z)\,
F_{\ve,2^{-k}}(x_*,z)\,F_{2^{-k+4},1}(x_*,z)\,F_{\ve,2^{-k}}(x_0,z)\Bigg).$$
Under the condition $z\notin B^{\bullet(x_0)}_{2^{-k}}(x_*)\cup B^{\bullet(x_*)}_{2^{-k}}(x_0)$, \eqref{equa-hull} holds,
which implies
$$\mathbf{1}_{\{D(x_0,x_*)>2^{-k}\}}\,F_{\ve,2^{-k}}(x_*,z)\leq F_{\ve,2^{-k}}(x_*,x_0)\,,\quad 
\mathbf{1}_{\{D(x_0,x_*)>2^{-k}\}}\,F_{\ve,2^{-k}}(x_0,z)\leq F_{\ve,2^{-k}}(x_0,x_*).$$
It follows that
\begin{align}
\label{key-tech2}
A'_{\ve,k} &\leq \N_0\Bigg(\mathbf{1}_{\{D(x_*,x_0)\in[2^{-k+2},2^{-k+3}]\}}\,F_{\ve,2^{-k}}(x_*,x_0)\,F_{\ve,2^{-k}}(x_0,x_*)\nonumber\\
&\qquad\times
\int_{(B^{\bullet(x_0)}_{2^{-k}}(x_*)\cup B^{\bullet(x_*)}_{2^{-k}}(x_0))^c}
\mathrm{Vol}(\dd z)\,F_{2^{-k+4},1}(x_*,z)\Bigg).
\end{align}
We now want to apply Corollary \ref{two-hulls} to the right-hand side. We observe that $F_{\ve,2^{-k}}(x_*,x_0)$ is
a function of the hull $B^{\bullet(x_0)}_{2^{-k}}(x_*)$, and $F_{\ve,2^{-k}}(x_0,x_*)$ is
the same function applied to the hull $B^{\bullet(x_*)}_{2^{-k}}(x_0)$. On the other hand, the quantity 
\begin{equation}
\label{key-tech222}
\mathbf{1}_{\{D(x_*,x_0)\in[2^{-k+2},2^{-k+3}]\}} \int_{(B^{\bullet(x_0)}_{2^{-k}}(x_*)\cup B^{\bullet(x_*)}_{2^{-k}}(x_0))^c}
\mathrm{Vol}(\dd z)\,F_{2^{-k+4},1}(x_*,z)
\end{equation}
is a function of $\overline{\mathcal{C}^{x_*,x_0}_{2^{-k}}}$, with the notation of Corollary \ref{two-hulls}. Let us explain this in the case where $2^{-k+4}<1$
(the case $2^{-k+4}\geq 1$ is easier and left to the reader). The indicator function $\mathbf{1}_{\{D(x_*,x_0)\in[2^{-k+2},2^{-k+3}]\}}$ is the indicator function of the event where the (intrinsic) distance between the two boundaries of $\overline{\mathcal{C}^{x_*,x_0}_{2^{-k}}}$ lies between $2^{-k+2}-2^{-k+1}$ and $2^{-k+3}-2^{-k+1}$. Furthermore, if $D(x_*,x_0)\in[2^{-k+2},2^{-k+3}]$
and $D(x_*,z)>1$, we have $(B^{\bullet(x_0)}_{2^{-k}}(x_*)\cup B^{\bullet(x_*)}_{2^{-k}}(x_0))\subset B^{\bullet(z)}_{2^{-k+4}}(x_*)$.
Let $\Delta_1$ stand for the first boundary of $\overline{\mathcal{C}^{x_*,x_0}_{2^{-k}}}$
(that is, $\Delta_1=\partial B^{\bullet(x_0)}_{2^{-k}}(x_*)$) and, using our notation 
$D^{\mathrm{intr},\cc}$ for the (extended) intrinsic distance on $\ov{\mathcal{C}^{x_*,x_0}_{2^{-k}}}$, define the hull 
of radius $r$ centered at $\Delta_1$ relative to $z$ as the complement of the connected component 
of $\{x\in \ov{\mathcal{C}^{x_*,x_0}_{2^{-k}}}:D^{\mathrm{intr},\cc}(x,\Delta_1)>r\}$ that contains $z$ (this makes sense
if $D^{\mathrm{intr},\cc}(\Delta_1,z)>r$). It follows
from the preceding considerations that the integral with respect 
to $\mathrm{Vol}(\dd z)$ in \eqref{key-tech222} can be rewritten as
$$\int_{\overline{\mathcal{C}^{x_*,x_0}_{2^{-k}}}} \mathrm{Vol}(\dd z)\,G_k(z),$$
where $G_k(z)\in\{0,1\}$ and $G_k(z)=1$ if and only if the (intrinsic) distance between $z$ and $\Delta_1$ is greater than $1-2^{-k}$, and if there are $m+1$ 
disjoint paths of (intrinsic) length $1-2^{-k+4}$ between the boundary of the hull of radius $1-2^{-k}$ centered at $\Delta_1$ relative to $z$
and the boundary of the same hull of radius $2^{-k+4}-2^{-k}$. 

Thanks to the previous observations, we can apply Corollary \ref{two-hulls} to the right-hand side of 
\eqref{key-tech2}, and we get that
\begin{align}
\label{key-tech3}
A'_{\ve,k}
&\leq \N_0\Bigg(\mathbf{1}_{\{D(x_*,x_0)\in[2^{-k+2},2^{-k+3}]\}}\,
\Theta_{2^{-k},\z^{x_*(x_0)}_{2^{-k}}}(H_\ve)\,\Theta_{2^{-k},\z^{x_0(x_*)}_{2^{-k}}}(H_\ve)\nonumber\\
&\qquad\times\int_{(B^{\bullet(x_0)}_{2^{-k}}(x_*)\cup B^{\bullet(x_*)}_{2^{-k}}(x_0))^c}
\mathrm{Vol}(\dd z)\,F_{2^{-k+4},1}(x_*,z)\Bigg),
\end{align}
where $H_\ve$ (applied to a space belonging to $\M^{\bullet b}$) denotes the event where there are $m+1$ geodesic paths from the boundary 
to the distinguished point, that are disjoint up to the time when there are at distance $\ve$ from the distinguished point.
We note that, thanks to Proposition \ref{distinct-geo}, we have
\begin{equation}
\label{bound1-geo}
\Theta_{2^{-k},z}(H_\ve)\leq C\,(\ve z 2^{3k})^m\wedge 1,
\end{equation}
with a universal constant $C$. To simplify notation, we write $\varphi_{\ve,k}(z)= \Theta_{2^{-k},z}(H_\ve)$.

Thanks to the symmetry properties of the Brownian sphere \eqref{inter-3}, we can interchange the roles of $z$
and $x_0$ in \eqref{key-tech3}, and we arrive at
\begin{align}
\label{key-tech4}
A'_{\ve,k}
&\leq \N_0\Bigg(F_{2^{-k+4},1}(x_*,x_0)\nonumber\\
&\qquad\times\int\mathrm{Vol}(\dd z)\,\mathbf{1}_{\{D(x_*,z)\in[2^{-k+2},2^{-k+3}]\}}\,
\varphi_{\ve,k}(\z^{x_*(z)}_{2^{-k}})\,\varphi_{\ve,k}(\z^{z(x_*)}_{2^{-k}})
\mathbf{1}_{\{x_0\notin B^{\bullet(z)}_{2^{-k}}(x_*)\cup B^{\bullet(x_*)}_{2^{-k}}(z)\}}\Bigg).
\end{align}

Under the condition $x_0\notin B^{\bullet(z)}_{2^{-k}}(x_*)\cup B^{\bullet(x_*)}_{2^{-k}}(z)$, we have 
$B^{\bullet(x_*)}_{2^{-k}}(z)=B^{\bullet(x_0)}_{2^{-k}}(z)$ and $B^{\bullet(z)}_{2^{-k}}(x_*)=B^{\bullet(x_0)}_{2^{-k}}(x_*)$.
It follows in particular that $\z^{x_*(z)}_{2^{-k}}=\z^{x_*(x_0)}_{2^{-k}}$. Hence, the right-hand side of
\eqref{key-tech4} is equal to
\begin{equation}
\label{key-tech44}
\N_0\Bigg(F_{2^{-k+4},1}(x_*,x_0)\,\varphi_{\ve,k}(\z^{x_*(x_0)}_{2^{-k}})\!\int\!\mathrm{Vol}(\dd z)\,\mathbf{1}_{\{D(x_*,z)\in[2^{-k+2},2^{-k+3}]\}}\,\varphi_{\ve,k}(\z^{z(x_*)}_{2^{-k}})
\mathbf{1}_{\{x_0\notin B^{\bullet(z)}_{2^{-k}}(x_*)\cup B^{\bullet(x_*)}_{2^{-k}}(z)\}}\Bigg).
\end{equation}

Let us assume that $2^{-k+4}<1$ and argue on the event where $D(x_*,x_0)>1$. We observe that the quantity
$$\varphi_{\ve,k}(\z^{x_*(x_0)}_{2^{-k}})\,\int\mathrm{Vol}(\dd z)\,\mathbf{1}_{\{D(x_*,z)\in[2^{-k+2},2^{-k+3}]\}}\,
\varphi_{\ve,k}(\z^{z(x_*)}_{2^{-k}})
\mathbf{1}_{\{x_0\notin B^{\bullet(z)}_{2^{-k}}(x_*)\cup B^{\bullet(x_*)}_{2^{-k}}(z)\}}$$
is a function of $B^{\bullet(x_0)}_{2^{-k+4}}(x_*)$. To this end, we first note that the condition $D(x_*,z)\in[2^{-k+2},2^{-k+3}]$
implies that $z\in B_{2^{-k+4}}(x_*)$, and also that $B^{\bullet(x_0)}_{2^{-k}}(z)\subset B^{\bullet(x_0)}_{2^{-k+4}}(x_*)$
(if $y\in B^{\bullet(x_0)}_{2^{-k}}(z)$, any path from $y$ to $x_0$ has to intersect $\partial  B^{\bullet(x_0)}_{2^{-k}}(z)$,
which is contained in $B_{2^{-k+4}}(x_*)$, and this exactly means that $y\in B^{\bullet(x_0)}_{2^{-k+4}}(x_*)$).
Using \eqref{approx-exit}, it follows that, under the conditions 
$x_0\notin B^{\bullet(z)}_{2^{-k}}(x_*)\cup B^{\bullet(x_*)}_{2^{-k}}(z)$ and $D(x_*,z)\in[2^{-k+2},2^{-k+3}]$, the quantity
 $\z^{z(x_*)}_{2^{-k}}=\z^{z(x_0)}_{2^{-k}}$ is determined by
the hull $B^{\bullet(x_0)}_{2^{-k+4}}(x_*)$, and clearly the same is true for $\z^{x_*(x_0)}_{2^{-k}}$.
In addition, the condition $x_0\notin B^{\bullet(z)}_{2^{-k}}(x_*)\cup B^{\bullet(x_*)}_{2^{-k}}(z)$ can also
be expressed in terms of
the hull $B^{\bullet(x_0)}_{2^{-k+4}}(x_*)$ since it is equivalent to saying that $z$ is  in the same
connected component of $B_{2^{-k}}(x_*)^c$ as the boundary of this hull, and similarly if $x_*$ and $z$ are interchanged.

On the other hand, $F_{2^{-k+4},1}(x_*,x_0)$
is a function of $\overline{\bm_\infty\backslash B^{\bullet(x_0)}_{2^{-k+4}}(x_*)}$
with the notation of Theorem \ref{deco-hull}. This theorem, or more precisely the remark following
the statement of the theorem, shows that $B^{\bullet(x_0)}_{2^{-k+4}}(x_*)$ and 
$\overline{\bm_\infty\backslash B^{\bullet(x_0)}_{2^{-k+4}}(x_*)}$ are conditionally independent given $\z^{x_*(x_0)}_{2^{-k+4}}$
under $\N_0^{[1]}$. Recalling that $F_{2^{-k+4},1}(x_*,x_0)$ can be nonzero
only if $D(x_*,x_0)>1$, it follows that the quantity \eqref{key-tech44} is equal to
\begin{align}
\label{key-tech5}
&\frac{3}{2}\,\N_0^{[1]}\Bigg(\N_0^{[1]}\Big(F_{2^{-k+4},1}(x_*,x_0)\,\Big|\, \z^{x_*(x_0)}_{2^{-k+4}}\Big)\nonumber\\
&\qquad\times \varphi_{\ve,k}(\z^{x_*(x_0)}_{2^{-k}})\int\mathrm{Vol}(\dd z)\,\mathbf{1}_{\{D(x_*,z)\in[2^{-k+2},2^{-k+3}]\}}\,
\varphi_{\ve,k}(\z^{z(x_*)}_{2^{-k}})
\mathbf{1}_{\{x_0\notin B^{\bullet(z)}_{2^{-k}}(x_*)\cup B^{\bullet(x_*)}_{2^{-k}}(z)\}}\Bigg).
\end{align}

We now need a lemma.

\begin{lemma}
\label{geodesic-exit}
There exists a constant $C_m$ such that, for every $\ve\in(0,1/2]$ and $z>0$
$$\N_0^{[1]}\Big(F_{\ve,1}(x_*,x_0)\,\Big|\, \z^{x_*(x_0)}_{\ve}=z\Big)\leq C_m\,z^{m/2}.$$
\end{lemma}

We postpone the proof of this lemma to the Appendix. Thanks to the lemma, the quantity
\eqref{key-tech5} can be bounded above by
\begin{equation}
\label{key-tech6}
C'_m\N_0^{[1]}\Bigg( (\z^{x_*(x_0)}_{2^{-k+4}})^{m/2}\,\varphi_{\ve,k}(\z^{x_*(x_0)}_{2^{-k}})
\int\mathrm{Vol}(\dd z)\,\mathbf{1}_{\{D(x_*,z)\in[2^{-k+2},2^{-k+3}]\}}\,
\varphi_{\ve,k}(\z^{z(x_*)}_{2^{-k}})
\Bigg)
\end{equation}
where $C'_m$ is a constant. 

We have assumed that $2^{-k+4}<1$, but if  $2^{-k+4}\geq 1$, replacing $F_{2^{-k+4},1}(x_*,x_0)$
by $\mathbf{1}_{\{D(x_*,x_0)>1\}}$ immediately shows that \eqref{key-tech44} is also
bounded by a quantity similar to \eqref{key-tech6} without the term $(\z^{x_*(x_0)}_{2^{-k+4}})^{m/2}$.
For simplicity, we assume until the end of the first step that $2^{-k+4}<1$, but clearly
the bounds that follow are also valid when $2^{-k+4}\geq 1$.

We use the Cauchy-Schwarz inequality to bound the quantity \eqref{key-tech6} by
\begin{equation}
\label{key-tech7}
C'_m\,\N_0^{[1]}\Big( (\z^{x_*(x_0)}_{2^{-k+4}})^m \,\varphi_{\ve,k}(\z^{x_*(x_0)}_{2^{-k}})^2\Big)^{1/2}
\times \N^{[1]}_0\Big( \Big(\int\mathrm{Vol}(\dd z)\,\mathbf{1}_{\{D(x_*,z)\in[2^{-k+2},2^{-k+3}]\}}\,
\varphi_{\ve,k}(\z^{z(x_*)}_{2^{-k}})\Big)^2\Big)^{1/2}.
\end{equation}
Consider the first term of the product, recalling that, by definition, $\z^{x_*(x_0)}_{2^{-k}}=\z_{W_*+2^{-k}}$.
Using the bound \eqref{bound1-geo} and again the Cauchy-Schwarz inequality, we get that 
$$\N_0^{[1]}\Big( (\z^{x_*(x_0)}_{2^{-k+4}})^m \,\varphi_{\ve,k}(\z^{x_*(x_0)}_{2^{-k}})^2\Big)^{1/2}
\leq C\,(\ve 2^{3k})^m\,\N_0^{[1]}\Big((\z^{x_*(x_0)}_{2^{-k}})^{4m}\Big)^{1/4}
\times \N_0^{[1]}\Big( (\z^{x_*(x_0)}_{2^{-k+4}})^{2m}\Big)^{1/4}.$$
Thanks to the bound \eqref{bd-mom-exit}, we arrive at
\begin{equation}
\label{key-tech8}
\N_0^{[1]}\Big( (\z^{x_*(x_0)}_{2^{-k+4}})^m \,\varphi_{\ve,k}(\z^{x_*(x_0)}_{2^{-k}})^2\Big)^{1/2}
\leq C\,(\ve  2^{3k})^m \times c_{4m}^{1/4}\,2^{-2mk}\times c_{2m}^{1/4}\,2^{-m(k-4)}=C''_m \,\ve^m.
\end{equation}

We then estimate the second term in the product of \eqref{key-tech7}. By the bound \eqref{bound1-geo}, we have
\begin{align*}
&\N^{[1]}_0\Big( \Big(\int\mathrm{Vol}(\dd z)\,\mathbf{1}_{\{D(x_*,z)\in[2^{-k+2},2^{-k+3}]\}}\,
\varphi_{\ve,k}(\z^{z(x_*)}_{2^{-k}})\Big)^2\Big)^{1/2}\\
&\quad\leq C\,(\ve 2^{3k})^m\,\N^{[1]}_0\Big( \Big(\int\mathrm{Vol}(\dd z)\,\mathbf{1}_{\{D(x_*,z)\in[2^{-k+2},2^{-k+3}]\}}\,
(\z^{z(x_*)}_{2^{-k}})^m\Big)^2\Big)^{1/2}\\
&\quad\leq C\,(\ve 2^{3k})^m\,\N^{[1]}_0\Big( \mathrm{Vol}(B_{2^{-k+3}}(x_*))\int\mathrm{Vol}(\dd z)\,\mathbf{1}_{\{D(x_*,z)\in[2^{-k+2},2^{-k+3}]\}}\,
(\z^{z(x_*)}_{2^{-k}})^{2m}\Big)^{1/2} \\
&\quad= \frac{2}{3}C\,(\ve 2^{3k})^m\,\N_0\Big(\mathbf{1}_{\{D(x_*,x_0)>1\}} \mathrm{Vol}(B_{2^{-k+3}}(x_*))\int\mathrm{Vol}(\dd z)\,\mathbf{1}_{\{D(x_*,z)\in[2^{-k+2},2^{-k+3}]\}}\,
(\z^{z(x_*)}_{2^{-k}})^{2m}\Big)^{1/2},
\end{align*}
using the Cauchy-Schwarz inequality and then the definition of $\N^{[1]}_0$. In the last integral under $\N_0$, we can 
use \eqref{inter-3} to
interchange the roles of $x_0$ and $z$, and then get that the right-hand side of the last display is equal to
\begin{equation}
\label{key-tech9}
\frac{2}{3}C\,(\ve 2^{3k})^m\,\N_0\Bigg(\Big(\int \mathrm{Vol}(\dd z)\,\mathbf{1}_{\{D(x_*,z)>1\}}\Big)
\mathrm{Vol}(B_{2^{-k+3}}(x_*))\,\mathbf{1}_{\{D(x_*,x_0)\in[2^{-k+2},2^{-k+3}]\}}\,
(\z^{x_0(x_*)}_{2^{-k}})^{2m}\Bigg)^{1/2}.
\end{equation}
Fix an integer $p\geq 4$ and set $q=p/(p-1)$. We assume that $p$ is chosen sufficiently
large so that $\frac{2}{q}-\frac{1}{pq}-\frac{1}{p} > 2-\delta$. By the H\"older inequality, we have
\begin{align}
\label{key-tech10}
&\N_0\Bigg(\Big(\int \mathrm{Vol}(\dd z)\,\mathbf{1}_{\{D(x_*,z)>1\}}\Big)
\mathrm{Vol}(B_{2^{-k+3}}(x_*))\,\mathbf{1}_{\{D(x_*,x_0)\in[2^{-k+2},2^{-k+3}]\}}\,
(\z^{x_0(x_*)}_{2^{-k}})^{2m}\Bigg)\nonumber\\
&\qquad\leq \N_0\Big(\mathbf{1}_{\{D(x_*,x_0)\in[2^{-k+2},2^{-k+3}]\}}\,(\mathrm{Vol}(B_{1}(x_*)^c))^q\Big)^{1/q}\nonumber\\
&\qquad\qquad\times \N_0\Big(\mathbf{1}_{\{D(x_*,x_0)\in[2^{-k+2},2^{-k+3}]\}}\,(\mathrm{Vol}(B_{2^{-k+3}}(x_*))^p
(\z^{x_0(x_*)}_{2^{-k}})^{2pm}\Big)^{1/p}\nonumber\\
&\qquad = A^{1/q} \times B^{1/p}.
\end{align}

Let us estimate $A$ and $B$ separately. Using the fact that (conditionally on $\sigma$) $x_0$ is uniformly
distributed over $\bm_\infty$, we get 
\begin{align}
\label{key-tech11}
A&\leq \N_0\Big(\frac{1}{\sigma}\mathrm{Vol}(B_{2^{-k+3}}(x_*)\, (\mathrm{Vol}(B_{1}(x_*)^c))^q\Big)\nonumber\\
&\leq \N_0\Big(\sigma^{-q}(\mathrm{Vol}(B_{1}(x_*)^c))^{q^2}\Big)^{1/q}
\times \N_0\Big(\mathrm{Vol}(B_{2^{-k+3}}(x_*))^p\Big)^{1/p}\nonumber\\
&\leq K_q\times \Big(c_{(p)}(2^{-k+3})^{4p-2}\Big)^{1/p}\nonumber\\
&= K'_q\, 2^{-4k+2k/p},
\end{align}
using the bounds on moments of the volume of balls (Lemma \ref{volume-ball}), together with the fact that
$$\N_0\Big(\sigma^{-q}(\mathrm{Vol}(B_{1}(x_*)^c))^{q^2}\Big)<\infty.$$
The latter estimate is easily obtained by writing, for every $s>0$,
$$
\N_0\Big(\sigma^{-q}(\mathrm{Vol}(B_{1}(x_*)^c))^{q^2}\,\Big|\,\sigma = s\Big)
\leq s^{q^2-q}\,\N_0(B_1(x_*)^c\not = \varnothing\midd \sigma=s)
\leq C\,s^{q^2-q}\exp(-c\,s^{-1/3})
$$
by an estimate of \cite[Proposition 14]{Serlet}. The right-hand side of the last display is
integrable with respect to $s^{-3/2}\dd s$ as soon as $q^2-q<1/2$, which holds here because $p\geq4$. 

Let us consider now the quantity $B$. By the Cauchy-Schwarz inequality, $B\leq\sqrt{B'B''}$, where
\begin{align*}
B'&=\N_0\Big(\mathbf{1}_{\{D(x_*,x_0)\in[2^{-k+2},2^{-k+3}]\}}\,(\mathrm{Vol}(B_{2^{-k+3}}(x_*))^{2p}\Big),\\
B''&=\N_0\Big(\mathbf{1}_{\{D(x_*,x_0)\in[2^{-k+2},2^{-k+3}]\}}\,(\z^{x_0(x_*)}_{2^{-k}})^{4pm}\Big).
\end{align*}
We have first 
$$B'\leq \N_0\Big((\mathrm{Vol}(B_{2^{-k+3}}(x_*))^{2p}\Big)\leq c_{(2p)}(2^{-k+3})^{8p-2}\leq \tilde c_p\,2^{-8pk+2k},$$
using again Lemma \ref{volume-ball}. Then,
$$B''\leq \N_0\Big((\z^{x_0(x_*)}_{2^{-k}})^{4pm}\Big)=\N_0\Big((\z^{x_*(x_0)}_{2^{-k}})^{4pm}\Big)=\N_0((\z_{W_*+2^{-k}})^{4pm}),$$
where we made the convention that all quantities $\z^{x_0(x_*)}_{2^{-k}}$, $\z^{x_*(x_0)}_{2^{-k}}$, $\z_{W_*+2^{-k}}$ are equal
to $0$ 
if $D(x_*,x_0)=-W_*\leq 2^{-k}$. From \eqref{bd-mom-exit}, we now get
$$B''\leq \N_0(W_*<-2^{-k})\times \N_0^{[2^{-k}]}\Big((\z_{W_*+2^{-k}})^{4pm}\Big)\leq \frac{3}{2}\,2^{2k}\times c_{4pm}\,2^{-8pmk}= c''_p\,2^{-8pmk+2k}.$$
By combining our estimates on $B'$ and $B''$, we arrive at
$$B\leq \sqrt{B'B''} \leq \ov{c}_p\,2^{-4p(m+1)k+2k}.$$

Using also \eqref{key-tech10} and \eqref{key-tech11}, we get that the quantity \eqref{key-tech9} is bounded above by
\begin{align*}
\frac{2}{3}C\,(\ve 2^{3k})^m\times \Big( (K'_q\, 2^{-4k+2k/p})^{1/q}\,( \ov{c}_p\,2^{-4p(m+1)k+2k})^{1/p}\Big)^{1/2}
&=\ov{C}_p\,\ve^m\,2^{-k(-3m+\frac{2}{q}-\frac{1}{pq}+2(m+1)-\frac{1}{p})}\\
&\leq \ov{C}'_p\,\ve^m\,2^{-(4-m)k+\delta k},
\end{align*}
by our choice of $p$.   By combining this estimate with \eqref{key-tech8}, we
obtain that the quantity \eqref{key-tech7} is bounded above by
$$C'_m\times C''_m\ve^m\times  \ov{C}'\,\ve^m\,2^{-(4-m)k+\delta k} = \wt C_m \ve^{2m}\,2^{-(4-m)k+\delta k}.$$
Since $A'_{\ve,k}$ was bounded above by the quantity \eqref{key-tech7}, we have obtained the
desired bound for $A'_{\ve,k}$. 

\bigskip
\noindent{\bf Second step.} We now need to get a similar bound for 
$$A''_{\ve,k}=\N_0\Bigg(\mathbf{1}_{\{D(x_*,x_0)\in[2^{-k+2},2^{-k+3}]\}}\,\int_{B^{\bullet(x_0)}_{2^{-k}}(x_*)\cup B^{\bullet(x_*)}_{2^{-k}}(x_0)} \mathrm{Vol}(\dd z)\,
F_{\ve,1}(x_*,z)\,F_{\ve,1}(x_0,z)\Bigg).$$
For obvious symmetry reasons (we can interchange $x_*$ and $x_0$), it is enough to bound
\begin{align}
\label{key-tech262}
A'''_{\ve,k}:=&\,\N_0\Bigg(\mathbf{1}_{\{D(x_*,x_0)\in[2^{-k+2},2^{-k+3}]\}}\,\int_{B^{\bullet(x_0)}_{2^{-k}}(x_*)} \mathrm{Vol}(\dd z)\,
F_{\ve,1}(x_*,z)\,F_{\ve,1}(x_0,z)\Bigg)\nonumber\\
\leq&\,\N_0\Bigg(\mathbf{1}_{\{D(x_*,x_0)\in[2^{-k+2},2^{-k+3}]\}}\,\int_{B^{\bullet(x_0)}_{2^{-k}}(x_*)} \mathrm{Vol}(\dd z)\,
F_{\ve,1}(x_*,z)\,F_{\ve,2^{-k}}(x_0,z)\Bigg).
\end{align}

\begin{figure}[!h]
\label{hull-1}
 \begin{center}
 \includegraphics[width=10cm]{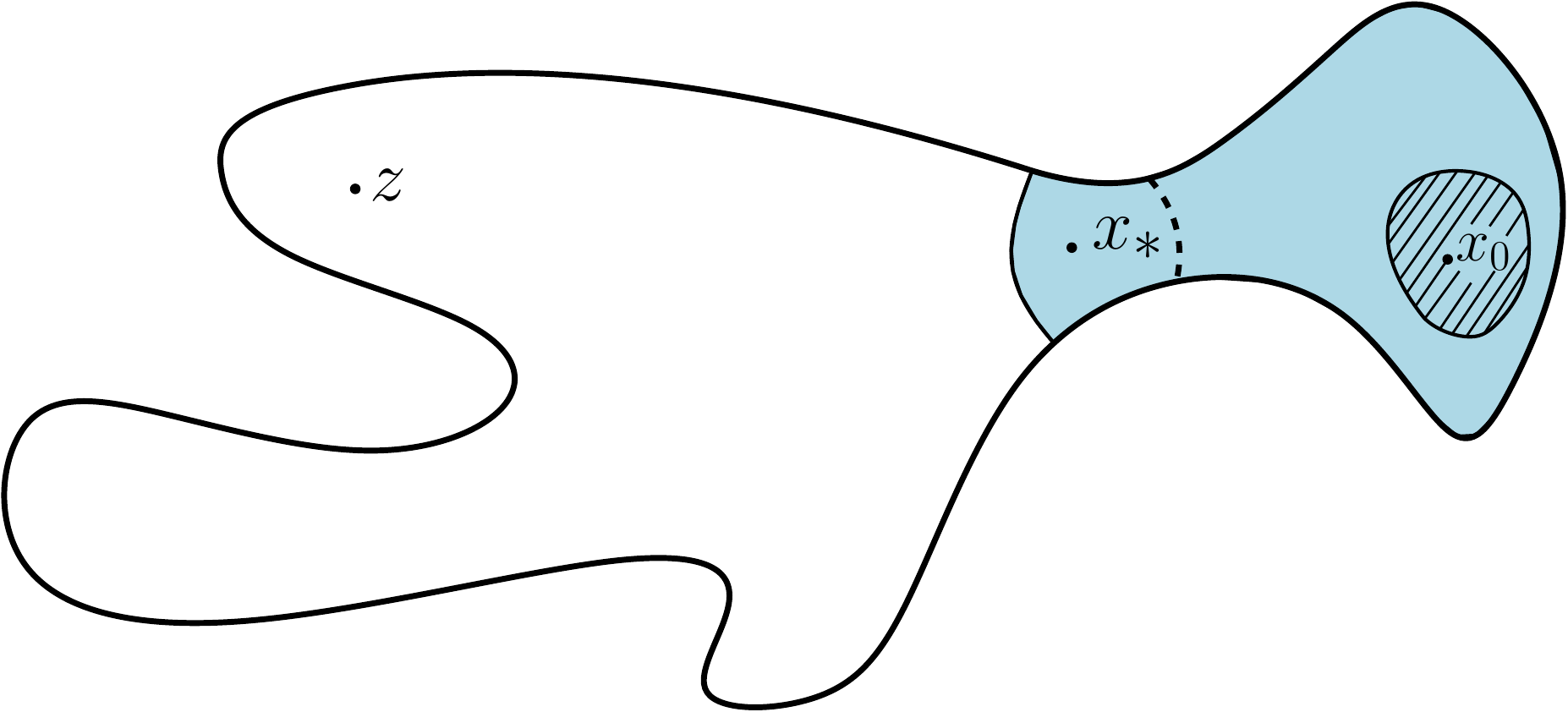}
 \caption{Illustration of a ``bottleneck'' case where $z\in B^{\bullet(x_0)}_{2^{-k}}(x_*)$
 and $B^{\bullet(x_*)}_{2^{-k}}(x_0)\subset B^{\bullet(z)}_{2^{-k}}(x_*)$. The shaded part 
 represents $B^{\bullet(z)}_{2^{-k}}(x_*)$, and the hatched part represents
 $B^{\bullet(x_*)}_{2^{-k}}(x_0)$. The dashed curve is meant to represent the boundary
 of $B^{\bullet(x_0)}_{2^{-k}}(x_*)$.}
 \end{center}
 \vspace{-5mm}
 \end{figure}

Let us argue on the event where $D(x_*,x_0)\in[2^{-k+2},2^{-k+3}]$, and note that this property
implies $B^{\bullet(x_0)}_{2^{-k}}(x_*)\cap B^{\bullet(x_*)}_{2^{-k}}(x_0)=\varnothing$. In the integral with respect to $\mathrm{Vol}(\dd z)$, 
we consider points $z$ such that $D(x_*,z)>1$ (otherwise $F_{\ve,1}(x_*,z)=0$)
and the fact that $z\in  B^{\bullet(x_0)}_{2^{-k}}(x_*)$ is equivalent to saying that $x_0$ and $z$ are different
components of $B_{2^{-k}}(x_*)^c$.  
Furthermore, we have then 
\begin{equation}
\label{key-tech263}
B^{\bullet(x_*)}_{2^{-k}}(x_0)\subset B^{\bullet(z)}_{2^{-k}}(x_*)
\end{equation}
 because the condition $D(x_*,x_0)\in[2^{-k+2},2^{-k+3}]$
ensures that every point of $\partial B^{\bullet(x_*)}_{2^{-k}}(x_0)$ belongs to the same connected component
of $B_{2^{-k}}(x_*)^c$ as $x_0$, and by the preceding observations this boundary is entirely contained in $B^{\bullet(z)}_{2^{-k}}(x_*)$. See Fig.~2 for an 
illustration. It now follows that
$B^{\bullet(z)}_{2^{-k}}(x_0)=B^{\bullet(x_*)}_{2^{-k}}(x_0)$
because clearly $z$ and $x_*$ are in the same connected component of $B_{2^{-k}}(x_0)^c$
(a geodesic path from $x_*$ to $\partial B^{\bullet(z)}_{2^{-k}}(x_*)$ does not intersect $B_{2^{-k}}(x_0)$
since it stays within distance $2^{-k}$ from $x_*$, and then \eqref{key-tech263} shows that $x_*$
is connected to $z$ by a path avoiding $B_{2^{-k}}(x_0)$). 

From the equality $B^{\bullet(z)}_{2^{-k}}(x_0)=B^{\bullet(x_*)}_{2^{-k}}(x_0)$ which holds for the relevant values of $z$, we conclude that, in the right-hand side of \eqref{key-tech262}, we can replace $F_{\ve,2^{-k}}(x_0,z)$ by $F_{\ve,2^{-k}}(x_0,x_*)$, 
and we have thus
\begin{equation}
\label{key-tech20}
A'''_k
\leq \N_0\Bigg(\mathbf{1}_{\{D(x_*,x_0)\in[2^{-k+2},2^{-k+3}]\}}\,F_{\ve,2^{-k}}(x_0,x_*)\int_{B^{\bullet(x_0)}_{2^{-k}}(x_*)} \mathrm{Vol}(\dd z)\,
F_{\ve,1}(x_*,z)\Bigg)
\end{equation}
At this point, we observe that $F_{\ve,2^{-k}}(x_0,x_*)$ is a function of the hull $B^{\bullet(x_*)}_{2^{-k}}(x_0)$, whereas one can verify that
$$\mathbf{1}_{\{D(x_*,x_0)\in[2^{-k+2},2^{-k+3}]\}}\int_{B^{\bullet(x_0)}_{2^{-k}}(x_*)} \mathrm{Vol}(\dd z)\,
F_{\ve,1}(x_*,z)$$
is a function of $\ov{\bm_\infty\backslash B^{\bullet(x_*)}_{2^{-k}}(x_0)}$ (in both cases, spaces are equipped with their
intrinsic distances). Let us explain this. First, we know that (under the condition $D(x_*,x_0)\in[2^{-k+2},2^{-k+3}]$)
we have $B^{\bullet(x_0)}_{2^{-k}}(x_*)\subset \bm_\infty\backslash B^{\bullet(x_*)}_{2^{-k}}(x_0)$, and moreover
the property $z\in B^{\bullet(x_0)}_{2^{-k}}(x_*)$ holds if and only if 
$z$ is not in the same component of $B_{2^{-k}}(x_*)^c$ as the boundary of $\bm_\infty\backslash B^{\bullet(x_*)}_{2^{-k}}(x_0)$. 
Then, assuming that $z\in B^{\bullet(x_0)}_{2^{-k}}(x_*)$ and $D(x_*,z)>1$, we observe that any geodesic
from $\partial B^{\bullet(z)}_1(x_*)$ to $x_*$ stays in the same connected component of $B_{2^{-k}}(x_*)^c$ as $z$ until it comes within distance $2^{-k}$ from $x_*$, and this component is contained in $\bm_\infty\backslash B^{\bullet(x_*)}_{2^{-k}}(x_0)$ by
\eqref{key-tech263}. 

We can thus apply Theorem \ref{deco-hull}, which (together with \eqref{inter-2}) implies that, under $\N_0(\cdot\midd D(x_0,x_*)>2^{-k})$ and conditionally on $\z^{x_0(x_*)}_{2^{-k}}$, 
$B^{\bullet(x_*)}_{2^{-k}}(x_0)$ and $\ov{\bm_\infty\backslash B^{\bullet(x_*)}_{2^{-k}}(x_0)}$ are independent, and the conditional distribution of 
$B^{\bullet(x_*)}_{2^{-k}}(x_0)$ is the law of a standard hull of radius $2^{-k}$ and perimeter $\z^{x_0(x_*)}_{2^{-k}}$. It follows that the right-hand side of
\eqref{key-tech20} is equal to
\begin{equation}
\label{key-tech21}
\N_0\Bigg(\mathbf{1}_{\{D(x_*,x_0)\in[2^{-k+2},2^{-k+3}]\}}\,\varphi_{\ve,k}(\z^{x_0(x_*)}_{2^{-k}})\,\int_{B^{\bullet(x_0)}_{2^{-k}}(x_*)} \mathrm{Vol}(\dd z)\,
F_{\ve,1}(x_*,z)\Bigg).
\end{equation}
Thanks to \eqref{inter-3}, we can now interchange the roles of $x_0$
and $z$ and get that the quantity \eqref{key-tech21} is also equal to
\begin{equation}
\label{key-tech22}
\N_0\Bigg(F_{\ve,1}(x_*,x_0) \int_{B^{\bullet(x_0)}_{2^{-k}}(x_*)} \mathrm{Vol}(\dd z)\,\mathbf{1}_{\{D(x_*,z)\in[2^{-k+2},2^{-k+3}]\}}\,\varphi_{\ve,k}(\z^{z(x_*)}_{2^{-k}})\Bigg)
\end{equation}
using once again the fact that $z\in  B^{\bullet(x_0)}_{2^{-k}}(x_*)$ is equivalent to $x_0\in B^{\bullet(z)}_{2^{-k}}(x_*)$. Now we observe that
$$F_{\ve,1}(x_*,x_0)\leq F_{2^{-k},1}(x_*,x_0)\,F_{\ve,2^{-k}}(x_*,x_0).$$
On one hand, $F_{2^{-k},1}(x_*,x_0)$ is a function of $\ov{\bm_\infty\backslash B^{\bullet(x_0)}_{2^{-k}}(x_*)}$, and on the other hand,
$$F_{\ve,2^{-k}}(x_*,x_0)\int_{B^{\bullet(x_0)}_{2^{-k}}(x_*)} \mathrm{Vol}(\dd z)\,\mathbf{1}_{\{D(x_*,z)\in[2^{-k+2},2^{-k+3}]\}}\,\varphi_{\ve,k}(\z^{z(x_*)}_{2^{-k}})$$
is a function of the hull $B^{\bullet(x_0)}_{2^{-k}}(x_*)$. For this last point, we observe that, for points $z$ that are relevant in the integral with respect
to $\mathrm{Vol}(\dd z)$, we have $B^{\bullet(x_*)}_{2^{-k}}(z)\subset B^{\bullet(x_0)}_{2^{-k}}(x_*)$
(by \eqref{key-tech263}, recalling that we interchanged the roles of $x_0$ and $z$), and we also note that, if $z$ belongs to the hull $B^{\bullet(x_0)}_{2^{-k}}(x_*)$, a geodesic (with respect to the distance $D$) from $z$ to $x_*$
may hit the boundary of the hull but cannot exit the hull, so that the intrinsic distance between $z$ and $x_*$ relative to the hull
indeed coincides with $D(x_*,z)$. Using Theorem \ref{deco-hull}, and replacing $\N_0$ by $\N_0^{[1]}$ (recall that
$F_{\ve,1}(x_*,x_0)=1$ implies $D(x_*,x_0)>1$), we get that \eqref{key-tech22} is bounded above by
\begin{align*}
&\frac{3}{2} \N^{[1]}_0\Bigg( \N^{[1]}_0(F_{2^{-k},1}(x_*,x_0)\midd \z^{x_*(x_0)}_{2^{-k}})\,
F_{\ve,2^{-k}}(x_*,x_0)\int_{B^{\bullet(x_0)}_{2^{-k}}(x_*)} \mathrm{Vol}(\dd z)\,\mathbf{1}_{\{D(x_*,z)\in[2^{-k+2},2^{-k+3}]\}}\,\varphi_{\ve,k}(\z^{z(x_*)}_{2^{-k}})\Bigg)\\
&\leq \frac{3C_m}{2} \N^{[1]}_0\Bigg( (\z^{x_*(x_0)}_{2^{-k}})^{m/2}\,
F_{\ve,2^{-k}}(x_*,x_0)\int_{B^{\bullet(x_0)}_{2^{-k}}(x_*)} \mathrm{Vol}(\dd z)\,\mathbf{1}_{\{D(x_*,z)\in[2^{-k+2},2^{-k+3}]\}}\,\varphi_{\ve,k}(\z^{z(x_*)}_{2^{-k}})\Bigg)
\end{align*}
by Lemma \ref{geodesic-exit}. Comparing with \eqref{key-tech6}, we see that the same arguments as in the first step would allow us to complete the
proof if we could replace  $F_{\ve,2^{-k}}(x_*,x_0)$ by $\varphi_{\ve,k}(\z^{x_*(x_0)}_{2^{-k}})$ in the right-hand side of the last display. Unfortunately, it is not so
easy to justify this replacement. 

Let $a>0$. To simplify notation, we write $\N^{[1],a}_0:=\N^{[1]}_0(\cdot\midd \z^{x_*(x_0)}_{2^{-k}}=a)$, so that, under $\N^{[1],a}_0$, $B^{\bullet(x_0)}_{2^{-k}}(x_*)$
is distributed as a standard hull of radius $2^{-k}$ and perimeter $a$, and we have to evaluate
$$M_{\ve,k}(a):=\N^{[1],a}_0\Bigg(F_{\ve,2^{-k}}(x_*,x_0)\int_{B^{\bullet(x_0)}_{2^{-k}}(x_*)} \mathrm{Vol}(\dd z)\,
\mathbf{1}_{\{D(x_*,z)\in[2^{-k+2},2^{-k+3}]\}}\,\varphi_{\ve,k}(\z^{z(x_*)}_{2^{-k}})\Bigg).$$
As we already explained in the proof of Theorem \ref{deco-hull}, $B^{\bullet(x_0)}_{2^{-k}}(x_*)$ corresponds 
in the Brownian snake representation to the excursions of $W$
below level $W_*+2^{-k}$. Let us write $(\omega^i)_{i\in I}$ and $\omega_*$ for these excursions, 
in a way similar to Section \ref{cons-hull}. Under $\N^{[1],a}_0$, $\omega_*$ is distributed 
according to $\N_0(\cdot\midd W_*=-2^{-k})$, and $\sum_{i\in I} \delta_{\omega^i}$ is an independent Poisson point measure with intensity $a\,\N_0(\cdot\midd W_*>-2^{-k})$.  As in Section \ref{sec:slice}, we can associate a slice $\mathbf{S}(\omega^i)$ 
with each excursion $\omega^i$, resp.~a slice $\mathbf{S}(\omega_*)$ with $\omega_*$, and this slice corresponds to
a subset of $B^{\bullet(x_0)}_{2^{-k}}(x_*)$ (this correspondence preserves the volume and is bijective except in the case 
of $\mathbf{S}(\omega_*)$, as we explained in Section \ref{slice-hull}), in such a way that the union of these subsets 
is the whole hull $B^{\bullet(x_0)}_{2^{-k}}(x_*)$, up to a set of zero volume. 
Then, we have, 
\begin{equation}
\label{key-tech23}
\int_{B^{\bullet(x_0)}_{2^{-k}}(x_*)}\mathrm{Vol}(\dd z)\,\mathbf{1}_{\{D(x_*,z)\in[2^{-k+2},2^{-k+3}]\}}
\,\varphi_{\ve,k}(\z^{z(x_*)}_{2^{-k}}) =\Psi_{\ve,k}(\omega_*) + \sum_{i\in I} \Psi_{\ve,k}(\omega^i),
\end{equation}
where, for $\omega\in\S_0$ such that $W_*(\omega)\geq -2^{-k}$,
$$\Psi_{\ve,k}(\omega):=\int_{\mathbf{S}(\omega)}\mathrm{Vol}(\dd \wt z)\,\mathbf{1}_{\{\tilde D(\tilde x_*,\tilde z)\in[2^{-k+2}-(W_*+2^{-k}),2^{-k+3}-(W_*+2^{-k})]\}}\,\varphi_{\ve,k}(\z^{\tilde z(\tilde x_*)}_{2^{-k}}).$$
We have used the notation of Section \ref{sec:slice} and the fact that, for every point $\wt z$ of the slice $\mathbf{S}(\omega)$ at distance greater than
$2^{-k}$ from  the distinguished point $\wt x_*$, we can define the hull
of radius $2^{-k}$ centered at $\wt z$ (relative to $\wt x_*$) and its boundary size $\z^{\tilde z(\tilde x_*)}_{2^{-k}}$
 via the analog of formula \eqref{exit-approx2}. 
The point in \eqref{key-tech23}
is to observe that any point $z$ of $B^{\bullet(x_0)}_{2^{-k}}(x_*)$ such that $D(x_*,z)\in[2^{-k+2},2^{-k+3}]$
corresponds to a point $\wt z$ of the slice $\SS(\omega^i)$ for some $i\in I$ (or of $\SS(\omega_*)$),
such that $\wt D(\wt x_*,\wt z)= D(x_*,z)-(2^{k}+W_*(\omega_i))$ (or $\wt D(\wt x_*,\wt z)= D(x_*,z)$), and moreover
the hull of radius $2^{-k}$ centered at $z$ in $\bm_\infty$ is contained in $B^{\bullet(x_0)}_{2^{-k}}(x_*)$ and identified 
with the same hull
centered at $\wt z$ in the slice, so that these two hulls have the same boundary size. To check the last property, we also
use the easy fact that the identification of the slice  $\SS(\omega^i)$ with a subset of $B^{\bullet(x_0)}_{2^{-k}}(x_*)$
is isometric on the set $\{\wt z \in \SS(\omega^i): \wt D(\wt x_*,\wt z)\geq 2^{-k+1}\}$ (and similarly for $\SS(\omega_*)$).

As in the proof of Proposition  \ref{distinct-geo}  (cf.~formula \eqref{nber-exc}), we have 
$F_{\ve,2^{-k}}(x_*,x_0)=\mathbf{1}_{\{N\geq m\}}$,
where $N:=\#\{i\in I:W_*(\omega^i)<-2^{-k}+\ve\}$. Therefore, we can write
\begin{equation}
\label{key-tech29}
M_{\ve,k}(a)= \N^{[1],a}_0\Bigg(\mathbf{1}_{\{N\geq m\}} \Big(\Psi_{\ve,k}(\omega_*)+\sum_{i\in I} \Psi_{\ve,k}(\omega^i)\Big)\Bigg).
\end{equation}
Under $\N^{[1],a}_0$, the variable $N$ follows the Poisson distribution with parameter $\frac{3a}{2}((2^{-k}-\ve)^{-2}-(2^{-k})^{-2})$ and thus,
as in \eqref{bound1-geo},
$$\N^{[1],a}_0(N\geq m)\leq C (\ve a 2^{3k})^m.$$
Since $\omega_*$ is independent of $\sum_{i\in I} \delta_{\omega^i}$, we have plainly
\begin{equation}
\label{key-tech30}
\N^{[1],a}_0\Big(\mathbf{1}_{\{N\geq m\}} \Psi_{\ve,k}(\omega_*)\Big)= \N^{[1],a}_0(N\geq m)\,\N^{[1],a}_0(\Psi_{\ve,k}(\omega_*))\leq C (\ve a 2^{3k})^m \N^{[1],a}_0(\Psi_{\ve,k}(\omega_*)).
\end{equation}
Set $I'=\{i\in I:W_*(\omega^i)<-2^{-k}+\ve\}$ (so that $N=\# I'$). Since $N$ is also independent of the point measure $\sum_{i\in I\backslash I'} \delta_{\omega^i}$, we get
similarly
\begin{equation}
\label{key-tech31}
\N^{[1],a}_0\Big(\mathbf{1}_{\{N\geq m\}} \sum_{i\in I\backslash I'}\Psi_{\ve,k}(\omega^i)\Big)\leq C (\ve a 2^{3k})^m \N^{[1],a}_0\Big(\sum_{i\in I\backslash I'}\Psi_{\ve,k}(\omega^i)\Big).
\end{equation}
The delicate part is to estimate 
\begin{align}
\label{key-tech32}
\N^{[1],a}_0\Big(\mathbf{1}_{\{N\geq m\}} \sum_{i\in I'}\Psi_{\ve,k}(\omega^i)\Big)&=\sum_{p=m}^\infty p\,\N^{[1],a}_0(N=p)\,\N_0(\Psi_{\ve,k}(\omega)\midd-2^{-k}<W_*<-2^{-k}+\ve)\nonumber\\
&\leq C'(\ve a 2^{3k})^m\,\N_0(\Psi_{\ve,k}(\omega)\midd-2^{-k}<W_*<-2^{-k}+\ve).
\end{align}
We use 
the spine decomposition of Section \ref{spine-decomp} to verify that, for every $r\in(2^{-k}-\ve,2^{-k})$, one can 
couple a snake trajectory $\omega_{(r)}$ distributed according to $\N_0(\cdot\midd W_*=-r)$
with a snake trajectory $\omega_{(2^{-k})}$ distributed according to $\N_0(\cdot\midd W_*=-2^{-k})$ in such a way that 
the following holds. The slice $\SS(\omega_{(r)})$ is identified isometrically (and in a manner preserving both
the volume measure and the first distinguished point $\wt x_*$) to a closed subset of $\SS(\omega_{(2^{-k})})$, and moreover,
if $z$ is a point of $\SS(\omega_{(r)})$ whose distance from $\wt x_*$ is greater than $2^{-k+1}$, the hull
of radius $2^{-k}$ centered at $z$ in $\SS(\omega_{(r)})$ is identified with the same hull in $\SS(\omega_{(1)})$,
and these two hulls have the same boundary size --- here we omit a few details that are left to
the reader. It follows that
$\Psi_{\ve,k}(\omega_{(r)})\leq \wt\Psi_{\ve,k}(\omega_{(2^{-k})})$, 
where $\wt\Psi_{\ve,k}\geq \Psi_{\ve,k}$ is defined as $\Psi_{\ve,k}$, except that the interval $[2^{-k+2},2^{-k+3}]$
is replaced by $[2^{-k+2}-2^{-k},2^{-k+3}]$. Hence we have
\begin{equation}
\label{key-tech33}
\N_0(\Psi_{\ve,k}(\omega)\midd-2^{-k}<W_*<-2^{-k}+\ve)\leq \N_0(\wt\Psi_{\ve,k}(\omega)\midd W_*=2^{-k})= \N^{[1],a}_0(\wt\Psi_{\ve,k}(\omega_*)). 
\end{equation}

Finally, using \eqref{key-tech29}, \eqref{key-tech30}, \eqref{key-tech31}, \eqref{key-tech32}, \eqref{key-tech33} and the
analog of \eqref{key-tech23} where $[2^{-k+2},2^{-k+3}]$
is replaced by $[2^{-k+2}-2^{-k},2^{-k+3}]$ and $\Psi_{\ve,k}$ is replaced by $\wt\Psi_{\ve,k}$, we get the existence
of a constant $C''$ such that
$$M_{\ve,k}(a)\leq C''(\ve a 2^{3k})^m\,\N^{[1],a}_0\Bigg(\int_{B^{\bullet(x_0)}_{2^{-k}}(x_*)} \mathrm{Vol}(\dd z)\,
\mathbf{1}_{\{D(x_*,z)\in[2^{-k+2}-2^{-k},2^{-k+3}]\}}\,\varphi_{\ve,k}(\z^{z(x_*)}_{2^{-k}})\Bigg).$$
It follows that $A'''_{\ve,k}$ is bounded by a constant times
 $$\N^{[1]}_0\Bigg( (\z^{x_*(x_0)}_{2^{-k}})^{m/2}\,
\times (C''(\ve \z^{x_*(x_0)}_{2^{-k}} 2^{3k})^m \int_{B^{\bullet(x_0)}_{2^{-k}}(x_*)} \mathrm{Vol}(\dd z)\,\mathbf{1}_{\{D(x_*,z)\in[2^{-k+2}-2^{-k},2^{-k+3}]\}}\,\varphi_{\ve,k}(\z^{z(x_*)}_{2^{-k}})\Bigg),$$
and  we get an upper bound by replacing the integral over $B^{\bullet(x_0)}_{2^{-k}}(x_*)$ by the same 
integral over $\bm_\infty$. The very same arguments that we used in the first step 
to bound the quantity \eqref{key-tech6}, now show that $A'''_{\ve,k}$ is bounded above by
a constant times $\ve^{2m}2^{-(4-m)k+\delta k}$.
This completes the proof of Lemma \ref{key-lem}. \endproof

\begin{lemma}
\label{key-lem2}
Let $\alpha\in (0,4-m)$. There exists a constant $C_\alpha$
such that, for every $\ve\in(0,1/2)$,
\begin{equation}
\label{key2-tech1}
\wt\N_0\Bigg(\int\!\!\int  \mathrm{Vol}(\dd x)\mathrm{Vol}(\dd y)\,\mathbf{1}_{\{D(x,y)<\ve\}}\,D(x,y)^{-\alpha}\,
F^{(m)}_{\ve,1}(x,x_*)\Bigg)
\leq C_\alpha\,\ve^{2m}.
\end{equation}
\end{lemma}

\proof As previously, we write $F_{\ve,r}(x,y)$ instead of $F^{(m)}_{\ve,r}(x,y)$ in the proof. Let $\wh A_{\ve,k}$ denote the left-hand side of \eqref{key2-tech1}. In a way similar
to the beginning of the proof of Lemma \ref{key-lem}, we can use the symmetry properties
of the Brownian sphere to write $\wh A_{\ve,k}$ in a different form. We write 
$$\wh A_{\ve,k}= \N_0\Big(\sigma\,\int\!\!\int \frac{\mathrm{Vol}(\dd x)}{\sigma}\,\frac{\mathrm{Vol}(\dd y)}{\sigma}\,
\wh\Gamma_{\ve,k}(x_*,x,y)\Big),$$
with an appropriate function $\wh\Gamma_{\ve,k}$, and observe that we have also
\begin{align*}
\wh A_{\ve,k}
&= \N_0\Big(\sigma\,\int \frac{\mathrm{Vol}(\dd z)}{\sigma}\,
\wh\Gamma_{\ve,k}(x_0,x_*,z)\Big)\\
&=\N_0\Big(F_{\ve,1}(x_*,x_0)\int\mathrm{Vol}(\dd z)\,\mathbf{1}_{\{D(x_*,z)<\ve\}}
D(x_*,z)^{-\alpha}\Big)\\
&=\frac{3}{2}\,\N^{[1]}_0\Big(F_{\ve,1}(x_*,x_0)\int \mathrm{Vol}(\dd z)\,\mathbf{1}_{\{D(x_*,z)<\ve\}}
D(x_*,z)^{-\alpha}\Big).
\end{align*}
The quantity $F_{\ve,1}(x_*,x_0)$ is a function of $\ov{\bm_\infty\backslash B^{\bullet(x_0)}_\ve(x_*)}$, whereas
$\int \mathrm{Vol}(\dd z)\,\mathbf{1}_{\{D(x_*,z)<\ve\}}
D(x_*,z)^{-\alpha}$
is a function of $B^{\bullet(x_0)}_\ve(x_*)$. We can thus apply Theorem \ref{deco-hull} to obtain that
the right-hand side of the last display is also equal to 
$$\frac{3}{2}\,\N^{[1]}_0\Big(\N^{[1]}_0(F_{\ve,1}(x_*,x_0)\midd \z^{x_*(x_0)}_\ve)\int \mathrm{Vol}(\dd z)\,\mathbf{1}_{\{D(x_*,z)<\ve\}}
D(x_*,z)^{-\alpha}\Big).$$
From Lemma \ref{geodesic-exit}, we get the bound
$$\wh A_{\ve,k}\leq \frac{3}{2}\,C_m\,\N^{[1]}_0\Big((\z^{x_*(x_0)}_\ve)^{m/2}\int \mathrm{Vol}(\dd z)\,\mathbf{1}_{\{D(x_*,z)<\ve\}}
D(x_*,z)^{-\alpha}\Big).$$
The remaining part of the argument is now easy. Write $k(\ve)\geq 1$ for the smallest integer such
that $2^{-k(\ve)}<\ve$. Fix $\kappa\in(0,1)$ such that $\alpha-(4-m)+\kappa<0$. For every $k\geq k(\ve)$, 
we have
\begin{align*}
&\N^{[1]}_0\Big((\z^{x_*(x_0)}_\ve)^{m/2}\int \mathrm{Vol}(\dd z)\,\mathbf{1}_{\{2^{-k}<D(x_*,z)\leq 2^{-k+1}\}}
D(x_*,z)^{-\alpha}\Big)\\
&\qquad\leq 2^{k\alpha}\,\N^{[1]}_0\Big((\z^{x_*(x_0)}_\ve)^{m/2}\,\mathrm{Vol}(B_{2^{-k+1}}(x_*))\Big)\\
&\qquad\leq 2^{k\alpha}\,\N^{[1]}_0\Big((\z^{x_*(x_0)}_\ve)^{m}\Big) ^{1/2}\,\N^{[1]}_0\Big((\mathrm{Vol}(B_{2^{-k+1}}(x_*)))^2\Big)^{1/2}\\
&\qquad\leq C_{m,\kappa}\,2^{k\alpha}\times \ve^m\times 2^{-(4-\kappa)k},
\end{align*}
using Lemma \ref{volume-ball} and the bound \eqref{bd-mom-exit}. By summing over $k\geq k(\ve)$,
we arrive at
\begin{align*}
\N^{[1]}_0\Big((\z^{x_*(x_0)}_\ve)^{m/2}\int \mathrm{Vol}(\dd z)\,\mathbf{1}_{\{D(x_*,z)<\ve\}}
D(x_*,z)^{-\alpha}\Big)&\leq C_{m,\kappa}\,\ve^m\,\sum_{k=k(\ve)}^\infty 2^{-k(4-\alpha-\kappa)}\\
&\leq C'_{m,\kappa}\,\ve^m \,2^{-k(\ve)(4-\alpha-\kappa)}\\
&\leq C'_{m,\kappa}\,\ve^{2m}.
\end{align*}
This completes the proof. \endproof

\section{Proof of Theorem \ref{main-th}}
\label{sec-proof-main}

As previously, $m\in\{1,2,3\}$ is fixed.
Recall the definition of $\wt F^{(m)}_{\ve,1}$ in Section \ref{sec:first-mom}.
For every $\ve\in(0,1/32)$, we introduce the measure $\nu_\ve$ on
$\bm_\infty$ defined by 
$$\nu_\ve(\dd x):= \ve^{-m}\,\wt F^{(m)}_{\ve,1}(x,x_*)\mathbf{1}_{\{D(x,x_*)< 2\}}\,\mathrm{Vol}(\dd x).$$
We use the notation
$$R^{\mathrm{max}}=\max\{D(x,x_*):x\in \bm_\infty\}$$
and note that $\nu_\ve$ is the zero measure if $R^{\mathrm{max}}<1$ (recall that $\wt F^{(m)}_{\ve,1}(x,x_*)=0$
if $D(x,x_*)<1$). We will then argue under the finite measure 
$$\wt\N^\star_0:=\wt \N_0(\cdot\cap\{ R^{\mathrm{max}}\geq 1\}).$$
As an immediate consequence of Proposition \ref{prop-first-mom}, we have
\begin{equation}
\label{first-mom-bd}
\wt \N_0^\star(\langle \nu_\ve,1\rangle) \geq c,
\end{equation}
with a positive constant $c$ independent of $\ve$. On the other hand, if $\delta\in(0,1)$ is fixed, we can
use Lemma \ref{key-lem} and Lemma \ref{key-lem2} to bound the integral under $\wt\N_0^\star$ of the 
quantity
\begin{align*}
&\int\!\!\int \nu_\ve(\dd x)\nu_\ve(\dd y)\,D(x,y)^{-(4-m-\delta)}\\
&\qquad \leq \ve^{-2m}\int\!\!\int \mathrm{Vol}(\dd x)\mathrm{Vol}(\dd y)\,F^{(m)}_{\ve,1}(x,x_*)\,F^{(m)}_{\ve,1}(y,x_*)\,\mathbf{1}_{\{D(x,y)< 4\}}\,D(x,y)^{-(4-m-\delta)},
\end{align*}
where we used the trivial bound $\wt F^{(m)}_{\ve,1}(x,x_*)\leq F^{(m)}_{\ve,1}(x,x_*)$, and the fact that
$D(x,x_*)< 2$ and $D(y,x_*)< 2$ imply $D(x,y)< 4$. Let $k(\ve)$ be the greatest integer
such that $2^{-k}>2\ve$. Using Lemma \ref{key-lem}, we have
\begin{align*}
&\wt\N_0^\star\Bigg(\ve^{-2m}\int\!\!\int \mathrm{Vol}(\dd x)\mathrm{Vol}(\dd y)\,F^{(m)}_{\ve,1}(x,x_*)\,F^{(m)}_{\ve,1}(y,x_*)\,\mathbf{1}_{\{2^{-k(\ve)+2}\leq D(x,y)< 4\}}\,D(x,y)^{-(4-m-\delta)}\Bigg)\\
&=\sum_{k=1}^{k(\ve)} \wt\N_0^\star\Bigg(\ve^{-2m}\int\!\!\int \mathrm{Vol}(\dd x)\mathrm{Vol}(\dd y)\,F^{(m)}_{\ve,1}(x,x_*)\,F^{(m)}_{\ve,1}(y,x_*)\,\mathbf{1}_{\{2^{-k+2}\leq D(x,y)< 2^{-k+3}\}}\,D(x,y)^{-(4-m-\delta)}\Bigg)\\
&\leq \sum_{k=1}^{k(\ve)} C_{(\delta/2)}\,2^{-(4-m)k+(\delta/2) k}\times 2^{-(-k+2)(4-m-\delta)}\\
&\leq C_{(\delta/2)}\sum_{k=1}^\infty 2^{-k\delta/2}\\
&=C'_{(\delta)}
\end{align*}
with some constant $C'_{(\delta)}$. On the other hand, using the trivial bound
$F^{(m)}_{\ve,1}(y,x_*)\leq 1$,  the fact that $2^{-k(\ve)}\leq 4\ve$, and Lemma \ref{key-lem2}, we have 
\begin{align*}
&\wt\N_0^\star\Bigg(\ve^{-2m}\int\!\!\int \mathrm{Vol}(\dd x)\mathrm{Vol}(\dd y)\,F^{(m)}_{\ve,1}(x,x_*)\,F^{(m)}_{\ve,1}(y,x_*)\,\mathbf{1}_{\{D(x,y)< 2^{-k(\ve)+2}\}}\,D(x,y)^{-(4-m-\delta)}\Bigg)\\
&\leq \wt\N_0^\star\Bigg(\ve^{-2m}\int\!\!\int \mathrm{Vol}(\dd x)\mathrm{Vol}(\dd y)\,F^{(m)}_{16\ve,1}(x,x_*)\,\mathbf{1}_{\{D(x,y)\leq 16\ve\}}\,D(x,y)^{-(4-m-\delta)}\Bigg)\\
&\leq C''_{(\delta)},
\end{align*}
with some constant $C''_{(\delta)}$. Summarizing, we have
\begin{equation}
\label{Frostman-bd}
\wt\N_0^\star\Bigg(\int\!\!\int \nu_\ve(\dd x)\nu_\ve(\dd y)\,D(x,y)^{-(4-m-\delta)}\Bigg)\leq K_{(\delta)}
\end{equation}
for a certain constant $K_{(\delta)}$ depending only on $\delta$. Since the measure $\nu_\ve(\dd x)\nu_\ve(\dd y)$ is supported 
on pairs $(x,y)$ such that $D(x,y)<4$, the bound \eqref{Frostman-bd} also implies that
\begin{equation}
\label{second-mom-bd}
\wt \N_0^\star(\langle \nu_\ve,1\rangle^2) \leq 64\,K_{(\delta)}.
\end{equation}
From \eqref{first-mom-bd} and \eqref{second-mom-bd}, a standard application of the
Cauchy-Schwarz inequality shows that we can find two positive constants $a$ and $c_0$
such that
$$\wt \N_0^\star(\langle \nu_\ve,1\rangle\geq a) \geq c_0.$$
Finally, using \eqref{Frostman-bd} and \eqref{second-mom-bd}, we can find $A>0$ large
enough such that
$$\wt \N_0^\star\Bigg( \Big\{\int\!\!\int \nu_\ve(\dd x)\nu_\ve(\dd y)\,D(x,y)^{-(4-m-\delta)}\leq A\Big\}
\cap \Big\{ a\leq \langle \nu_\ve,1\rangle \leq A\Big\}\Bigg) \geq c_0/2.$$

Then, let $(\ve_n)_{n\in\N}$ be a sequence in $(0,1/32)$ that converges to $0$. The event
$$\Theta:=\limsup_{n\to\infty} \Bigg( \Big\{\int\!\!\int \nu_{\ve_n}(\dd x)\nu_{\ve_n}(\dd y)\,D(x,y)^{-(4-m-\delta)}\leq A\Big\}
\cap \Big\{ a\leq \langle \nu_{\ve_n},1\rangle \leq A\Big\}\Bigg)$$
has $\wt\N_0^\star$-measure at least $c_0/2$. Let us argue on the event $\Theta$. On this event, we can
find a (random) subsequence $(\nu_{\ve_{n_p}})_{p\in\N}$ that converges weakly to a limiting
nonzero finite measure $\nu_0$ such that
$$\int\!\!\int \nu_0(\dd x)\nu_0(\dd y)\,D(x,y)^{-(4-m-\delta)}\leq A <\infty.$$
We claim that $\nu_0$ is supported on the set $\mathfrak{S}^{(m+1)}$ of $(m+1)$-geodesic stars. If our claim holds,
an application of the classical Frostman lemma shows that 
$\mathrm{dim}(\mathfrak{S}^{(m+1)})\geq 4-m-\delta$ on the event $\Theta$. 

Let us justify our claim. Let $x$ belong to the topological support of $\nu_0$. If $V$ is an open neighborhood
of $x$ in $\bm_\infty$, then, for $p$ large enough, we must have $\nu_{\ve_{n_p}}(V)>0$ and consequently
there exists a point $y$ of $V$ such that $\wt F^{(m)}_{\ve_{n_p},1}(y,x_*)=1$. It follows that
we can find a sequence $(x_n)_{n\in\N}$ in $\bm_\infty$ that converges to $x$, and a sequence
$(\ve'_n)_{n\in\N}$ of positive reals converging to $0$, such that, for every $n\in\N$, $\wt F^{(m)}_{\ve'_n,1}(x_n,x_*)=1$.
This means that there exist geodesics $(\xi^{(n)}_0(t))_{t\in[0,1]},(\xi^{(n)}_1(t))_{t\in[0,1]},\ldots,(\xi^{(n)}_m(t))_{t\in[0,1]}$
that terminate at $x_n$ and are such that, for every $0\leq i<j\leq m$, we have
$$D(\xi^{(n)}_i(t),\xi^{(n)}_j(t))\geq \delta_k$$
for every $t\in[1-2^{-k-1},1-2^{-k-2}]$ and $k\geq 1$ such that $2^{-k-4}\geq \ve'_n$. By a compactness argument,
up to extracting subsequences, we may assume that, for every $i\in\{0,\ldots,m\}$, 
$$\xi^{(n)}_i(t)\build{\la}_{n\to\infty}^{} \xi^{(\infty)}_i(t)\;,\qquad \hbox{uniformly in }t\in[0,1],$$
where the limit $(\xi^{(\infty)}_i(t))_{t\in[0,1]}$ must be a geodesic path that terminates at $x$. Furthermore,
we have for every $0\leq i<j\leq m$,
$$D(\xi^{(\infty)}_i(t),\xi^{(\infty)}_j(t))\geq \delta_k,$$
for every $t\in[1-2^{-k-1},1-2^{-k-2}]$ and every integer $k\geq 1$, and this ensures that the sets $\{\xi^{(\infty)}_i(t):t\in[3/4,1)\}$ are disjoint,
so that $x$ is an $(m+1)$-geodesic star, proving our claim.

At this point, we have proved that the dimension of the set of all $(m+1)$-geodesic stars is at least $4-m-\delta$
on an event of positive $\wt\N^\star_0$-measure. Clearly, we can replace $\wt\N^\star_0$ by $\wt\N_0$
or $\N_0$. To simplify notation, set $\N_0^{\{a\}}:=\N_0(\cdot\midd W_*=-a)$ for every $a>0$. Via a scaling argument,
we also get that the dimension of the set of $(m+1)$-geodesic star is at least $4-m-\delta$ on an event
of positive $\N_0^{\{a\}}$-probability. We now want to argue that the latter property even holds on
an event of full $\N_0^{\{a\}}$-probability, and we need to a kind of zero-one law argument,
for which it is more convenient to use the Brownian plane.

\begin{lemma}
\label{coupling-Bplane}
On the same probability space, we can construct both the Brownian plane $\mathcal{P}$ and a two-pointed random metric space
$\bm^{\{1\}}_\infty$ distributed according to the law of $\bm_\infty$ under $\N_0^{\{1\}}$, in such a way that, for every $\ve\in(0,1)$, there is an event $E_\ve$
of positive probability and independent of $\bm^{\{1\}}_\infty$, such that the following holds.
Write $x_*^{\{1\}}$ and $x_0^{\{1\}}$ for the two distinguished points of $\bm^{\{1\}}_\infty$, and 
let $B^\bullet_{1-\ve}(\bm^{\{1\}}_\infty)$ stand for the hull defined as the complement
of the connected component containing $x_0^{\{1\}}$ of the complement of the closed ball of radius $1-\ve$ centered 
at $x_*^{\{1\}}$ in $\bm^{\{1\}}_\infty$. Similarly, write $B^\bullet_{1-\ve}(\mathcal{P})$ for the complement of the 
unbounded component of the complement of the closed ball of radius $1-\ve$ centered at
the distinguished point of $\mathcal{P}$. On the event $E_\ve$, the interior
$\mathrm{Int}(B^\bullet_{1-\ve}(\bm^{\{1\}}_\infty))$ equipped with its intrinsic metric 
is isometric to $\mathrm{Int}(B^\bullet_{1-\ve}(\mathcal{P}))$ equipped with its intrinsic metric.
\end{lemma}

This lemma is obtained by comparing the construction of the Brownian plane
in \cite{CLG} with the spine decomposition of $\N_0^{\{1\}}$ in Section \ref{spine-decomp}.
We refer to the Appendix below for a detailed argument.

With the notation of the lemma, we have  a.s.,
\begin{equation}
\label{complement-hull}
\bm^{\{1\}}_\infty \backslash \bigcup_{\ve>0} B^\bullet_{1-\ve}(\bm^{\{1\}}_\infty)=\{x^{\{1\}}_0\}.
\end{equation}
To justify this, argue under $\N^{\{1\}}_0(\dd\omega)$ and observe that if 
$x=\Pi(a)\in\bm_\infty\backslash \{x_0\}$, labels along the line segment from $a$
to $\rho_{(\omega)}$ in $\t_{(\omega)}$ must take negative values, which ensures by the bound \eqref{cactus-bd}
that $x$ belongs to $B^\bullet_{1-\ve}(x_*)$ for $\ve>0$ small enough.

It now follows from \eqref{complement-hull} and the considerations 
preceding the lemma that, for $\ve>0$ small enough, the set of all $(m+1)$-geodesic stars of 
$\bm^{\{1\}}_\infty$ that lie in $\mathrm{Int}(B^\bullet_{1-\ve}(\bm^{\{1\}}_\infty))$ has dimension at least
$4-m-\delta$ with positive probability. Hence (here we use the fact that $E_\ve$
is independent of $\bm^{\{1\}}_\infty$), the set of all $(m+1)$-geodesic stars of $\mathcal{P}$
that lie in $\mathrm{Int}(B^\bullet_{1-\ve}(\mathcal{P}))$ has also dimension at least
$4-m-\delta$ with positive probability, for $\ve>0$ small enough. The scaling invariance
of $\mathcal{P}$ now shows that, for every $a>0$, the event where the set of all $(m+1)$-geodesic stars of $\mathcal{P}$
that lie in $\mathrm{Int}(B^\bullet_{a}(\mathcal{P}))$ has dimension at least
$4-m-\delta$ has the same (positive) probability. Writing $\mathfrak{S}^{(m+1)}(\mathcal{P})$
for the set of all $(m+1)$-geodesic stars of $\mathcal{P}$, we get that the event
$$\bigcap_{a>0} \Big\{\mathrm{dim}(\mathfrak{S}^{(m+1)}(\mathcal{P})\cap B^\bullet_{a}(\mathcal{P})) \geq 4-m-\delta\Big\}$$
has also positive probability. However, using the construction of $\mathcal{P}$ given in \cite{CLG}
(see below the proof of Lemma \ref{coupling-Bplane}), it is not hard to verify that the latter event belongs to an asymptotic $\sigma$-field 
which contains only events of probability $0$ or $1$. We thus get that the property 
$$\mathrm{dim}(\mathfrak{S}^{(m+1)}(\mathcal{P})\cap B^\bullet_{a}(\mathcal{P})) \geq 4-m-\delta$$
holds for every $a>0$, a.s. Since $\delta\in(0,1)$ was arbitrary, we conclude that
$$\mathrm{dim}(\mathfrak{S}^{(m+1)}(\mathcal{P})\cap B^\bullet_{a}(\mathcal{P})) \geq 4-m,$$
for every $a>0$, a.s. Finally, using the coupling between the Brownian sphere and the
Brownian plane found in \cite[Theorem 1]{Plane}, one gets that the same property
holds for the Brownian sphere. 

\section*{Appendix}

\subsection*{Proof of Lemma \ref{geodesic-exit}}

On the event $\{W_*<-1\}$, we define 
$M_\ve$ as the number of excursions below $W_*+1$
that hit $W_*+\ve$. As in the proof of Proposition \ref{distinct-geo}, we have $\N^{[1]}_0$ a.e.,
$$F^{(m)}_{\ve,1}(x_*,x_0)=\mathbf{1}_{\{M_\ve \geq m+1\}}.$$
So we have to bound $\N^{[1]}_0(M_\ve\geq m+1\midd \z_{W_*+\ve}=z)$
for $m\in\{1,2,3\}$ (recall that $\z^{x_*(x_0)}_\ve=\z_{W_*+\ve}$). To this end, we will rely on explicit
calculations. For every $a>0$, we write $h_a(z)$ for the density of the law
of $\z_{-a}$ under $\N_0(\cdot \cap\{\z_{-a}\not =0\})$, as given in \cite[Proposition 3]{Spine}:
\begin{equation}
\label{def-h}
h_a(z):= \Big(\frac{3}{2a^2}\Big)^2 \, \psi(\frac{3z}{2a^2})
\end{equation}
where 
\begin{equation}
\label{def-psi}
\psi(x)=\frac{2}{\sqrt{\pi}}(x^{1/2} + x^{-1/2}) - 2(x+\frac{3}{2})\,e^x\,\hbox{erfc}(\sqrt{x}),\quad x>0.
\end{equation}
We also recall from \cite[Corollary 13]{Disks} that, for every $a>0$, the density
of $\z_{W_*+a}$ under $\N_0^{[a]}$ is the function
\begin{equation}
\label{density1}
z\mapsto \frac{1}{a} \sqrt{\frac{3}{2\pi z}}\exp\Big(-\frac{3z}{2a^2}\Big).
\end{equation}

\begin{lemma}
\label{App-Lem1}
For every $\ve\in(0,1)$, we have
$$
\N^{[1]}_0(M_\ve=1)=1-\ve,\quad
\N^{[1]}_0(M_\ve=2)=\frac{1}{2}(1-\ve)\Big( 1 - (1-\ve)^2\Big),\quad
\N^{[1]}_0(M_\ve=3)=\frac{3}{8}(1-\ve) \Big(1-(1-\ve)^2\Big)^2.
$$
\end{lemma}

\proof By \cite[Proposition 12]{Disks} (see also the remark after this proposition), the conditional
distribution of $M_\ve-1$ under $\N^{[1]}_0$ knowing that $\z_{W_*+1}=z$ is Poisson with
parameter 
$$z\,\N_0(-1<W_*<-1+\ve)=z\,\Big(\frac{3}{2(1-\ve)^2}-\frac{3}{2}\Big).$$
Since the distribution of $\z_{W_*+1}$ is given by \eqref{density1}, the formulas 
of the lemma follow by straightforward calculations. \endproof

\begin{lemma}
\label{App-Lem2}
Let $\ve\in(0,1)$. The law of $\z_{-1+\ve}$ under $\N_0(\cdot\midd W_*=-1)$ has density
$$f_\ve(z):=\ve^{-3}\,z\,\exp\Big(-\frac{3z}{2\ve^2}\Big) h_{1-\ve}(z).$$
The law of $\z_{-1+\ve}$ under $\N_0(\cdot\midd -1<W_*<-1+\ve)$ has density
$$\tilde f_\ve(z):=\Big(\frac{3}{2(1-\ve)^2}-\frac{3}{2}\Big)^{-1}\,\exp\Big(-\frac{3z}{2\ve^2}\Big) h_{1-\ve}(z).$$
\end{lemma}

\proof
Let $v>0$ and $a>v$. An application of the special Markov property gives, for any nonnegative measurable function $\varphi$ on $[0,\infty)$ such that $\varphi(0)=0$,
$$\N_0\Big(\mathbf{1}_{\{W_*>-a\}}\varphi(\z_{-v})\Big)=
\N_0\Big(\varphi(\z_{-v})\,\exp\Big(-\frac{3\z_{-v}}{2(a-v)^2}\Big)\Big)
=\int_0^\infty \dd z\,h_v(z)\varphi(z)\,\exp\Big(-\frac{3z}{2(a-v)^2}\Big).$$
Hence the joint density of the pair $(-W_*,\z_{-v})$ under 
$\N_0(\cdot\cap\{W_*<-v\})$ is the function
$$(a,z)\mapsto \mathbf{1}_{\{a>v\}}\,\frac{3z}{(a-v)^3}\exp\Big(-\frac{3z}{2(a-v)^2}\Big)\,h_v(z).$$
On the other hand, the density of $-W_*$ under $\N_0(\cdot\cap\{W_*<-v\})$ is the function
$a\mapsto \mathbf{1}_{\{a>v\}}\,3a^{-3}$. It follows that, for $a>v$, the density of $\z_{-v}$ under $\N_0(\cdot\midd -W_*=a)$ 
is 
$$z\mapsto \frac{a^3}{(a-v)^3}\,z\,\exp\Big(-\frac{3z}{2(a-v)^2}\Big) h_v(z).$$
The case $a=1$, $v=1-\ve$ gives the first assertion of the lemma.

The proof of the second assertion is straightforward. For a function $\varphi$ as above,
\begin{align*}
\N_0\Big(\varphi(\z_{-1+\ve})\mathbf{1}_{\{-1<W_*<-1+\ve\}}\Big)
&=\N_0\Big(\varphi(\z_{-1+\ve})\,\exp(-\z_{-1+\ve}\N_{-1+\ve}(W_*\leq -1))\Big)\\
&=\N_0\Big(\varphi(\z_{-1+\ve})\,\exp(-\frac{3}{2\ve^2}\z_{-1+\ve})\Big)\\
&=\int_0^\infty \dd z\,h_{1-\ve}(z)\,\varphi(z)\,\exp(-\frac{3z}{2\ve^2})
\end{align*}
and the desired result follows since $\N_0(-1<W_*<-1+\ve)=\frac{3}{2(1-\ve)^2}-\frac{3}{2}$. \endproof

\begin{lemma}
\label{App-Lem3}
Let $\ve\in(0,1)$.
The law of $\z_{W_*+\ve}$ under $\N^{[1]}_0$ has density
$$g_\ve(z):=2\,\ve^{-3}\,z\exp\Big(-\frac{3z}{2\ve^2}\Big) \int_{1-\ve}^\infty \dd a\,h_a(z).$$
\end{lemma}

\proof We rely on a formula found in \cite[Proposition 12]{Disks}, which 
gives for any nonnegative measurable function $\varphi$ on $[0,\infty)$,
\begin{align*}
\N_0(\mathbf{1}_{\{W_*<-1\}}\,\varphi(\z_{W_*+\ve}))
&=3\ve^{-3}\int_{-\infty}^{-1+\ve} \dd b\,\N_0\Big(\z_b\,\exp\Big(-\frac{3\z_b}{2\ve^2}\Big)\,\varphi(\z_b)\Big)\\
&=3\ve^{-3}\int_{1-\ve}^\infty \dd a\int_0^\infty \dd z\,z\varphi(z)\,\exp\Big(-\frac{3z}{2\ve^2}\Big)\,h_a(z)\\
&=3\ve^{-3}\int_0^\infty \dd z\,z\varphi(z)\Big(\int_{1-\ve}^\infty \dd a\,h_a(z)\Big)\exp\Big(-\frac{3z}{2\ve^2}\Big).
\end{align*}
The deired result follows since $\N_0(W_*<-1)=3/2$. \endproof

Let us use the preceding lemmas to evaluate $\N^{[1]}_0(M_\ve=1\midd \z_{W_*+\ve}=z)$. For
any nonnegative measurable function $\varphi$ on $[0,\infty)$, we have
$$\N^{[1]}_0(\mathbf{1}_{\{M_\ve=1\}}\,\varphi(\z_{W_*+\ve}))
=\N^{[1]}_0(M_\ve=1)\times \N^{[1]}_0(\varphi(\z_{W_*+\ve})\midd M_\ve=1)= (1-\ve)\,\int_0^\infty \dd z\,f_\ve(z)\varphi(z),$$
by Lemma \ref{App-Lem1} and Lemma \ref{App-Lem2}, using also the fact that the law of $\z_{W_*+\ve}$
under $\N^{[1]}_0(\cdot\midd M_\ve=1)$ coincides with the law of $\z_{-1+\ve}$ under $\N_0(\cdot\midd W_*=-1)$ 
(see \cite[Proposition 12]{Disks}). On the other hand,
\begin{align*}
\N^{[1]}_0(\mathbf{1}_{\{M_\ve=1\}}\,\varphi(\z_{W_*+\ve}))&=\N^{[1]}_0\Big(\varphi(\z_{W_*+\ve})\,\N^{[1]}_0(M_\ve=1\midd\z_{W_*+\ve})\Big)\\
&=\int_0^\infty \dd z\,g_{\ve}(z)\,\varphi(z)\,\N^{[1]}_0(M_\ve=1\midd \z_{W_*+\ve}=z),
\end{align*}
by Lemma \ref{App-Lem3}.
By comparing the last two displays, we get
$$\N^{[1]}_0(M_\ve=1\midd \z_{W_*+\ve}=z)=(1-\ve)\,\frac{f_\ve(z)}{g_\ve(z)}=\frac{(1-\ve)\,h_{1-\ve}(z)}{2\int_{1-\ve}^\infty \dd a\,h_a(z)}
=\frac{h_1((1-\ve)^{-2}z)}{2 \int_1^\infty \dd a\,h_a((1-\ve)^{-2}z)},$$
where the last equality is a consequence of \eqref{def-h}. 

We can similarly compute $\N^{[1]}_0(M_\ve=2\midd \z_{W_*+\ve}=z)$. Observing that the law of $\z_{W_*+\ve}$
under $\N^{[1]}_0(\cdot\midd M_\ve=2)$ has density $f_\ve*\tilde f_\ve$ (use again \cite[Proposition 12]{Disks}), and recalling 
the formula for $\N^{[1]}_0(M_\ve=2)$ in Lemma \ref{App-Lem2}, the same argument shows that 
\begin{align*}
\N^{[1]}_0(M_\ve=2\midd \z_{W_*+\ve}=z)=\N^{[1]}_0(M_\ve=2)\,\frac{f_\ve*\tilde f_\ve(z)}{g_\ve(z)}
&=\frac{(1-\ve)^3\int_0^z \dd y\, yh_{1-\ve}(y)h_{1-\ve}(z-y)}{6z\int_{1-\ve}^\infty \dd a\,h_a(z)}\\
&=\frac{\int_0^{(1-\ve)^{-2}z} \dd y\, yh_{1}(y)h_{1}((1-\ve)^{-2}z-y)}{6(1-\ve)^{-2}z\int_{1}^\infty \dd a\,h_a((1-\ve)^{-2}z)}.
\end{align*}
Similarly, since the law of $\z_{W_*+\ve}$
under $\N^{[1]}_0(\cdot\midd M_\ve=3)$ has density $f_\ve*\tilde f_\ve*\tilde f_\ve$, we get
\begin{align*}
\N^{[1]}_0(M_\ve=3\midd \z_{W_*+\ve}=z)&=\N^{[1]}_0(M_\ve=3)\,\frac{f_\ve*\tilde f_\ve*\tilde f_\ve(z)}{g_\ve(z)}\\
&=\frac{\int_0^{(1-\ve)^{-2}z} \dd y\, y\,h_{1}(y)\,h_1*h_1((1-\ve)^{-2}z-y)}{12 (1-\ve)^{-2}z\int_{1}^\infty \dd a\,h_a((1-\ve)^{-2}z)}.
\end{align*}
We note that $\N^{[1]}_0(M_\ve=m\midd \z_{W_*+\ve}=z)$ (for $m=1,2,3$) only depends on the
quantity $(1-\ve)^{-2}z$, which we could have seen from a scaling argument.

At this stage, we use the explicit formula for the functions $h_a$ in \eqref{def-h} and \eqref{def-psi} to
get asymptotic expansions as $z\to 0$. We have first
$$h_1(z)=\frac{3^{3/2}2^{-1/2}}{\sqrt{\pi}}\,\frac{1}{\sqrt{z}} - \frac{27}{4} + \frac{3^{5/2}2^{1/2}}{\sqrt{\pi}}\,\sqrt{z}+ O(z)$$
and
\begin{equation}
\label{expan-inte}
2 \int_1^\infty \dd a\,h_a(z)=\frac{3^{3/2}2^{-1/2}}{\sqrt{\pi}}\,\frac{1}{\sqrt{z}} -\frac{9}{2} +\frac{3^{5/2}2^{-1/2}}{\sqrt{\pi}}\,\sqrt{z}+O(z).
\end{equation}
From the preceding formula for $\N^{[1]}_0(M_\ve=1\midd \z_{W_*+\ve}=z)$, it follows that
$$\N^{[1]}_0(M_\ve=1\midd \z_{W_*+\ve}=z)=1+O(\sqrt{z}),\quad \hbox{as }z\to 0,$$
where the remainder $O(\sqrt{z})$ is uniform in $\ve\in (0,1/2]$. Hence
$\N^{[1]}_0(M_\ve\geq2\midd \z_{W_*+\ve}=z)=O(\sqrt{z})$, which gives the case $m=1$ of Lemma \ref{geodesic-exit}.

Similarly, tedious but straightforward calculations show that
\begin{align*}
\int_0^{z} \dd y\, y\,h_{1}(y)\,h_{1}(z-y) &= \frac{27}{4} z - \frac{3^{9/2}2^{-3/2}}{\sqrt{\pi}}\,z^{3/2} + O(z^2)\\
\int_0^{z} \dd y\, y\,h_{1}(y)\,h_1*h_1(z-y)&=\frac{3^{7/2}2^{-1/2}}{\sqrt{\pi}}\,z^{3/2} + O(z^2).
\end{align*}
To simplify notation, write $z'=(1-\ve)^{-2}z$. Then we get
\begin{align*}
\Big(2 \int_1^\infty \dd a\,h_a(z')\Big)\,\N^{[1]}_0(M_\ve\leq 2\midd \z_{W_*+\ve}=z)
&= h_1(z')+ \frac{1}{3z'}\int_0^{z'} \dd y\, yh_{1}(y)h_{1}(z'-y)\\
&=\frac{3^{3/2}2^{-1/2}}{\sqrt{\pi}}\,\,\frac{1}{\sqrt{z'}} -\frac{9}{2} +\frac{3^{5/2}2^{-3/2}}{\sqrt{\pi}}\,z'^{1/2} + O(z').
\end{align*}
Comparing with \eqref{expan-inte}, we obtain that $\N^{[1]}_0(M_\ve\leq 2\midd \z_{W_*+\ve}=z)=1+O(z)$ as $z\to 0$, and thus 
$\N^{[1]}_0(M_\ve\geq 3\midd \z_{W_*+\ve}=z)=O(z)$, giving the case $m=2$ of Lemma \ref{geodesic-exit}.
Finally, we have also
\begin{align*}
&\Big(2 \int_1^\infty \dd a\,h_a(z')\Big)\,\N^{[1]}_0(M_\ve\leq 3\midd \z_{W_*+\ve}=z)\\
&\qquad= h_1(z')+ \frac{1}{3z'}\int_0^{z'} \dd y\, yh_{1}(y)h_{1}(z'-y)+ \frac{1}{6z'}\int_0^{z'} \dd y\, yh_{1}(y)\,h_1*h_1(z'-y) \\
&\qquad=\frac{3^{3/2}2^{-1/2}}{\sqrt{\pi}}\,\frac{1}{\sqrt{z'}} -\frac{9}{2} +\frac{3^{5/2}2^{-1/2}}{\sqrt{\pi}}\,z'^{1/2} + O(z').
\end{align*}
Comparing again with \eqref{expan-inte}, we get that $\N^{[1]}_0(M_\ve\leq 3\midd \z_{W_*+\ve}=z)=1+O(z^{3/2})$ as $z\to 0$, and thus 
$\N^{[1]}_0(M_\ve\geq 4\midd \z_{W_*+\ve}=z)=O(z^{3/2})$, which gives the case $m=3$ of Lemma \ref{geodesic-exit}
and completes the proof of this lemma. 

\subsection*{Proof of Lemma \ref{coupling-Bplane}}

We start by recalling the construction of the Brownian plane as described in \cite{CLG}.
We consider a nine-dimensional Bessel process $R=(R_t)_{t\geq 0}$ started at $0$, and,
conditionally on $R$, two independent Poisson point measures $\n$ and $\n'$ on $[0,\infty)\times \S$
with the same intensity
$$\dd t\,\N_{R_{t}}(\dd \omega\cap\{W_*(\omega)>0\}).$$
It is convenient to write
$$\n=\sum_{i\in I} \delta_{(t_i,\omega_i)}\,,\ \n'=\sum_{i\in J}  \delta_{(t_i,\omega_i)},$$
where the indexing sets $I$ and $J$ are disjoint. We then consider the (non-compact) $\R$-tree
$\t_\infty$ defined by
$$\t_\infty:=[0,\infty)\cup \Big( \bigcup_{i\in I\cup J} \t_{(\omega_i)}\Big),$$
where for every $i\in I\cup J$, the root of $\t_{(\omega_i)}$ is identified with the point 
$t_i$ of $[0,\infty)$: we view $\t_{(\omega_i)}$ as grafted on the ``spine'' $[0,\infty)$
at height $t_i$. In fact, we consider the trees  $\t_{(\omega_i)}$ for $i\in I$ as grafted
to the left side of the spine, and the trees  $\t_{(\omega_i)}$ for $i\in J$ as grafted 
to the right side of the spine. This is reflected in the exploration process $(\ee^\infty_s)_{s\in\R}$
of $\t_\infty$, which is such that $\{\ee^\infty_s:s\leq 0\}$ is exactly the union of 
the spine and of the trees $\t_{(\omega_i)}$ for $i\in I$, whereas 
$\{\ee^\infty_s:s\geq 0\}$ is the union of 
the spine and of the trees $\t_{(\omega_i)}$ for $i\in J$ (we refer to \cite[Section 2.4]{Spine}
for a more precise definition of $(\ee^\infty_s)_{s\in\R}$). The exploration process allows us 
to define intervals on the tree $\t_\infty$. We make the convention that, if 
$s,s'\in \R$ and $s>s'$, the ``interval'' $[s,s']$ is equal to $[s,\infty)\cup(-\infty,s']$. Then 
if $a,b\in\t_\infty$, there is a smallest ``interval'' $[s,s']$ such that $\ee^\infty_s=a$
and $\ee^\infty_{s'}=b$, and we take $[a,b]_\infty=\{\ee^\infty_r:r\in[s,s']\}$. 

We also assign labels $(\Lambda_a)_{a\in\t_\infty}$ to the points of $\t_\infty$. If 
$a=t$ belongs to the spine $[0,\infty)$, we take $\Lambda_a=R_t$. If $a\in\t_{(\omega_i)}$
for some $i\in I\cup J$, we let $\Lambda_a$ be the label of $a$ in $\t_{(\omega_i)}$. We then
define, for every $a,b\in \t_\infty$,
$$D^{\infty,\circ}(a,b):=\Lambda_a+\Lambda_b-2\max\Big(\min_{c\in[a,b]_\infty}\Lambda_c,\min_{c\in[b,a]_\infty}\Lambda_c\Big),$$
and we let $D^\infty(a,b)$ be the maximal symmetric function of the pair $(a,b)$ that is
bounded above by $D^{\infty,\circ}(a,b)$ and satisfies the triangle inequality. It turns out that
the property $D^\infty(a,b)=0$ holds if and only if $D^{\infty,\circ}(a,b)=0$. The Brownian 
plane $\pp$ can then be defined as the quotient space $\t_\infty/\{D^\infty=0\}$, which is equipped
with the distance induced by $D^\infty$ and with a distinguished point which is the 
point $0$ of the spine. We write $\Pi_\infty$ for the canonical projection from 
$\t_\infty$ onto $\pp$.

Thanks to Section \ref{spine-decomp}, the Bessel process $R$ and the point measures $\n$ and $\n'$ can also be used to construct a
random snake trajectory distributed according to $\N_0(\cdot\midd W_*=-1)$, whose genealogical tree
is identified to
$$\t_1:=[0,L_1]\cup \Big( \bigcup_{i\in I\cup J,t_i\leq L_1} \t_{(\omega_i)}\Big),$$
where we make the same identifications as for $\t_\infty$ and, for every $r>0$, we have set
$$L_r:=\sup\{t\geq 0:R_t=r\}.$$
We may and will view
$\t_1$ as the subset of $\t_\infty$ obtained by
removing the part of the spine above height $L_1$ (and of course the trees $\t_{(\omega_i)}$
grafted to this part). Write $[a,b]_1$ for the intervals on the tree $\t_1$. Following the construction
of the Brownian sphere in Section \ref{Brown-sphere}, we define $D^{1,\circ}(a,b)$ 
for $a,b\in\t_1$ by the very
same formula as $D^{\infty,\circ}(a,b)$ above, but replacing the intervals 
$[a,b]_\infty$ and $[b,a]_\infty$ by $[a,b]_1$ and $[b,a]_1$ respectively. 
We note that $D^{1,\circ}(a,b)\leq D^{\infty,\circ}(a,b)$
for $a,b\in\t_1$ since we have clearly $[a,b]_1\subset [a,b]_\infty$.
Finally, we let $D^1(a,b)$ be the maximal symmetric function of $a,b\in\t_1$ that is
bounded above by $D^{1,\circ}(a,b)$ and satisfies the triangle inequality, and
we define $\bm_\infty^{\{1\}}$ as the quotient space
$\t_1/\{D^1=0\}$, which is equipped with the metric induced by $D^1$
and the two distinguished points which are the points $0$ and $L_1$ (bottom
and top of the spine).
Then $\bm_\infty^{\{1\}}$ has the distribution of $\bm_\infty$
under $\N_0(\cdot\midd W_*=-1)$. We write $\Pi_1$ for the canonical projection from 
$\t_1$ onto $\bm_\infty^{\{1\}}$.

We then claim that the conclusion of Lemma \ref{coupling-Bplane} holds 
if we take
$$E_\ve:=\{W_*(\omega_i)>1-\frac{\ve}{2}\hbox{ for every }i\in I\cup J\hbox{ such that }t_i>1\}.$$
We note that $\P(E_\ve)>0$ (as a simple consequence of the formula for
$\N_0(W_*<-r)$) and that $E_\ve$ is independent of $\bm_\infty^{\{1\}}$. We then
verify that, if $E_\ve$ holds, $B^\bullet_{1-\ve}(\bm^{\{1\}}_\infty)$ and $B^\bullet_{1-\ve}(\mathcal{P})$
can be identified as sets. For every $r>0$, let $F^\infty_r$
be the set of all $a\in\t_\infty$ such that the minimal label along the
geodesic from $a$ to $\infty$ in $\t_\infty$ is smaller than
or equal to $r$. Then $B^\bullet_r(\pp)=\Pi_\infty(F^\infty_r)$
(see formula (16) in \cite{CLG}). Similarly, for every $r\in(0,1)$, let $F^1_r$
be the set of all $a\in\t_1$ such that the minimal label along the
geodesic from $a$ to $L_1$ in $\t_1$ is smaller than
or equal to $r$. Then, we have
$B^\bullet_{r}(\bm^{\{1\}}_\infty)=\Pi_1(F^1_r)$ as a consequence of the bound \eqref{cactus-bd}. Next, on the event $E_\ve$,
one immediately gets that $F^\infty_{1-\ve}=F^1_{1-\ve}$. Furthermore,
still on the event $E_\ve$,
for $a,b\in F^\infty_{1-\ve}$, we have $\Pi_\infty(a)=\Pi_\infty(b)$
if and only $\Pi_1(a)=\Pi_1(b)$. The fact that $\Pi_\infty(a)=\Pi_\infty(b)$
implies $\Pi_1(a)=\Pi_1(b)$ is trivial since $D^1\leq D^\infty$. Conversely, if $\Pi_1(a)=\Pi_1(b)$, the only case 
where we do not immediately get $\Pi_\infty(a)=\Pi_\infty(b)$ is when
$[a,b]_1\not =[a,b]_\infty$ and $\Lambda_a=\Lambda_b=\min\{\Lambda_c:c\in [a,b]_1\}$.
In that case however, the interval $[a,b]_1$ must contain the point $L_1$
(top of the spine), and also the geodesic from $a$ to $L_1$ in $\t_1$, so that
 $\Lambda_a=\Lambda_b\leq 1-\ve$ (by the definition of $F^1_{1-\ve}$), and then 
$\min\{\Lambda_c:c\in [a,b]_1\}=\min\{\Lambda_c:c\in [a,b]_\infty\}$
(because labels on $\t_\infty\backslash \t_1$ are greater than $1-\ve/2$).
Finally we have also $\Pi_\infty(a)=\Pi_\infty(b)$ in that case. 

The preceding considerations show that, on the event $E_\ve$, 
the sets $B^\bullet_{1-\ve}(\bm^{\{1\}}_\infty)$ and $B^\bullet_{1-\ve}(\mathcal{P})$ are
identified. Very similar arguments (using formula (17) in \cite{CLG}, and its analog
for $\bm^{\{1\}}_\infty$) show that the boundary of $B^\bullet_{1-\ve}(\bm^{\{1\}}_\infty)$ in
$\bm^{\{1\}}_\infty$ is also identified to the boundary of $B^\bullet_{1-\ve}(\mathcal{P})$
in $\pp$. Note that the topology induced by $D^\infty$ on $B^\bullet_{1-\ve}(\bm^{\{1\}}_\infty)=
B^\bullet_{1-\ve}(\mathcal{P})$ must be the same as the one induced by $D^1$
since both are compact and $D^1\leq D^\infty$. 
Finally, using the definitions of the distances $D^\infty$ and $D^1$,
one checks that, for every compact subset $K$ of 
$\mathrm{Int}(B^\bullet_{1-\ve}(\bm^{\{1\}}_\infty))=\mathrm{Int}(B^\bullet_{1-\ve}(\mathcal{P}))$,
for every $x,y\in K$ such that $D^\infty(x,y)$ is small enough,
resp. such that $D^1(x,y)$ is small enough, one has $D^\infty(x,y)=D^1(x,y)$
(we omit a few details here). It follows that the intrinsic distance on 
$\mathrm{Int}(B^\bullet_{1-\ve}(\bm^{\{1\}}_\infty))$ coincides
with the intrinsic distance on $\mathrm{Int}(B^\bullet_{1-\ve}(\mathcal{P}))$, and
this completes the proof. 

\subsection*{Proof of Lemma \ref{coupling-exit}}

We assume that the Brownian plane $\pp$ is constructed as explained in the
preceding proof. Let $r>0$. According to formula (18) in \cite{CLG}, the random 
variable $Z_r$ can be obtained as 
$$Z_r=\sum_{i\in I\cup J, t_i>L_r} \z_r(\omega_i).$$

If $0<r<u$, the spine decomposition of $\N_0(\cdot\midd W_*=-u)$ in Section \ref{spine-decomp} shows that
\begin{equation}
\label{exit-formula-spine}
\sum_{i\in I\cup J,L_r<t_i<L_u} \z_r(\omega_i)
\end{equation}
has the distribution of $\z_{W_*+r}$ under $\N_0(\cdot\midd W_*=-u)$
(compare \eqref{exit-formula-spine} with the right-hand side of \eqref{def-exit2}). 
From the last two displays, we immediately obtain that
the random variable $\z_{W_*+r}$ under $\N_0(\cdot\midd W_*=-u)$
is stochastically dominated by $Z_r$. The lemma follows.

\end{document}